\input amstex
\font\smallbf=ptmbo at 10pt

\font\tensans=cmss10

\define\lx{\text{\tensans[}}
\define\rx{\text{\tensans]}}

\define\syd{\sy.2.d}
\define\sye{\sy.2.e}

\define\bbP{{\bold P}}
\define\bbR{{\bold R}}

\define\bbC{{\bold C}}

\define\bbCP{{\bold C}\text{\rm P}}

\define\ren{\text{\bf R}\hskip-.7pt^n}
\define\renr{\text{\bf R}\hskip-.7pt^{n-r}}

\define\dimr{\dim_{\hskip.1pt\bbR\hskip-1.7pt}}
\define\dimc{\dim_{\hskip.4pt\bbC\hskip-1.2pt}}

\define\Lie{\text{\smallbf L}\hskip.5pt}

\define\sgn{\text{\rm sgn}\hs}

\define\y{c}

\define\bz{b\hs}
\define\dg{r}
\define\fy{\sigma}

\define\ax{A}
\define\bx{B}
\define\cx{C}
\define\fe{f}
\define\iy{\Cal I}
\define\iyp{\Cal I\hs'}
\define\jy{\Cal J}

\define\df{d\hskip-.8ptf}

\define\e{\text{\rm exp}}
\define\uexp{U^{\text{\rm Exp}}}
\define\uo{U^{\text{\rm o}}}

\define\hs{\hskip.7pt}
\define\nh{\hskip-1pt}
\define\hn{\Cal H^{N\!}}
\define\hnyz{\Cal H_*}
\define\vn{v^{N\!}}
\define\un{u^{N\!}}
\define\ptmi{\phantom{i}}

\define\qd{\hskip7pt$\blacksquare$}


\define\dsq{d^{\hskip.6pt2}\hskip-.7pt}

\define\hrz{^{\hskip.5pt\text{\rm hrz}}}
\define\vrt{^{\hskip.2pt\text{\rm vrt}}}
\define\nrm{^{\hskip.2pt\text{\rm nrm}}}
\define\a{}
\define\f{\thetag}
\define\ff{\tag}

\define\vp{{\tau\hskip-4.55pt\iota\hskip.6pt}}
\define\vpsq{{\tau\hskip-4.55pt\iota\hskip1.3pt^2}}
\define\vpsu{{\tau\hskip-3.8pt\iota\hskip1pt}}
\define\vpab{{\tau\hskip-3.7pt\iota\hskip.2pt}}
\define\navpab{{\nabla\hskip-.2pt\tau\hskip-3.9pt\iota\hskip.2pt}}
\define\mgmt{M,g,m\hs,\hskip.2pt\vp}

\define\si{\phi}
\define\ta{\psi}
\define\la{\lambda}
\define\my{\mu}

\define\sa{s}
\define\xe{\zeta}
\define\ps{t}

\define\zx{z}
\define\zy{\zeta}
\define\srf{S}
\define\hsf{\hat S}
\define\gx{\gamma}

\define\diml{-di\-men\-sion\-al}

\define\sky{skew-sym\-me\-try}

\define\kip{Killing potential}
\define\om{\omega}

\define\omfs{\omega_{\hs\text{\sixrm FS}}}
\define\ve{\varepsilon}

\define\kx{\kappa}

\define\ri{\text{\rm r}}
\define\rih{\ri^{(h)}}

\define\navp{\nabla\hskip-.2pt\vp}
\define\ts{\eta}
\define\proj{\pi} 

\define\krp{K\"ahler-Ricci potential}
\define\rmk{\example}
\define\endrmk{\endexample}
\define\kp{5}
\define\ck{6}
\define\sr{7}
\define\mc{9}
\define\qs{11}
\define\cm{12}
\define\rr{13}
\define\dx{16}
\define\lc{17}

\define\ty{22}
\define\ec{23}
\define\ls{24}
\define\tg{25}


\define\ir{0}
\define\pn{1}
\define\cn{2}
\define\tb{3}
\define\kg{4}
\define\cz{5}
\define\gd{6}
\define\mb{7}
\define\sm{8}
\define\xa{9}
\define\am{10}
\define\tm{11}
\define\mm{12}
\define\ib{13}
\define\os{14}
\define\md{15}
\define\ob{16}
\define\mw{17}
\define\cp{18}
\define\dc{19}
\define\cg{20}
\define\pf{21}
\define\bo{22}
\define\db{23}
\define\nd{24}
\define\cu{25}
\define\cv{26}
\define\fv{27}
\define\nx{28}
\define\gc{29}
\define\cf{30}
\define\rs{31}
\define\bc{32}
\define\tf{33}
\define\sy{34}
\define\sd{35}
\define\tl{36}
\define\at{37}


\define\fp{13}


\define\agc{1}
\define\acg{2}
\define\ber{3}
\define\btt{4}
\define\bry{5}
\define\cal{6}
\define\cao{7}
\define\dml{8}
\define\dmg{9}
\define\dms{10}
\define\kno{11}
\define\koo{12}
\define\mor{13}
\define\ptv{14}
\define\tiz{15}
\documentstyle{amsppt}
\magnification=1200
\NoBlackBoxes
\topmatter
\title Special K\"ahler-Ricci potentials\\
on compact K\"ahler manifolds
\endtitle
\rightheadtext{Special K\"ahler-Ricci potentials}
\author A. Derdzinski and G. Maschler\endauthor
\address Dept. of Mathematics, The Ohio State University,
Columbus, OH 43210, USA\endaddress
\email andrzej\@math.ohio-state.edu\endemail
\address Department of Mathematics, University of Toronto, Canada M5S 3G3
\endaddress
\email maschler\@math.toronto.edu\endemail
\abstract{A special K\"ahler-Ricci potential on a K\"ahler manifold is any
nonconstant $\,C^\infty$ function $\,\vpab\,$ such that $\,J(\navpab)\,$ is a
Killing vector field and, at every point with $\,\,d\vpab\ne0$, all nonzero
tangent vectors orthogonal to $\,\navpab\,$ and $\,J(\navpab)\,$ are
eigenvectors of both $\,\nabla d\vpab\,$ and the Ricci tensor. For instance,
this is always the case if $\,\vpab\,$ is a nonconstant $\,C^\infty$ function
on a K\"ahler manifold $\,(M,g)\,$ of complex dimension $\,m>2\,$ and the
metric $\,\tilde g=g/\vpab^2$, defined wherever $\,\vpab\ne0$, is Einstein.
(When such $\,\vpab\,$ exists, $\,(M,g)\,$ may be called {\it
almost-everywhere conformally Einstein}.) We provide a complete classification
of compact K\"ahler manifolds with special K\"ahler-Ricci potentials and use
it to prove a structure theorem for compact K\"ahler manifolds of any complex
dimension $\,m>2\,$ which are almost-everywhere conformally Einstein.}
\endabstract
\endtopmatter
\voffset=-35pt
\document
\head\S\ir. Introduction\endhead
This paper, although self-contained, can also be viewed as the second in a
series of three papers that starts with \cite{\dml} and ends with \cite{\dmg}.

We call $\,\vp\,$ a {\it special K\"ahler-Ricci
potential\/} on a K\"ahler manifold $\,(M,g)\,$ if
\vskip4pt
\settabs\+\noindent&\f{\ir.1}\hskip11pt&\cr
\+&&$\vp\hs\,$ is a nonconstant Killing potential on $\,\hs(M,g)\hs\,$
and,\hskip6ptat every point\cr
\+&\f{\ir.1}&with $\,d\vp\ne0$,\hskip3.8ptall nonzero tangent vectors
orthogonal to $\hs v=\navp\hs$ and\cr
\+&&to $\,\,u\hs=\hs Jv\,\,$ are\hs\ eigenvectors\hs\ of\hs\ both
$\,\,\nabla d\vp\,\,$ and\hs\ the\hs\ Ricci\hs\ tensor $\hs\,\,\ri\hs$.\cr
\vskip4pt
\noindent(Cf. \cite{\dml}, \S\sr; for more on Killing potentials, see \S\kg\
below.) The word `potential' reflects the fact that \f{\ir.1} is closely
related, although not equivalent, to the requirement that
$\,\nabla d\vp+\chi\,\ri\,=\fy g\,$ for some $\,C^\infty$ functions
$\,\chi,\fy\,$ (see \cite{\dml}, beginning of \S\sr). This requirement is
reminiscent of K\"ahler-Ricci solitons (some of which, in fact, do satisfy
\f{\ir.1}, cf. \cite{\dms} and Remark \a\am.1 below); while, in complex
dimensions $\,m>2$, it implies that $\,\vp\,$ arises from a {\it Hamiltonian
$\,2$-form\/} on the underlying K\"ahler manifold (\cite{\acg}, \S1.4). See
also \cite{\bry}. What further sparked our interest in \f{\ir.1} was its
being, in cases such as \f{\ir.4} below, a consequence of the following
assumption:
\vskip4pt
\settabs\+\noindent&\f{\ir.2}\hskip24pt&\cr
\+&&$(M,g)\hs$ is a K\"ahler manifold of complex dimension $\hs m\hs$
and $\hs\vp\hs$ is\cr
\+&\f{\ir.2}&a nonconstant $\hs\,C^\infty\hs$ function on $\hs\,M\,\hs$ such
that the conformally\cr
\+&&related metric $\,\,\tilde g=g/\vpsq$,\hskip4ptdefined wherever
$\,\,\vp\ne0$,\hskip4ptis Einstein.\cr
\vskip4pt
\noindent
When $\,m>2$, \f{\ir.2} implies the seemingly stronger condition
$$\text{\rm$\mgmt\,$ satisfy \f{\ir.2} and
$\,d\vp\wedge\hs d\Delta\vp=0\,$ everywhere in $\,M$}\ff\ir.3$$
(see \cite{\dml}, Proposition \a\ck.4), so that locally, at points with
$\,d\vp\ne0$, the Laplacian of $\,\vp\,$ is a function of $\,\vp$. Therefore,
\f{\ir.3} is of independent interest, as opposed to just \f{\ir.2}, only for
K\"ahler {\it surfaces\/} ($m=2$). In \cite{\dml}, Corollary \a\mc.3, we found
that
$$\text{\rm Condition \f{\ir.2} with $\,m\ge3$, or \f{\ir.3} with $\,m=2$,
implies \f{\ir.1}.}\ff\ir.4$$

The first main result of this paper is a complete classification of compact
K\"ahler manifolds $\,(M,g)$, in all complex dimensions $\,m\ge1$, with
functions $\,\vp\,$ satisfying \f{\ir.1}. Aside from the trivial case
$\,m=1\,$ (see Example \a\xa.1), we show in Theorem \a\gc.2 that, for each
fixed $\,m\ge2$, such manifolds form two separate families: in one, described
in \S\mw, $\,M\,$ is a holomorphic $\,\bbCP^1$ bundle (of a very specific
kind) over a compact K\"ahler manifold which is Einstein unless $\,m=2$, while
the other has $\,M\,$ biholomorphic to $\,\bbCP^m$ (see \S\cp). Pairs
$\,g,\vp\,$ with \f{\ir.1} on a fixed compact complex manifold $\,M\,$ turn
out to form an infinite\diml\ moduli space. Namely, $\,Q=g(\navp,\navp)\,$
then is a $\,C^\infty$ function of the real variable
$\,\vp\in[\hs\text{\rm min}\,\vp\hs,\,\text{\rm max}\,\vp]$, satisfying the
positivity and boundary conditions \f{\mw.1} in \S\mw; conversely, any given
assignment $\,\vp\mapsto Q\,$ with \f{\mw.1} is realized in this way on each
of the complex manifolds $\,M\,$ just mentioned. See \S\mw, \S\cp.

Theorem \a\gc.2 combined with \f{\ir.4} leads to our second main result,
consisting of four statements that together form a structure theorem for, and
a partial classification of, those quadruples $\,\mgmt\,$ with compact $\,M\,$
which satisfy \f{\ir.2} with $\,m\ge3$, or \f{\ir.3} with $\,m=2$.
Specifically, in \S\tf\ we divide all such quadruples into four
disjoint ``types'' (a), (b), (c1), (c2). We then prove that type (b) is empty
(Theorem \a\tf.2), and verify (in Corollary \a\sd.1) that type (c2) leads to a
certain necessary condition which, as we show in a subsequent paper
\cite{\dmg}, is never satisfied; therefore, type (c2) eventually turns out to
be empty as well. Next, in Theorem \a\tf.3, we completely classify type (a);
every $\,M\,$ occurring in it is a flat holomorphic $\,\bbCP^1$ bundle.
Finally, our Theorem \a\sy.3 reduces the classification of type (c1) (in which
$\,M\,$ always is a {\it nonflat\/} holomorphic $\,\bbCP^1$ bundle) to the
question of finding all rational functions that lie in a specific three\diml\
vector space depending on $\,m\,$ and satisfy an analogue of \f{\mw.1}. An
answer to this last question, with examples (not limited to those of
\cite{\ber}) is, again, postponed to \cite{\dmg}, as it requires extensive
additional arguments based on entirely different methods.

The text is organized as follows.
Sections \pn\ -- \mb\ cover preliminary material. A discussion of basic
properties of special \krp s on K\"ahler manifolds in \S\sm\ is followed by
constructions of examples in sections \xa\ -- \cp.

Critical manifolds of special \krp s, their
dimensions, geodesics normal to them, as well as their normal connections,
curvature properties, and normal exponential mappings are studied in sections
\dc, \db\ -- \nx\ and \at.

Further properties of critical manifolds established in \S\pf\ allow us to
show, in sections \cg\ and \bo, that the assignment $\,\vp\mapsto Q\,$
mentioned above satisfies conditions \f{\mw.1}, which we then use in \S\gc\ to
prove Theorem \a\gc.2. A similar local result is obtained in \S\tl. The
remaining sections \cf\ -- \sd\ deal with quadruples $\,\mgmt\,$ having the
property \f{\ir.2} with $\,m\ge3$, or \f{\ir.3} with $\,m=2$.

\head\S\pn. Preliminaries\endhead
Except in \S\cn, the symbol $\,\nabla\,$ will stand either for the Levi-Civita
connection of a given Riemannian metric $\,g$, or for the $\,g$-gra\-di\-ent.
For a $\,C^1$ vector field $\,v\,$ on a Riemannian manifold $\,(M,g)\,$ we
will write
$$\nabla v:TM\to TM\qquad\text{\rm with}\quad(\nabla v)w\,=\,\nabla_{\!w}v\,,
\ff\pn.1$$
treating the covariant derivative $\,\nabla v\,$ as a vector-bundle morphism
sending each $\,w\in T_xM$, $\,x\in M$, to $\,\nabla_{\!w}v\in T_xM$. For the
second covariant derivative $\nabla d\vp\,$ of a $\,C^2$ function $\,\vp\,$
and vector fields $\,u,w\,$ on a Riemannian manifold we have
$$(\nabla d\vp)(u,w)\,=\,g(u,\nabla_{\!w}v)\,=\,g(\nabla_{\!u}v,w)\,,\qquad
\text{\rm where}\quad v\,=\,\navp\,.\ff\pn.2$$
Thus, if $\,\vp\,$ is a $\,C^2$ function on a Riemannian manifold,
\vskip4pt
\hbox{\hskip-.8pt
\vbox{\hbox{\f{\pn.3}}\vskip3.2pt}
\hskip16pt
\vbox{
\hbox{the\hskip4.2pteigenvalues\hskip4.2ptand\hskip4.2pteigenvectors\hskip4.2ptof\hskip4.2ptthe\hskip4.2ptsymmetric $\,2$-tensor $\,\,\nabla d\vp$}
\vskip1pt
\hbox{are, at any point, the same as those of $\,\nabla v\,$ with
\f{\pn.1}, for $\,v=\navp$.}}}
\vskip4pt
\noindent The tensor product of $\,1$-forms $\,\xi,\xi'$ acts on tangent
vectors $\,u,v\,$ by
$$(\xi\hskip-.7pt\otimes\xi')(u,v)\,=\,\xi(u)\hs\xi'(v)\,.\ff\pn.4$$
We use the symbol $\,J\,$ for the complex-structure tensor of any complex
manifold $\,M$. Thus, $\,J\,$ is a real vector-bundle morphism $\,TM\to TM\,$
with $\,J^2=-1$. At the same time, we treat $\,TM\,$ as a complex vector
bundle with the multiplication by $\,\,i\,\,$ provided by $\,J$. A K\"ahler
manifold $\,(M,g)\,$ is, as usual, a complex manifold $\,M\,$ with a
Riemannian metric $\,g\,$ which makes the complex-structure tensor $\,J\,$
skew-adjoint and parallel. The K\"ahler form $\,\hs\om\,$ of $\,(M,g)\,$ then
is given by
$$\om(u,v)\,=\,g(Ju,v)\,,\hskip12pt\text{\rm for}\hskip6ptu,v\in T_xM,
\hskip5ptx\in M.\ff\pn.5$$
\rmk{Remark \a\pn.1}Let $\,\varphi\,$ be a $\,C^{k+1}$ function,
$\,0\le k\le\infty$, of a real variable $\,\sa$, defined on an interval
containing $\,0\,$ (possibly as an endpoint), and such that $\,\varphi(0)=0$.
Then $\,\varphi(\sa)/\sa\,$ can be extended to a $\,C^k$ function of $\,\sa\,$
defined on the same interval, including $\,\sa=0$. In fact, integrating
$\,d\hs[\varphi(\sa\fy)]/d\fy\,$ we obtain the Taylor formula
$\,\varphi(\sa)=\sa H(\sa)$, where
$\,H(\sa)=\int_0^1\dot \varphi(\sa\fy)\hs d\fy\,$ with
$\,\dot \varphi=\hs d\varphi/d\sa$.
\endrmk
\rmk{Remark \a\pn.2}If $\,\vp,Q\,$ be $\,C^\infty$-differentiable even
functions of a real variable $\,\sa$, defined on a neighborhood of $\,0\,$ in
$\,\bbR\hs$, and $\,\dsq\vp/d\sa^2\ne0\,$ at $\,\sa=0$, then $\,Q\,$
restricted to some neighborhood of $\,0\,$ is a $\,C^\infty$ function of
$\,\vp$, that is, a composite in which $\,\vp\,$ is followed by a $\,C^\infty$
function of the variable $\,\vp$, defined on a suitable interval.

In fact, let $\,\xe=\sa^2$. By induction on $\,k\ge0$, any even $\,C^{2k}$
function of $\,\sa\,$ is a $\,C^k$ function of $\,\xe$. Namely, for $\,k=0\,$
this is true since $\,\sqrt\xe\,$ is a continuous function of
$\,\xe\ge0$. Assuming our claim for a given $\,k\ge0$, let $\,f\,$ be an even
$\,C^{2k+2}$ function of $\,\sa$. The $\,C^{2k+1}$ function
$\,\dot f=\hs\df\nh/d\sa\,$ then is odd, and so $\,\dot f(\sa)/\sa\,$ is an
even $\,C^{2k}$ function of $\,\sa$, also at $\,\sa=0\,$ (Remark \a\pn.1).
Hence, by the inductive assumption, $\,2\,\df\nh/d\xe=\dot f(\sa)/\sa\,$ is a
$\,C^k$ function of $\,\xe$, i.e., $\,f\,$ is of class $\,C^{k+1}$ in $\,\xe$.

Consequently, $\,\vp\,$ and $\,Q$, are $\,C^\infty$ functions of
$\,\xe=\sa^2$. However, the limit of $\,2\hs d\vp/d\xe=\dot\vp/\sa\,$ as
$\,\sa\to0\,$ (or, $\,\xe\to0$) equals $\,\ddot\vp(0)\ne0$, so that the
assignment $\,\xe\mapsto\vp\,$ is a $\,C^\infty$ diffeomorphism for
$\,\xe\ge0\,$ near $\,0$.
\endrmk

\head\S\cn. Connections and curvature\endhead
Our sign convention for the curvature tensor $\,R\,$ of a (linear) connection
$\,\nabla\,$ in any real/com\-plex vector bundle $\,\Cal E\,$ over a manifold
is
$$R(u,v)w\,=\,\nabla_{\!v}\nabla_{\!u}w\,-\,
\nabla_{\!u}\nabla_{\!v}w\,+\,\nabla_{[u,v]}w\,,\ff\cn.1$$
where $\,u,v\,$ are $\,C^2$ vector fields tangent to the base and $\,w\,$ is a
$\,C^2$ section of $\,\Cal E$.
\rmk{Remark \a\cn.1}The {\it curvature form\/} of a connection $\,\nabla\,$
in a complex line bundle over any manifold is the complex-val\-ued $\,2$-form
$\,\varOmega\,$ with $\,R(u,v)w=i\hs\varOmega(u,v)w\,$ for $\,u,v,w,R\,$ as in
\f{\cn.1}. Any local $\,C^\infty$ section $\,w\,$ without zeros gives rise to
the {\it connection form\/} $\,\varGamma\,$ defined by
$\,\nabla_{\!v}w=\varGamma(v)\hs w$, and \f{\cn.1} easily yields
$\,\varOmega=i\hskip1ptd\varGamma$.
\endrmk
Given a vector bundle $\,\Cal L\,$ over a manifold $\,N$, we will use
the same symbol $\,\Cal L\,$ for its total space, so that
$$\Cal L\,=\,\{(y,\zx):y\in N,\hskip6pt\zx\in\Cal L_y\}\qquad\text{\rm and}
\qquad N\,\subset\,\Cal L\,,\ff\cn.2$$
where $\,N\,$ is identified with the zero section formed by all $\,(y,0)\,$
with $\,y\in N$. We similarly treat each fibre $\,\Cal L_y$ as a subset of
$\,\Cal L\hs$, identifying it with $\,\{y\}\times\Cal L_y$. In this way $\,N\hs$
and all $\,\Cal L_y$ are submanifolds of $\,\Cal L\,$ with its obvious
manifold structure. Being a vector space, every fibre $\,\Cal L_y$ is
naturally identified, for any $\,\zx\in\Cal L_y$, with the tangent space
$\,\Cal V_{(y,\zx)}\subset T_{(y,\zx)}\Cal L\,$ of the submanifold
$\,\{y\}\times\Cal L_y$ at $\,(y,\zx)$. These $\,\Cal V_{(y,\zx)}$ form a
vector subbundle $\,\Cal V\,$ of $\,T\Cal L\,$ called the {\it vertical
distribution\/} of $\,\Cal L\hs$. If $\,\Cal L\,$ is a complex line bundle
over $\,N\nh$, any fixed real number $\,a\ne0\,$ gives rise to the vertical
vector fields $\,v,u\,$ on $\,\Cal L\hs$, that is, sections of $\,\Cal V\nh$,
with
$$v(y,\zx)\,=\,a\zx\,,\qquad u(y,\zx)\,=\,i\hs a\zx\qquad\text{\rm for\ all}
\quad(y,\zx)\in\Cal L\,.\ff\cn.3$$
\rmk{Remark \a\cn.2}Any Riem\-ann\-i\-an/Her\-mit\-i\-an fibre metric
$\,\langle\,,\rangle\,$ in a real/com\-plex vector bundle $\,\Cal L\,$ over a
manifold $\,N\hs$ is determined by its {\it norm function\/}
$\,\Cal L\to[\hs0,\infty)$, denoted $\,r\,$ or $\,\sa$, which assigns
$\,|\zx|=\langle\zx,\zx\rangle^{1/2}$
to each $\,(y,\zx)\in\Cal L\hs$. As $\,\Cal V_{(y,\zx)}=\Cal L_y$ (see above),
$\,\langle\,,\rangle\,$ may also be treated as a fibre metric in the vertical
subbundle $\,\Cal V\hs$ of $\,T\Cal L\hs$, and then, if $\,\Cal L\,$ is a
complex line bundle, for $\,v,u\,$ given by \f{\cn.3} we have
$\,\langle v,v\rangle=\langle u,u\rangle=a^2r^2$ and
$\,\,\text{\rm Re}\hskip1pt\langle v,u\rangle=0$.
\endrmk
\rmk{Remark \a\cn.3}Suppose that $\,M,\hat M\,$ are locally trivial fibre
bundles over a manifold $\,N\hs$ and a $\,C^\infty$ diffeomorphism
$\,\varPhi:M\to\hat M\,$ is fibre-preserving, i.e.,
$\,\hat\proj\circ\varPhi=\proj$, where $\,\proj,\hat\proj\,$ denote the bundle
projections $\,M\to N\,$ and $\,\hat M\to N\nh$. Also, let
$\,\Cal H,\hat{\Cal H}\,$ be some fixed ``horizontal'' distributions in
$\,M\,$ and $\,\hat M$, that is, vector subbundles of $\,TM\,$ and
$\,T\hat M\,$ with $\,TM=\Cal H\oplus\Cal V\,$ and
$\,T\hat M=\hat{\Cal H}\oplus\hat{\Cal V}\,$ for the vertical distributions
$\,\Cal V,\hat{\Cal V}\,$ (tangent to the fibres). If $\,\varPhi\,$ sends
$\,\Cal H\,$ onto $\,\hat{\Cal H}$, then the differential $\,d\varPhi_x$ of
$\,\varPhi\,$ at any $\,x\in M$, restricted to $\,\Cal H_x$, preserves any
fibre metric or complex-bundle structure obtained in both $\,\Cal H\,$ and
$\,\hat{\Cal H}\,$ by pulling back a fixed analogous object in $\,TN\,$ via
$\,\proj:M\to N\,$ and $\,\hat\proj:\hat M\to N\nh$.

In fact, $\,d\varPhi_x$ acts as the identity mapping between $\,\Cal H_x$ and
$\,\hat{\Cal H}_{\varPhi(x)}$ identified, via the differential of $\,\proj\,$
or $\,\hat\proj$, with $\,T_yN\,$ at $\,y=\proj(x)$.
\endrmk
\rmk{Remark \a\cn.4}Let $\,\Cal L\,$ be a $\,C^\infty$ complex line bundle
over a complex manifold $\,N\nh$, and let $\,\Cal H\,$ be the horizontal
distribution of a fixed $\,C^\infty$ linear connection in $\,\Cal L\,$ whose
curvature form $\,\varOmega\,$ (Remark \a\cn.1) is real-val\-ued and
skew-Her\-mit\-i\-an in the sense that
$\,\varOmega(Jv,v')=-\hs\varOmega(v,Jv')\,$ for all $\,y\in N\,$ and
$\,v,v'\in T_yN\nh$. Then $\,\Cal L\,$ admits a unique structure of a
holomorphic line bundle over $\,N\hs$ such that $\,\Cal H\,$ is
$\,J$-invariant as a subbundle of $\,T\Cal L\hs$, where
$\,J:T\Cal L\to T\Cal L\,$ now denotes the complex structure tensor on the
total space $\,\Cal L\hs$.

In fact, let $\,\varGamma\,$ be the connection form (Remark \a\cn.1)
corresponding to a $\,C^\infty$ local trivializing section $\,w\,$ of
$\,\Cal L\hs$, defined on a contractible open set $\,N'\subset N\nh$. Using
$\,w\,$ to identify the portion $\,\Cal L'$ of $\,\Cal L\,$ lying over $\,N'$
with $\,N'\nh\times\bbC$, and writing down the parallel-transport equation in
terms of $\,\varGamma$, we see that, for any $\,(y,\zx)\in\Cal L'$ and
$\,(w,\zeta)\in T_{(y,\zx)}\Cal L'$, the $\,\Cal H\,$ component of
$\,(w,\zeta)\,$ relative to the decomposition $\,T\Cal L=\Cal H\oplus\Cal V\,$
equals $\,(w,-\varGamma(w)\zx)$. Thus, $\,w\,$ is holomorphic for a
holomorphic-bundle structure in $\,\Cal L'$ for which $\,\Cal H\,$ is
$\,J$-invariant if and only if $\,\varGamma\,$ is of type $\,(1,0)$, i.e., the
bundle morphism $\,\varGamma:TN'\to N'\nh\times\bbC\,$ is complex-linear.

Our $\,\varGamma$, with $\,d\varGamma=-\hs i\hs\varOmega\,$ (Remark \a\cn.1)
need not be of type $\,(1,0)$. However, a $\,(1,0)\,$ form
$\,\tilde\varGamma\,$ with $\,d\tilde\varGamma=-\hs i\hs\varOmega\,$ exists on
$\,N'$ since $\,\varOmega\,$ is a closed real-val\-ued form of type $\,(1,1)$,
and so, choosing a function $\,\varphi:N'\to\bbR\,$ with
$\,i\hs\varOmega=\partial\overline{\partial}\varphi$, we may set
$\,\tilde\varGamma=\partial\varphi$. As $\,d\hs(\varGamma-\tilde\varGamma)=0$,
we have $\,\tilde\varGamma=\varGamma+\hs d\varPhi\,$ for some $\,C^\infty$
function $\,\varPhi:N'\to\bbC\hs$. Now $\,\Cal H\,$ is $\,J$-invariant for the
holomorphic-bundle structure in $\,\Cal L'$ obtained by declaring the section
$\,\tilde w=e^\varPhi w\,$ holomorphic. (Note that the connection form
corresponding to $\,\tilde w\,$ is $\,\tilde\varGamma$.) Any other
$\,C^\infty$ section of $\,\Cal L'$ without zeros having a $\,(1,0)\,$
connection form equals $\,e^\varPsi\tilde w$, where $\,\varPsi:N'\to\bbC\,$ is
holomorphic since $\,d\varPsi\,$ is of type $\,(1,0)$. The structure in
question is therefore unique.
\endrmk

\head\S\tb. Tautological bundles\endhead
\rmk{Remark \a\tb.1}Given a Hermitian inner product $\,\langle\,,\rangle\,$
in a complex vector space $\,V$ with $\,\hs\dim V<\infty\,$ and a real
constant $\,a\ne0$, we will denote $\,v,u\,$ the vector fields on $\,V$
given by $\,x\mapsto ax\,$ and $\,x\mapsto aix$, and define two distributions
$\,\Cal V,\Cal H\,$ on $\,V\smallsetminus\{0\}\,$ by
$\,\Cal V=\,\text{\rm Span}\,\{v,u\}\,$ and $\,\Cal H=\Cal V^\perp$. Clearly,
they do not depend on the choice of $\,a$.
\endrmk
For a complex vector space $\,V$ with $\,1\le\hs\dimc V=m<\infty$, let
$\,N\approx\bbCP^{m-1}$ be the projective space of $\,V\nh$, and let
$\,\Cal L\,$ be the tautological bundle over $\,N\nh$. Thus, $\,N\hs$ is a
complex manifold with the underlying set formed by all complex lines through
$\,0\,$ in $\,V\nh$, and $\,\Cal L\,$ is the complex line bundle over $\,N\hs$
whose fibre over any point (i.e., line) is the line itself. We will write
$\,\Cal L\smallsetminus N=V\smallsetminus\{0\}$, that is, identify
$\,\Cal L\smallsetminus N\,$ (which is an open set in the total space
$\,\Cal L$, cf. \f{\cn.2}) with $\,V\smallsetminus\{0\}\,$ using the natural
biholomorphism given, in the notation of \f{\cn.2}, by $\,(y,\zx)\mapsto\zx$.

A fixed Hermitian inner product $\,\langle\,,\rangle\,$ in $\,V$ gives
rise to two objects in $\,\Cal L\hs$. The first is a Hermitian fibre metric,
also denoted $\,\langle\,,\rangle$, and obtained by restricting the inner
product to the fibres of $\,\Cal L\hs$. The second is a canonical connection
in $\,\Cal L\,$ making the fibre metric $\,\langle\,,\rangle\,$ parallel,
whose horizontal distribution, restricted to
$\,\Cal L\smallsetminus N=V\smallsetminus\{0\}$, is $\,\Cal H$, defined in
Remark \a\tb.1. This canonical connection is obtained by projecting the
standard flat connection in the product bundle $\,\Cal E=N\times V\,$ onto
the $\,\Cal L\,$ summand of the direct-sum decomposition
$\,\Cal E=\Cal L\oplus\Cal L^\perp$. Note that its horizontal distribution is
contained in $\,\Cal H\,$ (and hence coincides with $\,\Cal H$, for
dimensional reasons), since any horizontal $\,C^1$ curve
$\,t\mapsto x(t)\in V\smallsetminus\{0\}\,$ has, by definition,
$\,dx/dt\in\Cal H_{x(t)}$ at every $\,t$.

For $\,m\ge2\,$ the curvature form $\,\varOmega\,$ of this canonical
connection (see Remark \a\cn.1) equals $\,-\hs2\hs\omfs$, where $\,\omfs$ is
the K\"ahler form, defined as in \f{\pn.5}, of the quotient (Fubini-Study)
metric on the projective space $\,N\nh$.

In fact, both forms are invariant under the obvious action on
$\,N\approx\bbCP^{m-1}$ of the group $\,G\approx\hs\text{\rm U}\hs(m)\,$ of
all unitary automorphisms of $\,V\nh$. (The original action of $\,G\,$ on
$\,\Cal L\smallsetminus N=V\smallsetminus\{0\}\,$ amounts to a lift of its
action on $\,N\hs$ to $\,\Cal L$, preserving both the fibre metric
$\,\langle\,,\rangle\,$ and the canonical connection, as $\,\Cal H\,$ in
Remark \a\tb.1 is $\,G$-invariant.) Moreover, both forms are real-val\-ued and
skew-Her\-mit\-i\-an (cf. Remark \a\cn.4), so that $\,\omfs\,$ (or,
$\,\varOmega$) is related, as in \f{\pn.5}, to the Fubini-Study metric (or,
respectively, to some  twice-co\-var\-i\-ant symmetric tensor field $\,\bz\,$
on $\,N$). Therefore, $\,\varOmega\,$ is a constant times $\,\omfs$, since
the same is true for $\,\bz\,$ and the metric (as $\,G\,$ leaves them
invariant and acts on $\,N\hs$ with an irreducible isotropy representation).
This constant is $\,-\hs2$, as one sees integrating both forms over a fixed
complex projective line $\,\srf\subset N\nh$, formed by all complex lines
through $\,0\,$ in some complex plane $\,W\subset V\nh$. Namely,
$\,\int_\srf\varOmega\,$ is $\,2\pi\,$ times the integral of $\,c_1(\Cal L)\,$
over the cycle $\,[\srf]$, and, since $\,\Cal L\,$ restricted to $\,\srf\,$ is
the tautological bundle of $\,\srf$, its Chern number, i.e., the latter
integral, equals $\,-\hs1$. (Note that the dual $\,\Cal L^*$ of $\,\Cal L\,$
restricted to $\,\srf\,$ has the Chern number $\,+\hs1$, as its admits a
holomorphic section with one simple zero, obtained by restricting a nonzero
linear functional $\,W\to\bbC\,$ to the fibres.) Finally, as
$\,\int_\srf\omfs$ is the area of $\,\srf\,$ with its own Fubini-Study metric,
it equals $\,\pi$, since $\,\srf\,$ is an orientable Riemannian surface having
a positive constant Gaussian curvature and the diameter $\,\pi/2$.
\rmk{Remark \a\tb.2}The restriction to $\,\Cal H\,$ of the pull\-back of
the Fubini-Study metric on $\,N\nh$ under the standard projection
$\,V\smallsetminus\{0\}\to N\,$ coincides with the Euclidean metric
$\,\,\text{\rm Re}\hskip1pt\langle\,,\rangle\,$ divided by the norm-squared
function on $\,V\smallsetminus\{0\}$. Namely, for any complex-linear operator
$\,\varTheta:V\to V\nh$, the linear vector field $\,x\mapsto\varTheta x\,$ on
$\,V\smallsetminus\{0\}\,$ is projectable onto $\,N\hs$ (since so is its flow,
consisting of linear automorphisms), and such fields, or their projections
onto $\,N\nh$, realize all vectors tangent at any point to
$\,V\smallsetminus\{0\}\,$ (or, $\,N$). The function
$\,x\mapsto\langle\varTheta x,\varTheta x\rangle/\langle x,x\rangle\,$ then is
also projectable onto $\,N\nh$, i.e., homogeneous of degree zero. Fixing
$\,x\in V\smallsetminus\{0\}\,$ and $\,w\in\Cal H_x$, and then choosing
$\,\varTheta\,$ with $\,w=\varTheta x$, we now see that our claim holds at
$\,x$, since, according to the definition of the Fubini-Study metric, it holds
at all $\,x\,$ with $\,\langle x,x\rangle=1$.
\endrmk

\head\S\kg. Killing fields and Killing potentials\endhead
This section contains well-known facts, listed here for easy reference. More
details can be found, for instance, in \cite{\dml}, \S\kp.

A real-val\-ued $\,C^\infty$ function $\,\vp\,$ on a K\"ahler manifold
$\,(M,g)\,$ is called a {\it Killing potential\/} if $\,u=J(\navp)\,$ is a
Killing field on $\,(M,g)$.

As usual, a {\it Killing vector field\/} on a Riemannian
manifold $\,(M,g)\,$ is any $\,C^\infty$ vector field $\,u\,$ such that
$\,\nabla u\,$ is skew-adjoint at every point (cf. \f{\pn.1}).
\rmk{Remark \a\kg.1}A Killing field $\,u\,$ on $\,(M,g)\,$ (or, an isometry
$\,\varPhi:M\to M$) is uniquely determined by $\,u(x)\,$ and
$\,(\nabla u)(x)\,$ (or, by $\,\varPhi(x)\,$ and $\,d\varPhi_x$) at any given
point $\,x\in M$. In fact, using normal coordinates one sees that
$\,\varPhi(x)\,$ and $\,d\varPhi_x$ determine $\,\varPhi\,$ along any broken
geodesic emanating from $\,x$, and the same is true for the local isometries
forming the local flow of $\,u$. This implies a {\it unique continuation
property\/}: an isometry, or a Killing field, is uniquely determined by its
restriction to any nonempty open set.
\endrmk
We call a (real) $\,C^\infty$ vector field $\,v\,$ on a complex manifold {\it
hol\-o\-mor\-phic\/} if $\,\Lie_vJ=0$, where $\,\Lie\,$ is the Lie derivative.
For a $\,C^\infty$ vector field $\,v\,$ on a {\it K\"ahler\/} manifold,
$$\alignedat2
&\text{\rm a)}\quad&&
\text{\rm$v\,$ is holomorphic if and only if $\,[J,\nabla v]=0$,}\\
&\text{\rm b)}\quad&&
\text{\rm$\nabla u=\hs J\circ(\nabla v)$, with the convention \f{\pn.1}, if
$\,u=Jv$.}\hskip32pt\endalignedat\ff\kg.1$$
\proclaim{Lemma \a\kg.2}For a\/ $\,C^\infty$ function\/ $\,\vp\,$ on a
K\"ahler manifold\/ $\,(M,g)$, the following three conditions are
equivalent\/{\rm:}\hskip7pt{\rm(i)}\hskip5pt$\vp\,$ is a \kip\/{\rm;}
\hskip7pt{\rm(ii)}\hskip5ptThe gradient\/ $\,v=\navp\,$ is a holomorphic
vector field\/{\rm;}\hskip7pt{\rm(iii)}\hskip5pt$\bz=\nabla d\vp\,$ is
Hermitian, that is, $\,\bz(Jw,w')=-\hs\bz(w,Jw')\,$ for all\/ $\,x\in M\,$
and\/ $\,w,w'\in T_xM$.{\hfill\qd}
\endproclaim
Next, there is a well-known local one-to-one correspondence between \kip s
defined up to an additive constant, and holomorphic Killing vector fields:
\proclaim{Lemma \a\kg.3}Let\/ $\,(M,g)\,$ be a K\"ahler manifold. For every
\kip\/ $\,\vp\,$ on\/ $\,(M,g)$, the Killing field\/ $\,J(\navp)\,$ is
hol\-o\-mor\-phic. Conversely, if\/ $\,H^1(M,\bbR)=\{0\}$, then every
holomorphic Killing vector field on\/ $\,(M,g)\,$ has the form\/
$\,J(\navp)\,$ for a \kip\/ $\,\vp$, which is unique up to an additive
constant.{\hfill\qd}
\endproclaim
\rmk{Remark \a\kg.4}Let $\,\vp\,$ be a nonconstant \kip\ on a K\"ahler
manifold $\,(M,g)$. Then $\,\nabla d\vp\ne0\,$ wherever $\,d\vp=0\,$ (and
hence $\,\,d\vp\ne0\,$ on a dense open subset of $\,M\nh$, which also follows
from Lemma \a\kg.2(ii)). In fact, if $\,\nabla d\vp\,$ and $\,d\vp\,$ both
vanished at some point, so would $\,v=\navp\,$ and $\,\nabla v\,$ (by
\f{\pn.2}), as well as $\,u=Jv\,$ and $\,\nabla u\,$ (by \f{\kg.1.b}). The
Killing field $\,u=J(\navp)\,$ thus would vanish identically on $\,M\,$ (see
Remark \a\kg.1), contradicting nonconstancy of $\,\vp$.
\endrmk

\head\S\cz. Critical manifolds\endhead
The following two lemmas are well-known. For details, see, e.g., \cite{\dml},
\S\cm.
\proclaim{Lemma \a\cz.1}Let\/ $\,u\,$ be a Killing field on a Riemannian
manifold\/ $\,(M,g)$, and let\/ $\,y\in M\,$ be a point such that\/
$\,u(y)=0$. Then, for every sufficiently small\/ $\,\,\text{\rm d}>0$, the
flow of\/ $\,u\,$ restricted to the radius\/ $\,\,\text{\rm d}\,$ open ball\/
$\,\,U\,$ centered at\/ $\,y\,$ consists of ``global'' isometries\/
$\,\,U\to\,U\nh$, while the exponential mapping\/ $\,\,\e_{\hskip.4pty}$ is
defined everywhere in the open ball\/ $\,\,U'$ of radius\/
$\,\,\text{\rm d}\,$ in\/ $\,T_yM$, centered at\/ $\,0$, and\/
$\,\,\e_{\hskip.4pty}$ maps\/ $\,\,U'$ diffeomorphically onto\/
$\,\,U\nh$.{\hfill\qd}
\endproclaim
\proclaim{Lemma \a\cz.2}For a Killing vector field\/ $\,u\,$ on a Riemannian
manifold\/ $\,(M,g)$, let\/ $\,N(u)=\{y\in M:u(y)=0\}\,$ be the set of all
zeros of\/ $\,u$. If\/ $\,u\ne0\,$ somewhere in\/ $\,M\nh$, then, for every
connected component\/ $\,N\hs$ of\/ $\,N(u)$, with\/ $\,\nabla u\,$ as in\/
\f{\pn.1},
\widestnumber\item{(d)}\roster
\item"(a)"$N\,$ is contained in an open set that does not intersect any other
component.
\item"(b)"$N\subset M\,$ is a closed set and a submanifold with the subset
topology.
\item"(c)"The submanifold\/ $\,N\hs$ is totally geodesic in\/ $\,(M,g)\,$
and\/ $\,\dim M-\hs\dim N\ge2$.
\item"(d)"For any\/ $\,y\in N\,$ we have
$\,T_yN=\,\text{\rm Ker}\,[(\nabla u)(y)]=\{w\in T_yM:\nabla_{\!w}u=0\}$.
\endroster
Furthermore, the set\/ $\,M'=M\smallsetminus N(u)\,$ is connected, open and
dense in\/ $\,M$.{\hfill\qd}
\endproclaim
Let $\,\vp:M\to\bbR\,$ be a $\,C^1$ function on a manifold $\,M$.
If all connected components $\,N\hs$ of the set
$\,\,\text{\rm Crit}\hs(\vp)\,$ of its critical points happen to satisfy
conditions (a), (b) of Lemma \a\cz.2, we will refer to them as the {\it
critical manifolds\/} of $\,\vp$.
\rmk{Remark \a\cz.3}Let $\,\vp:M\to\bbR\,$ be a nonconstant \kip\ (\S\kg)
on a K\"ahler manifold $\,(M,g)$, and let $\,M'\subset M\,$ be the open set on
which $\,d\vp\ne0$. Then
\widestnumber\item{(iii)}\roster
\item"(i)"$M'$ is connected and dense in $\,M$.
\item"(ii)"The connected components of $\,\,\text{\rm Crit}\hs(\vp)\,$ are
totally geodesic submanifolds of $\,(M,g)\,$ (the {\it critical manifolds\/}
of $\,\vp$), satisfying (a) -- (d) in Lemma \a\cz.2.
\item"(iii)"Every critical manifold $\,N\hs$ of $\,\vp\,$ is a complex
submanifold of $\,M\nh$, and
$\,T_yN=\,\text{\rm Ker}\,[(\nabla v)(y)]=\{w\in T_yM:\nabla_{\!w}v=0\}\,$
for any $\,y\in N\nh$.
\endroster
(All three conclusions are well-known.) In fact, $\,u=J(\navp)\,$ is a Killing
field, and so (i), (ii) are obvious from Lemma \a\cz.2 with
$\,N(u)=\,\text{\rm Crit}\hs(\vp)$, i.e., $\,M'=M\smallsetminus N(u)$.
Finally, $\,N\hs$ in (iii) is a complex submanifold since every $\,T_yN\nh$,
$\,y\in N\nh$, is $\,J$-invariant (by Lemma \a\cz.2(d), as $\,J\,$ and
$\,\nabla u\,$ commute, cf. Lemma \a\kg.3 and \f{\kg.1.a}), while, for
$\,v=\navp$, \f{\kg.1.b} gives $\,\nabla u=J\circ(\nabla v)$, so that
$\,\,\text{\rm Ker}\,[(\nabla u)(y)]=\,\text{\rm Ker}\,[(\nabla v)(y)]\,$ for
$\,y\in M$, and the formula for $\,T_yN\,$ follows from Lemma \a\cz.2(d).
\endrmk

\head\S\gd. Geodesic vector fields\endhead
Let $\,\nabla\,$ be a fixed connection in the tangent bundle $\,TM\,$ of a
manifold $\,M$, for instance, the Levi-Civita connection of some Riemannian
metric on $\,M$. By a {\it geodesic vector field\/} on $\,M\,$ we mean any
$\,C^\infty$ vector field $\,v\,$ on $\,M\,$ such that
$$\nabla_{\!v}v\,=\,\ta\hskip.4ptv\hskip10pt\text{\rm for\ some\ function}
\hskip8pt\ta:M\to\bbR\hs.\ff\gd.1$$
The function $\,\ta\,$ is not even required to be continuous: its values at
points where $\,v=0\,$ are not determined by \f{\gd.1}, and may be completely
arbitrary.
\rmk{Remark \a\gd.1}Condition \f{\gd.1} holds if and only if every
integral curve $\,\sa\mapsto x(\sa)\,$ of $\,v\,$ is a (re-parameterized)
geodesic of $\,\nabla$. In fact, since $\,\dot x(\sa)=v(x(\sa))\,$ (where
$\,\dot x=\hs dx/d\sa$), it follows that $\,\nabla_{\!\dot x}\dot x\,$ at
any $\,\sa\,$ equals $\,\nabla_{\!v}v\,$ at $\,x(\sa)$, while curves
$\,\sa\mapsto x(\sa)\,$ obtained from geodesics by diffeomorphic changes of
parameter are characterized by
$\,\nabla_{\!\dot x}\dot x=\ta\hskip.4pt\dot x\,$ with $\,\ta$, this time,
denoting a function of $\,\sa$.
\endrmk
\rmk{Remark \a\gd.2}Let $\,\nabla\,$ be a connection in the tangent bundle
$\,TM$. If $\,v\,$ is a $\,C^\infty$ vector field with \f{\gd.1} and
$\,X\subset M\,$ is a geodesic segment such that $\,v(x)\,$ is tangent to
$\,X\hs$ at some point $\,x\in X\hs$ and $\,v\ne0\,$ at all points of $\,X$,
then $\,v\,$ is tangent to $\,X\hs$ at every point of $\,X$. In fact, both
$\,X\hs$ and the underlying set $\,\tilde X\hs$ of the maximal integral curve
of $\,v\,$ containing $\,x\,$ are geodesics (Remark \a\gd.1), tangent to each
other at $\,x$, and so $\,x\in X'\subset\tilde X\hs$ for some nontrivial
subsegment $\,X'$ of $\,X$. Choosing $\,X'$ to be the maximal subsegment with
this property, we must have $\,X'=X$, for otherwise an endpoint $\,x'$ of
$\,X'$ would be an interior point of $\,X\hs$ and $\,v(x')\ne0\,$ would be
tangent to $\,X\hs$ at $\,x'$, thus allowing $\,X'$ to be extended past
$\,x'$, contrary to maximality.
\endrmk
\rmk{Remark \a\gd.3}Given a connection $\,\nabla\,$ in the tangent bundle
$\,TM$, we use the standard symbol $\,\,\e_{\hskip.4ptx}:\hs U_x\to M\,$ for
the geodesic exponential mapping of $\,\nabla\,$ at any point $\,x\in M$. Here
$\,\,U_x$ is a neighborhood of the zero vector in $\,T_xM$, obtained as a
union of maximal line segments emanating from zero on which
$\,\,\e_{\hskip.4ptx}$ is defined. Thus,
$\,\sa\mapsto x(\sa)=\,\e_{\hskip.4ptx}\hs\sa w\,$ is the geodesic with
$\,x(0)=x\,$ and $\,\dot x(0)=w\in T_xM$. A related mapping is
$\,\,\text{\rm Exp}:\hs\uexp\to M\,$ with
$\,\,\text{\rm Exp}\hs(x,w)=\e_{\hskip.4ptx}\hs w\hs$, defined on the subset
$\,\,\uexp=\hs\bigcup_{x\in M}[\{x\}\times\,U_x]\,$ of the total space
$\,TM=\{(x,w):x\in M,\hskip6ptw\in T_xM\}$, containing the zero section
$\,M\subset TM\,$ (cf. \f{\cn.2}). The set $\,\,\uexp$ is open in $\,TM$, and
$\,\,\text{\rm Exp}\,\,$ is of class $\,C^\infty$ (see, e.g., \cite{\kno}, p.
147).
\endrmk
\proclaim{Lemma \a\gd.4}Suppose that\/ $\,\nabla\,$ is a connection in the
tangent bundle\/ $\,TM\,$ of a manifold\/ $\,M\,$ and\/ $\,v\,$ is a\/
$\,C^\infty$ vector field on\/ $\,M\,$ with\/ \f{\gd.1}, while\/
$\,X\subset M\,$ is a geodesic segment containing an endpoint\/
$\,y\,$ with\/ $\,v(y)=0$. If\/ $\,\nabla_{\!w}v=aw\,$ for some nonzero
vector\/ $\,w\,$ tangent to\/ $\,X\hs$ at\/ $\,y\,$ and some\/
$\,a\in\bbR\smallsetminus\{0\}$, then
\widestnumber\item{(b)}\roster
\item"(a)"There exists a nontrivial compact subsegment\/ $\,X'$ of $\,X$,
containing\/ $\,y$, and such that\/ $\,v(x)\ne0\,$ for all\/
$\,x\in X'\smallsetminus\{y\}$.
\item"(b)"For any subsegment\/ $\,X'\subset X\hs$ with the properties listed
in\/ {\rm(a)} we have\/ $\,v(x)\in T_xX\hs$ at every\/ $\,x\in X'$.
\endroster
\endproclaim
\demo{Proof}Let $\,\sa\mapsto x(\sa)\,$ be a geodesic parameterization of
$\,X\hs$ with $\,x(0)=y$, defined on a subinterval of $\,[\hs0,\infty)$. Thus,
$\,\nabla_{\!\dot x}\dot x=0$, where $\,\dot x=\hs dx/d\sa$. A fixed
$\,1$-form $\,\xi\,$ of class $\,C^\infty$ on a neighborhood $\,\,U\,$ of
$\,y\,$ such that $\,a\xi(w)>0\,$ at $\,y$, for $\,w=\dot x(0)$, gives
rise to a $\,C^\infty$ function $\,\varphi=\xi(v):\hs U\to\bbR\,$ with
$\,\varphi=0\,$ wherever $\,v=0\,$ in $\,\,U\,$ and
$\,\hs d\hs[\varphi(x(\sa))]/d\sa>0\,$ for all $\,\sa\ge0\,$ near $\,0\,$ (as
$\,d_w\varphi=a\xi(w)$), which proves (a).

For $\,X'$ as in (a), let $\,\ell>0\,$ be such that $\,x(\ell)\,$ is an
endpoint of $\,X'$, and let $\,\sa\mapsto w(\sa)\in T_{x(\sa)}M\,$ be the
vector field along $\,X'$ given by $\,w(0)=\dot x(0)\,$ and
$\,w(\sa)=v(x(\sa))/f(\sa)\,$ for $\,\sa\in(0,\ell\hs]$, where
$\,f:[\hs0,\ell\hs]\to\bbR\,$ is any fixed $\,C^1$ function with $\,f(0)=0$,
$\,\dot f(0)=a$ and $\,|f|>0\,$ on $\,(0,\ell\hs]$. Thus, $\,w(\sa)\ne0\,$ for
all $\,\sa\in[\hs0,\ell\hs]\,$ due to our choice of $\,X'$ and $\,\ell$. Also,
setting $\,\tilde v(\sa)=v(x(\sa))\,$ we have $\,\tilde v(0)=0$, while
$\,\nabla_{\!\dot x}\tilde v\,$ at $\,\sa=0\,$ equals $\,\nabla_{\!w}v=aw$,
with $\,w=\dot x(0)$. This, along with l'Hospital's rule, shows that the
mapping $\,[\hs0,\ell\hs]\ni\sa\mapsto(x(\sa),w(\sa))$, valued in the total
space $\,TM\,$ (see Remark \a\gd.3), is continuous, also at $\,\sa=0$.

For any fixed $\,\sa\in[\hs0,\ell\hs]$, let $\,\dg\mapsto x_\sa(\dg)\in M\,$
be the geodesic with $\,x_\sa(\sa)=x(\sa)\,$ and
$\,d\hs[x_\sa(\dg)/d\dg]_{\dg=\sa}=w(\sa)$, defined on the maximal possible
interval containing $\,\sa$. Then, for any sufficiently small
$\,\ve\in(0,\ell\hs]$,
\widestnumber\item{(ii)}\roster
\item"(i)"$\,\dg\mapsto x_\sa(\dg)\,$ is defined on an interval containing
$\,[\hs0,\ell\hs]$, for every $\,\sa\in[\hs0,\ve]$.
\item"(ii)"$\,v\ne0\,$ at $\,x_\sa(\dg)\,$ for any $\,\sa,\dg\,$ with
$\,0<\sa\le\dg\le\ve$.
\endroster 
In fact, if there were no $\,\ve\in(0,\ell\hs]\,$ with (i), we could find
values of $\,\sa\in(0,\ell\hs]\,$ arbitrarily close to $\,0\,$ such that one
of the points $\,(x(\sa),-\hs\sa w(\sa))$, $\,(x(\sa),(\ell-\sa)\hs w(\sa))\,$
lies in the complement $\,TM\smallsetminus\,\uexp$, with $\,\,\uexp$ as in
Remark \a\gd.3. Since $\,TM\smallsetminus\,\uexp$ is a closed set, it would
then also contain the limit of one of these points as $\,s\to0$, i.e.,
$\,(y,0)\,$ or $\,(y,\ell\hs\dot x(0))$, contradicting either the fact that
$\,M\subset\,\uexp$, or our choice of $\,\ell$.

Also, $\,\hs d\hs[\varphi(x_\sa(\dg))]/d\dg]>0\,$ for all sufficiently small
$\,\dg,\sa\in[\hs0,\ell\hs]\,$ since, due to our choice of $\,\varphi$, this
is the case for $\,\sa=0\,$ and $\,\dg\ge0\,$ close to $\,0$. (Note that
$\,x_0(\dg)=x(\dg)$.) As $\,\varphi>0\,$ on a nontrivial subsegment of $\,X'$
containing $\,y$, except for the point $\,y\,$ at which $\,\varphi=0$,
by making $\,\ve>0\,$ with (i) smaller we thus have $\,\varphi>0\,$ (and hence
$\,v\ne0$) at $\,x_\sa(\dg)\,$ for any $\,\sa,\dg\,$ with
$\,\sa\in(0,\ve]\,$ and $\,\sa\le\dg\le\sa+\ve$, which gives (ii).

By (i), (ii) and Remark \a\gd.2, for every $\,\sa\in(0,\ve)$, the geodesic
$\,[\sa,\ell\hs]\ni\dg\mapsto x_\sa(\dg)\,$ is a (re-parameterized) integral
curve of $\,v\,$ and, in particular, $\,v\,$ is tangent to it at the point
$\,x_\sa(\ve)$. Taking the limit as $\,\sa\to0$, we now see that $\,v\,$ is
tangent to the limiting geodesic, i.e., to $\,X'$, at $\,x(\ve)$. Applying
Remark \a\gd.2 to $\,x=x(\ve)$, we obtain (b) for our $\,X'$, which completes
the proof.{\hfill\qd}
\enddemo
\proclaim{Lemma \a\gd.5}Let a Killing field\/ $\,u\,$ on a Riemannian
manifold\/ $\,(M,g)\,$ vanish at a point\/ $\,y\,$ and let\/
$\,\zx\in T_yM\,$ lie in the set\/ $\,\,U_y$ defined as in\/ {\rm Remark
\a\gd.3} for the Levi-Civita connection\/ $\,\nabla$. Then, at the point\/
$\,x=\,\text{\rm Exp}\hs(y,\zx)$, the vector\/ $\,u(x)\,$ is the image of\/
$\,\nabla_{\!\zx}u\in T_yM=T_\zx(T_yM)\,$ under the differential of\/
$\,\,\e_{\hskip.4pty}$ at\/ $\,\zx$.
\endproclaim
In fact, the local isometries $\,\varPhi^t$ forming the local flow of $\,u\,$
are all defined, for $\,t\,$ near $\,0\,$ in $\,\bbR\hs$, on some open set in
$\,M\,$ containing the compact geodesic segment
$\,X=\{\hs\e_{\hskip.4pty}\hs\sa \zx:0\le\sa\le1\}$. Since they keep $\,y\,$
fixed and map geodesics onto geodesics, we have
$\,\varPhi^t(\e_{\hskip.4pty}\hs \zx)=\,\e_{\hskip.4pty}\hs(d\varPhi^t_y\zx)$.
Our claim follows if we apply $\,d/dt\,$ and let $\,t\to0$, since
$\,(\nabla u)(y):T_yM\to T_yM\,$ (cf. \f{\pn.1}) is the infinitesimal
generator of the one-parameter group $\,t\mapsto\hs d\varPhi^t_y$ in
$\,T_yM$.{\hfill\qd}

\head\S\mb. Morse-Bott functions\endhead
By a {\it Morse-Bott function\/} on a manifold $\,M\,$ (cf. \cite{\btt}) we
mean any $\,C^\infty$ function $\,\vp:M\to\bbR\,$ such that every connected
component $\,N\hs$ of the set $\,\,\text{\rm Crit}\hs(\vp)\,$ of its critical
points satisfies conditions (a), (b) of Lemma \a\cz.2 and, for every
$\,y\in N\nh$, the nullspace of the Hessian
$\,\,\text{\rm Hess\hskip.3pt}_y\vp\,$ coincides with $\,T_yN\nh$. (Since the
nullspace contains $\,T_yN\,$ for any submanifold
$\,N\subset\,\text{\rm Crit}\hs(\vp)$, the last requirement amounts to
$\,\,\text{\rm rank}\hskip3pt\text{\rm Hess\hskip.3pt}_y\vp
=\hs\dim M-\hs\dim N\nh$.) The connected components of
$\,\,\text{\rm Crit}\hs(\vp)\,$ then are called the {\it critical manifolds\/}
of $\,\vp$, cf. \S\cz.
\example{Example \a\mb.1}Every \kip\ $\,\vp\,$ on a K\"ahler manifold
$\,(M,g)\,$ (\S\kg) is a Morse-Bott function: this is clear when $\,\vp\,$ is
constant, while for nonconstant $\,\vp\,$ it follows from Remark \a\cz.3(ii),
(iii) along with \f{\pn.2}. (Note that
$\,\,\text{\rm Hess\hskip.3pt}_y\vp=(\nabla d\vp)(y)\,$ whenever
$\,y\in\,\text{\rm Crit}\hs(\vp)$.)
\endexample
\proclaim{Lemma \a\mb.2}Let\/ $\,\vp\,$ be a Morse-Bott function on a
manifold\/ $\,M$. Every point\/ $\,y\,$ of any critical manifold\/ $\,N\hs$
of\/ $\,\vp\,$ at which\/ $\,\,\text{\rm Hess\hskip.3pt}_y\vp\,$ is
positive/negative semidefinite then has a neighborhood\/ $\,\,U$ such that\/
$\,\vp>\vp(y)\,$ everywhere in\/ $\,\,U\smallsetminus N\,$ or, respectively,
$\,\vp<\vp(y)\,$ everywhere in\/ $\,\,U\smallsetminus N\nh$.
\endproclaim
\demo{Proof}We may assume that $\,M=\hs\ren$, $\,y=(0,\dots,0)\,$ and
$\,N=\{(0,\dots,0)\}\times\hs\renr$. Everywhere in $\,N\hs$ we thus have
$\,\hs\partial_a\vp=\hs\partial_\lambda\vp=0\,$ and hence
$\,\hs\partial_\lambda\partial_\mu\vp=\hs\partial_\lambda\partial_a\vp
=\hs\partial_a\partial_\lambda\vp=0\,$ for all $\,a,b\in\{1,\dots,r\}\,$
and $\,\lambda,\mu\in\{r+1,\dots,n\}$, where
$\,\hs\partial_j=\hs\partial/\partial x^j$. Since
$\,\,\text{\rm rank}\hskip3pt\text{\rm Hess\hskip.3pt}_y\vp=r\,$ and
$\,\dim N=n-r$, semidefiniteness of $\,\,\text{\rm Hess\hskip.3pt}_y\vp\,$ now
implies that, at $\,y=(0,\dots,0)$, the
$\,r\times r\hs$ matrix $\,[\hs\partial_a\partial_b\vp]\,$ is
positive/negative definite.

For any fixed vectors $\,\,\text{\bf u}\,=(u^{r+1},\dots,u^n)\,$ close to
$\,\,\text{\bf0}\,=(0,\dots,0)\,$ and $\,\,\text{\bf w}\,=(w^1,\dots,w^r)\,$
with $\,|\text{\bf w}|=1$, we will write
$\,\dot\vp=\hs d\hs[\vp(x(\sa))]/d\sa$, where the curve
$\,\sa\mapsto x(\sa)\in\,U\,$ is given by
$\,x(\sa)=(\sa\text{\bf w},\,\text{\bf u})$, that is, $\,x^a(\sa)=\sa w^a$,
$\,\,x^\lambda(\sa)=u^\lambda$ for $\,\sa\,$ near $\,0\in\bbR\,$ and
$\,a,b,\lambda,\mu\,$ as before. If no neighborhood $\,\,U\,$
of $\,\,\text{\bf0}\,\,$ in $\,\ren$ (i.e., of $\,y\,$
in $\,M$) had the required property, we could find sequences
$\,\,\text{\bf u}_k,\hs\text{\bf w}_k$ of such vectors
with $\,\hs\text{\bf u}_k\to\hs\text{\bf0}\hs\,$ and
$\,\hs\text{\bf w}_k\to\hs\text{\bf w}\hs\ne\hs\text{\bf0}\hs\,$ as
$\,k\to\infty$, with $\,\pm\hs[\vp(x_k(\sa_k))-\vp(y)]\le0\,$ for
$\,x_k(\sa)=(\sa\hs\text{\bf w}_k,\,\text{\bf u}_k)\,$ and some sequence
$\,\sa_k$ converging to $\,0\,$ in $\,\bbR\hs$. However, as
$\,\dot\vp=w^a\hs\partial_a\vp\,$ (summed over $\,a=1,\dots,r$), at
$\,\sa=0\,$ we have $\,\vp=\vp(y)\,$ and $\,\dot\vp=0$. The mean value theorem
for $\,\sa\mapsto\vp(x_k(\sa))\,$ and
$\,\dot\vp_k=\hs d\hs[\vp(x_k(\sa))]/d\sa\,$ then would give
$\,\pm\hs\ddot\vp_k\le0\,$ at $\,\sa=\tilde\sa_k$ for some
$\,\tilde\sa_k$ with $\,\tilde\sa_k\to0\,$ as $\,k\to\infty$. Since
$\,\ddot\vp=w^aw^b\hs\partial_a\partial_b\vp$, it would follow that, in the
limit, $\,\pm\hs w^aw^b\hs\partial_a\partial_b\vp\le0\,$ at
$\,y=\,\text{\bf0}$, contrary to positive/negative definiteness of
$\,[\hs\partial_a\partial_b\vp]$. This completes the proof.{\hfill\qd}
\enddemo
\proclaim{Lemma \a\mb.3}Let there be given a\/ $\,C^\infty$ vector field\/
$\,v\,$ without zeros on a manifold\/ $\,M'\nh$, an interval\/ $\,\iy\,$ of
the variable\/ $\,\vp$, a sequence\/ $\,\varPhi_k:\iy\to M'$, $\,k=1,2,\dots$,
of integral curves of\/ $\,v$, and a family\/ $\,\{P[\vp]:\vp\in\iy\}\,$ of
compact sets\/ $\,P[\vp]\subset M'$ with\/ $\,\varPhi_k(\vp)\in P[\vp]\,$ for
every\/ $\,k=1,2,\dots\,$ and\/ $\,\vp\in\iy$. If, for some\/ $\,\vp\in\iy$,
the sequence\/ $\,\varPhi_k(\vp)\,$ converges in\/ $\,M'\nh$, then there
exists an integral curve\/ $\,\varPhi:\iy\to M'$ of\/ $\,v\,$ such that\/
$\,\varPhi_k(\vp)\to\varPhi(\vp)\,$ as\/ $\,k\to\infty\,$ for every\/
$\,\vp\in\iy$.
\endproclaim
\demo{Proof}The set $\,\jy\,$ of those $\,\vp\in\iy\,$ for which
$\,\varPhi(\vp)=\,\lim_{\,k\to\infty}\varPhi_k(\vp)\,$ exists is relatively
open in $\,\iy\,$ and the restriction of the mapping $\,\varPhi:\jy\to M'$
just defined to any interval contained in $\,\jy\,$ is an integral curve of
$\,v$. This becomes clear if one fixes $\,\vp\in\jy\,$ and chooses local
coordinates at $\,\varPhi(\vp)\,$ in which $\,v\,$ is a coordinate vector
field.

All we still need to show is that $\,\jy_0=\iy$, where the interval $\,\jy_0$
is any fixed connected component of $\,\jy$. To this end, let us assume, on
the contrary, that $\,\vp\hs'\in\iy\smallsetminus\jy_0$ is an endpoint of
$\,\jy_0$. Since $\,P[\vp\hs']\,$ is compact, the $\,\varPhi_k(\vp\hs')\,$
have a convergent subsequence. Let $\,\tilde{\jy}\,$ be the analogue of
$\,\jy\,$ for that subsequence, and let $\,\tilde{\jy}_0$ be the
connected component of $\,\tilde{\jy}\,$ containing $\,\vp\hs'$. Applying the
last paragraph to our subsequence and $\,\tilde{\jy}$, we see that some
integral curve $\,\varPhi\hs':\tilde{\jy}_0\to M'$ of $\,v\,$ is the
pointwise limit, on $\,\tilde{\jy}_0$, of the $\,\varPhi_k$ in our
subsequence. Hence $\,\varPhi\hs'$ must agree with $\,\varPhi$, the pointwise
limit of the original sequence $\,\varPhi_k$, on the (obviously nonempty) open
interval $\,\jy_0\cap\tilde{\jy}_0$. 
All convergent subsequences of the original sequence
$\,\varPhi_k(\vp\hs')\,$ therefore have the same limit
$\,\varPhi\hs'(\vp\hs')=\,\lim_{\,\vpsu\to\vpsu\hs'}\varPhi(\vp)$, i.e. (as
$\,\varPhi_k(\vp\hs')\in P[\vp\hs']\,$ and $\,P[\vp\hs']\,$ is compact), the
sequence $\,\varPhi_k(\vp\hs')\,$ itself converges. Hence $\,\vp\hs'$ lies in
$\,\jy_0$ instead of $\,\iy\smallsetminus\jy_0$. This contradiction completes
the proof.{\hfill\qd}
\enddemo
\proclaim{Lemma \a\mb.4}Let\/ $\,\vp\,$ be a\/ $\,C^\infty$ function on a
manifold\/ $\,M'$ such that the\/ $\,\vp$-pre\-im\-age of every real number
is compact and\/ $\,\vp\,$ has no critical points.
\widestnumber\item{(ii)}\roster
\item"(i)"There exist a compact manifold\/ $\,P\,$ and a diffeomorphic
identification\/ $\,M'=P\times(\vp_-,\vp_+)\,$ under which\/ $\,\vp\,$ appears
as the projection onto the $\,(\vp_-,\vp_+)\,$ factor, $\,\,\vp_-$ and
$\,\vp_+$ being the infimum and supremum of\/ $\,\vp$.
\item"(ii)"The\/ $\,\vp$-pre\-im\-age of every real number is both compact and
connected.
\endroster
\endproclaim
\demo{Proof}Let us choose a Riemannian metric $\,g\,$ on $\,M'$ with
$\,g(v,v)=1$, where $\,v=\navp\,$ is the $\,g$-gra\-di\-ent of $\,\vp$. (A
unique metric with this property exists in every conformal class.) Also, let
$\,\,U\,$ be the union of all maximal integral curves of $\,v=\navp\,$ that
intersect a fixed connected component $\,P\,$ of a (nonempty)
$\,\vp$-pre\-im\-age of a given real number $\,\hat\vp$.

For any $\,x\in P$, we denote
$\,\varPhi(x,\,\cdot\,):(\vp_-(x),\vp_+(x))\to\,U\,$ the maximal
(pa\-ram\-e\-triz\-ed) integral curve of $\,v\,$ with
$\,\varPhi(x,\hat\vp)=x$, for $\,\hat\vp\,$ as above.
Since $\,d_v\vp=g(v.\navp)=1$, the natural parameter of every integral curve
of $\,v\,$ coincides, up to an additive constant, with $\,\vp$. Thus,
$\,\varPhi(x,\,\cdot\,)\,$ is parameterized by $\,\vp$, i.e., for any real
number $\,\vp\hs'\in(\vp_-(x),\vp_+(x))$, the value of $\,\vp\,$ at
$\,\varPhi(x,\vp\hs')\,$ is $\,\vp\hs'$. If we now denote $\,\,U'$ the union
of all $\,\{x\}\times(\vp_-(x),\vp_+(x))\,$ with $\,x\in P$, then
$\,\,U'$ is an open set in
$\,P\times(\vp_{\text{\rm inf}},\vp_{\text{\rm sup}})\,$ and the
mapping $\,\varPhi\,$ defined above is a $\,C^\infty$ diffeomorphism
$\,\,U'\!\to\,U\nh$. Note that $\,\,U\,$ is an open subset of $\,M'$, as it is the
union of the images under the local flow of $\,v\,$ of sufficiently small
neighborhoods in $\,M'$ of points of $\,P$.

For any given closed interval $\,[\vp_0,\vp_1]\,$ containing $\,\hat\vp$, the
set of those $\,x\in P\,$
for which $\,[\vp_0,\vp_1]\subset(\vp_-(x),\vp_+(x))\,$ is both open and
closed in $\,P$. In fact, its openness follows from that of $\,\,U'$ (see
above), while its closedness is immediate from Lemma \a\mb.3 applied to any
open interval $\,\iy\subset[\vp_0,\vp_1]$, the compact sets $\,P[\vp\hs']\,$
being the $\,\vp$-pre\-im\-ages of real numbers $\,\vp\hs'$.

Consequently, as $\,P\,$ is connected, $\,\vp_\pm=\vp_\pm(x)\,$ do not depend
on $\,x\in P$, i.e., $\,\,U'=P\times(\vp_-,\vp_+)$. Thus, $\,\,U=M'$, since
the open set $\,\,U=\varPhi(U')\,$ is also closed in $\,M'$. Namely, if a
sequence $\,\varPhi(x_k,\vp_k)\,$ converges to $\,x\in M'$, we have
$\,\vp_k=\vp(\varPhi(x_k,\vp_k))\to\vp(x)\,$ as $\,k\to\infty\,$ and, since
$\,P\,$ is compact, the $\,x_k$ contain a subsequence converging to some
$\,x_\infty\in P$, with $\,x=\varPhi(x_\infty,\vp(x))\in\varPhi(U')=\hs U\nh$.
This proves (i). Now (ii) follows, as the pre\-im\-ages in question are
connected (being either empty or diffeomorphic to $\,P$), which completes the
proof.{\hfill\qd}
\enddemo
\proclaim{Corollary \a\mb.5}Let\/ $\,\vp\,$ be a Morse-Bott function on a
compact manifold\/ $\,M\,$ such that\/ $\,\,\text{\rm Hess\hskip.3pt}_y\vp\,$
is semidefinite for every\/ $\hs\,y\in\,\text{\rm Crit}\hs(\vp)$, and the real
codimensions of all critical manifolds of\/ $\,\hs\vp\,$ are greater than one.
Then\/ $\,\vp\,$ has exactly two critical manifolds, which are the\/
$\,\vp$-pre\-im\-ages of its extremum values\/ $\,\vp_+=\vp_{\text{\rm max}}$
and\/ $\,\vp_-=\vp_{\text{\rm min}}$, and the\/ $\,\vp$-pre\-im\-age of every
real number is both compact and connected.
\endproclaim
\demo{Proof}Since $\,M\,$ is compact, $\,\vp\,$ has only finitely many
critical manifolds, due to their discreteness property analogous to (a) in
Lemma \a\cz.2. None of them dis\-con\-nects $\,M\nh$, even locally, in view of
the codimension condition, and $\,\vp\,$ restricted to each of them is
constant. Therefore, the open set
$\,M'=M\smallsetminus\,\text{\rm Crit}\hs(\vp)\,$ of all noncritical points of
$\,\vp\,$ is connected and dense in $\,M\nh$, and its $\,\vp$-im\-age
$\,\vp(M')\,$ is connected, open in $\,\bbR\hs$, and dense in
$\,[\vp_-,\vp_+]$, so that
$\,\vp(M')=(\vp_-,\vp_+)$.

Moreover, the function $\,\vp:M'\to\bbR\,$ satisfies the hypotheses of Lemma
\a\mb.4. In fact, any sequence of points in $\,M'$ that lies in the
$\,\vp$-pre\-im\-age of a given real number has a subsequence converging to a
limit $\,y\in M$, and then $\,y\notin\,\text{\rm Crit}\hs(\vp)$, for otherwise
our semidefiniteness assumption, combined with Lemma \a\mb.2, would lead to a
contradiction. Hence the assertion of Lemma \a\mb.4(ii) holds for
$\,\vp:M'\to\bbR\hs$.

The only critical values of $\,\vp:M\to\bbR\,$ are $\,\vp_\pm$. In fact, let
$\,y\in\,\text{\rm Crit}\hs(\vp)$. Denseness of $\,M'$ in $\,M\,$ gives
$\,x_k\to y\,$ as $\,k\to\infty\,$ for some sequence $\,x_k$ in $\,M'$. If we
had $\,\vp(y)\in(\vp_-,\vp_+)$, the sequence in $\,P\times(\vp_-,\vp_+)\,$
corresponding to the $\,x_k$ under the identification of Lemma \a\mb.4(i)
(applied to $\,\vp:M'\to\bbR\hs$, with $\,\vp(M')=(\vp_-,\vp_+)$) would have a
convergent subsequence, i.e., a subsequence of the $\,x_k$ would have a limit
in $\,M'$, even though $\,x_k\to y\notin M'$.

We thus have connectedness of the $\,\vp$-pre\-im\-ages of all real
$\,\vp\hs'\notin\{\vp_+,\vp_-\}$. To see that the $\,\vp$-pre\-im\-ages
$\,P[\vp_\pm]\,$ of $\,\vp_\pm$ are connected as well, let us fix
$\,\vp_0\in\{\vp_-,\vp_+\}\,$ and denote $\,N_1,\dots,N_l$ the connected
components of $\,P[\vp_0]$. Also, let $\,\,U_1,\dots,\hs U_l$ be pairwise
disjoint open sets in $\,M\,$ with
$\,N_j=\hs U_j\cap\,\text{\rm Crit}\hs(\vp)\,$ for $\,j=1,\dots,l$. The
$\,\vp$-pre\-im\-age $\,P[\vp\hs']\,$ of every $\,\vp\hs'\in(\vp_-,\vp_+)\,$
sufficiently close to $\,\vp_0$ must now be
contained in the union $\,\,U=\,U_1\cup\dots\cup\hs U_l$, or else there
would be a sequence $\,x_k$ in $\,M'\smallsetminus\hs U\,$ with
$\,\vp(x_k)\to\vp_0$ as $\,k\to\infty$, a subsequence of which would have a
limit that lies in $\,P[\vp_0]$, yet not in the set $\,\,U\,$ containing
$\,P[\vp_0]$. However, $\,P[\vp\hs']\,$ obviously intersects {\it each\/} of
the sets $\,\,U_1,\dots,\hs U_l$ whenever $\,\vp\hs'\in(\vp_-,\vp_+)\,$ is
sufficiently close to $\,\vp_0$. Since such $\,P[\vp\hs']\,$ are connected
(see above) and $\,\,U_1,\dots,\hs U_l$ are pairwise disjoint and open, we
must have $\,l=1$. This completes the proof.{\hfill\qd}
\enddemo

\head\S\sm. Special K\"ahler-Ricci potentials\endhead
Except for Remarks \a\sm.3 -- \a\sm.4, the material in this section also
appears in \cite{\dml}.

Given a special \krp\ $\,\vp:M\to\bbR\,$ on a K\"ahler manifold $\,(M,g)$, as
in \f{\ir.1}, let $\,M'\subset M\,$ be the open set on which $\,d\vp\ne0$, and
let the vector fields $\,v,u\,$ on $\,M\,$ and distributions
$\,\Cal H,\Cal V\hs$ on $\,M'$ be defined by
$$\Cal V\,=\,\,\text{\rm Span}\,\{v,u\}\hskip9pt\text{\rm and}\hskip9pt
\Cal H\,=\,\Cal V^\perp\,,\hskip6pt\text{\rm with}\hskip6pt v\,=\,\navp
\hskip6pt\text{\rm and}\hskip6pt u\,=\,Jv\,.\ff\sm.1$$
Setting $\,Q=g(\navp,\navp)$, we thus have
$$g(v,v)\,=\,g(u,u)\,=\,Q\,,\qquad g(v,u)\,=\,0\qquad\text{\rm everywhere\ in}
\quad M.\ff\sm.2$$
In view of the eigenvector clause of \f{\ir.1}, and since \f{\sm.6.b} implies
\f{\sm.6.a} (see Remark \a\sm.1 below), there exist
$\,C^\infty$ functions $\,\si,\ta,\la,\my:M'\to\bbR\,$ with
$$\alignedat2
&\ri\,\,=\,\la\hs g\quad\text{\rm on}\quad\Cal H\hs,\qquad
&&\ri\,\,=\,\my\hs g\quad\text{\rm on}\quad\Cal V\hs,\\
&\nabla d\vp\,=\,\si g\quad\text{\rm on}\quad\Cal H\hs,\qquad
&&\nabla d\vp\,=\,\ta g\quad\text{\rm on}\quad\Cal V\hs,\hskip45pt\\
&\ri\hskip1pt(\Cal H,\Cal V)\,=\,(\nabla d\vp)(\Cal H,\Cal V)\,=\,&&\{0\}
\qquad\text{\rm for}\quad\Cal H,\Cal V\quad\text{\rm as\ in}\quad
\text{\rm\f{\sm.1}.}\endalignedat\ff\sm.3$$
The last line states that $\,\Cal H\,$ is both $\,\,\ri$-or\-thog\-o\-nal and
$\,\nabla d\vp$-or\-thog\-o\-nal to $\,\Cal V\nh$.  We set $\,\si=\la=0\,$
if $\,\dimc M=1$, For a vector field $\,w\,$ on $\,M'\nh$, \f{\pn.3} now gives
$$\text{\rm$\nabla_{\!w}v\,$ equals $\,\si\hskip.4ptw\,$ (or,
$\,\ta\hskip.4ptw$) whenever $\,w\,$ is a section of $\,\Cal H\,$ (or, of
$\,\Cal V$).}\ff\sm.4$$
Also, according to \cite{\dml}, Lemmas \a\sr.5 and \a\qs.1(b), on $\,M'$ we
have
$$\alignedat2
&\text{\rm\ptmi i)}\quad&&
dQ\,=\,2\ta\,d\vp\,,\hskip7pt\text{\rm i.e.,}\hskip6pt\nabla Q\,
=\,2\ta\hskip.4ptv\,,\hskip6pt\text{\rm and}\hskip6pt
\nabla\nh\si\,=\,2(\ta-\si)\hs\si\hskip.4ptv/Q\,.\hskip13pt\\
&\text{\rm ii)}\quad&&
Y\,=\,2\ta\,+\,2(m-1)\hs\si\,,\hskip7pt\text{\rm where}\hskip6ptY\,
=\,\Delta\vp\,.\endalignedat\ff\sm.5$$
Relations $\,dQ=2\ta\,d\vp\,$ and \f{\sm.5.ii} also follow from \f{\sm.3}, as
$\,d_wQ=d_w[g(v,v)]=2\hs g(\nabla_{\!w}v,v)=2\hs\ta g(w,v)=2\hs\ta\hs d_w\vp$,
by \f{\sm.2}, \f{\sm.4}, \f{\sm.1}, for any vector field $\,w$.

As $\,dY=-\hs2\,\ri\hs(\navp,\,\cdot\,)\,$ for any \kip\ $\,\vp$, the Ricci
tensor $\,\,\ri\hs$, and $\,Y=\Delta\vp=\langle\hs g,\nabla d\vp\rangle\,$
(\cite{\cal}; see also \cite{\dml}, formula \f{\kp.4}), \f{\sm.1} and
\f{\sm.3} give $\,dY=-\hs2\my\,d\vp$.
\rmk{Remark \a\sm.1}For a distribution $\,\Cal V\hs$ on a Riemannian
manifold $\,(M,g)\,$ and a symmetric twice-co\-var\-i\-ant tensor $\,\bz\,$
at a point $\,x\in M$, consider the conditions
\vskip4pt
\hbox{\hskip-.8pt
\vbox{\hbox{\f{\sm.6}}\vskip3.2pt}
\hskip16pt
\vbox{
\hbox{a)\hskip9ptAll nonzero vectors in $\,\Cal V_x$ and
$\,\Cal H_x=\Cal V_x^\perp$ are eigenvectors of $\,\bz$.}
\vskip1pt
\hbox{b)\hskip9ptAll nonzero vectors in $\,\Cal H_x=\Cal V_x^\perp$ are
eigenvectors of $\,\bz$.}}}
\vskip4pt
\noindent Let $\,\,\Cal V\hs$ now be a $\,J${\it-invariant\/} distribution of
complex dimension one on a K\"ahler manifold $\,(M,g)$, and let $\,x\in M$.
For a symmetric twice-co\-var\-i\-ant tensor $\,\bz\,$ at $\,x\,$ which is
also Hermitian (cf. Lemma \a\kg.2(iii)), condition \f{\sm.6.b} then implies
\f{\sm.6.a}. In fact, the operator $\,B:T_xM\to T_xM\,$ with
$\,\bz(w,w')=g(Bw,w')\,$ for all $\,w,w'\in T_xM\,$ is self-adjoint, commutes
with $\,J$, and $\,B\Cal V_x^\perp\subset\Cal V_x^\perp$. Hence
$\,B\Cal V_x\subset\Cal V_x$. Choosing $\,v\in\Cal V_x\smallsetminus\{0\}\,$
and $\,\la\in\bbR\,$ with $\,Bv=\la v$, we thus have $\,Bu=\la u\,$ for
$\,u=Jv\,$ (as $\,BJv=JBv$), which yields \f{\sm.6.a} since
$\,\,\dimr\Cal V_x=2$.
\endrmk
\proclaim{Lemma \a\sm.2}Given\/ $\,\vp\,$ with\/ \f{\ir.1} on a K\"ahler
manifold\/ $\,(M,g)$, let\/ $\,Q=g(\navp,\navp)\,$ and let\/ $\,\si\,$ be as
in\/ \f{\sm.3}. Then either\/ $\,\si=0\,$ identically on\/ $\,M'$, or\/
$\,\si\ne0\,$ everywhere in\/ $\,M'\nh$, where\/ $\,M'\subset M\,$ is the open
set on which\/ $\,d\vp\ne0$. In the latter case, there exists a constant\/
$\,\,\y\,$ with\/ $\,Q/\si=2(\vp-\y)\,$ and\/ $\,\vp\ne\y\,$ everywhere in\/
$\,M'$.
\endproclaim
For a proof, see \cite{\dml}, Lemma \a\cm.5. Note that relation $\,\vp\ne\y\,$
on $\,M'$ is obvious from $\,Q/\si=2(\vp-\y)$.{\hfill\qd}
\rmk{Remark \a\sm.3}Let a function $\,\vp\,$ satisfy \f{\ir.1} on a
K\"ahler manifold $\,(M,g)$. We define $\,\ve\in\{-\hs1,0,1\}\,$ by
$\,\ve=0\,$ when $\,\si=0\,$ identically on $\,M'$ and
$\,\ve=\,\text{\rm sgn}\,(\vp-\y)\,$ when $\,\si\ne0\,$ everywhere in $\,M'$
(with $\,M'\nh,\si,\y\,$ as in Lemma \a\sm.2). Note that in the latter case
$\,\ve=\pm\hs1\,$ is uniquely defined, since $\,\vp\ne\y\,$ everywhere in
$\,M'$ (Lemma \a\sm.2) and $\,M'$ is connected (Remark \a\cz.3(i)).
\endrmk
\rmk{Remark \a\sm.4}Let $\,Q=g(\navp,\navp)\,$ and let $\,\si,\ta\,$ be as
in \f{\sm.3} for a special \krp\ $\,\vp\,$ on a K\"ahler manifold $\,(M,g)\,$
(see \f{\ir.1}), and let $\,v=\navp$.
\widestnumber\item{(ii)}\roster
\item"(i)"Writing $\,\dot f=\hs d\hs[f(x(\sa))]/d\sa\,$ for a fixed $\,C^1$
curve $\,\sa\mapsto x(\sa)\in M\,$ and a $\,C^1$ function $\,f\,$ defined in
$\,M\nh$, we have $\,\dot f=d_{\dot x}f=g(\nabla\!f,\dot x)$, where
$\,\dot x=\hs dx/d\sa$. Consequently, $\,g(v,\dot x)=\dot\vp\,$ and, by
\f{\sm.5.i}, $\,\dot Q=2\ta\dot\vp$.
\item"(ii)"Given $\,y\in M\,$ with $\,v(y)=0$, let
$\,\sa\mapsto x(\sa)\in M\,$ be a $\,C^2$ curve such that $\,x(0)=y\,$ and
$\,\dot x(0)\,$ is an eigenvector of $\,(\nabla d\vp)(y)\,$ for an eigenvalue
$\,a\ne0$. (Such $\,a$ exists since $\,\nabla d\vp\ne0\,$ at $\,y\,$
by \f{\ir.1} and Remark \a\kg.4.) Then $\,\dot\vp=0\,$ and
$\,\ddot\vp=a|\dot x|^2\ne0\,$ at $\,\sa=0\,$ (notation of (i)), so that
$\,\dot\vp\ne0\,$ for all $\,\sa\ne0\,$ close to $\,0$. In fact,
$\,g(v,\dot x)=\dot\vp\,$ (see (i)) and applying $\,d/d\sa\,$ we obtain
$\,\ddot\vp=g(\nabla_{\!\dot x}v,\dot x)+g(v,\nabla_{\!\dot x}\dot x)$, which
at $\,\sa=0\,$ equals $\,(\nabla d\vp)(\dot x,\dot x)\,$ (by \f{\pn.2} with
$\,v(y)=0$).
\endroster
\endrmk
\rmk{Remark \a\sm.5}Let $\,Q=g(\navp,\navp)\,$ for a function $\,\vp\,$
with \f{\ir.1} on a K\"ahler manifold $\,(M,g)$, and let $\,M'\subset M\,$ be
the open set given by $\,d\vp\ne0$. Then
\widestnumber\item{(b)}\roster
\item"(a)"Every point of $\,M'$ has a neighborhood on which $\,Q\,$ is a
$\,C^\infty$ function of $\,\vp\,$ such that $\,dQ/d\vp=2\ta$, with $\,\ta\,$
given by \f{\sm.3}.
\item"(b)"$Q\,$ is constant along every connected component of the
pre\-im\-age of any real number under $\,\vp:M'\to\bbR\hs$.
\endroster
In fact, (a) is obvious as $\,dQ=2\ta\,d\vp\,$ (cf. \f{\sm.5.i}); using local
coordinates in $\,M'$ having $\,\vp\,$ as one of the coordinate functions, one
sees that the assignment $\,\vp\mapsto Q\,$ is (locally, in $\,M'$) of class
$\,C^\infty$. Thus, $\,Q\,$ is locally constant on each of the connected
components in (b), and assertion (b) follows.
\endrmk
\head\S\xa. The simplest examples\endhead
\example{Example \a\xa.1}A nonconstant $\,C^\infty$ function $\,\vp\,$ on a
K\"ahler manifold $\,(M,g)\,$ of complex dimension $\,1\,$ satisfies \f{\ir.1}
if and only if it is a \kip, as the remainder of \f{\ir.1} then is vacuously
true. Consequently, in complex dimension $\,1\,$ there is, locally, a
one-to-one correspondence between special \krp s (defined up to an additive
constant), and nontrivial Killing fields. This is clear from Lemma \a\kg.3 and
the fact that, in complex dimension $\,1$, every Killing field $\,u\,$ is
holomorphic (by \f{\kg.1.a}, since skew-adjoint\-ness of $\,\nabla u$, cf.
\S\kg, then gives $\,\nabla u=\ta J\,$ for some function $\,\ta$).
\endexample
Similarly, for a nonconstant $\,C^\infty$ function $\,\vp\,$ on a K\"ahler
manifold $\,(M,g)\,$ of any complex dimension $\,m$, one has \f{\ir.1}
whenever $\,g\,$ is an Einstein metric and $\,\nabla d\vp=\ta g\,$ for some
function $\,\ta$. For instance, this is the case when $\,g\,$ is the standard
Euclidean metric on $\,M=\bbC^m$, with the norm-squared function $\,\vp\,$ and
$\,\ta=2$, or with a real-linear function $\,\vp\,$ and $\,\ta=0$. (Cf.
Example \a\xa.4 below.)
\example{Example \a\xa.2}Let $\,(\srf,\gx)\,$ be an oriented Riemannian
surface admitting a nonconstant \kip\ $\,\vp:\srf\to\bbR\hs$, and let
$\,(N,h)\,$ be a K\"ahler manifold of complex dimension $\,m-1\ge1\,$ with the
Ricci tensor $\,\,\rih=\la\hs h\,$ for some function $\,\la:N\to\bbR\hs$, so
that $\,(N,h)\,$ is K\"ahler-Einstein (and $\,\la\,$ is constant) unless
$\,m=2$. Note that, for an oriented surface $\,(\srf,\gx)\,$ with
$\,H^1(\srf,\bbR)=\{0\}$, such $\,\vp\,$ exists if and only if
$\,(\srf,\gx)\,$ admits a nontrivial Killing field (Lemma \a\kg.3 and Example
\a\xa.1).

Treated as a function on $\,M=N\nh\times\srf\,$ constant along the $\,N\hs$
factor, $\,\vp\,$ then is a special \krp\ on the K\"ahler manifold $\,(M,g)\,$
of complex dimension $\,m\ge2\,$ obtained as the Riemannian product of
$\,(N,h)\,$ and $\,(\srf,\gx)$.

In fact, let $\,\Cal H,\Cal V\hs$ be the $\,N\hs$ and $\,\srf\,$ factor
distributions on $\,M$. Conditions \f{\sm.3} are obviously satisfied by the
function $\,\la\,$ with $\,\,\rih=\la\hs h\,$ and the Gaussian curvature
$\,\my\,$ of $\,\gx$, along with $\,\si=0\,$ and $\,\ta=\Delta\vp/2\,$ (the
Laplacian of $\,\vp/2\,$ in $\,(\srf,\gx)$), as relation $\,\nabla u=\ta J\,$
in Example \a\xa.1 gives $\,\nabla d\vp=\ta\gx$, by \f{\pn.2} with $\,g\,$
replaced by $\,\gx$. This implies \f{\ir.1}, as $\,\vp\,$ is a \kip\ on
$\,(M,g)\,$ due to being one when viewed as a function on $\,(\srf,\gx)$.
\endexample
\proclaim{Lemma \a\xa.3}Let two vector fields\/ $\,v,u\,$ on a Riemannian
manifold\/ $\,(M,g)\,$ be linearly independent at every point, and let\/
$\,\Cal H\,$ be a distribution on\/ $\,M\,$ with\/
$\,TM=\Cal H\oplus\Cal V$, where\/ $\,\Cal V=\,\text{\rm Span}\,\{v,u\}$. If\/
$\,G\,$ is a group of isometries of\/ $\,(M,g)\,$ such that, at every\/
$\,x\in M$, the action on\/ $\,T_xM$, via differentials, of the isotropy
subgroup of\/ $\,G\,$ at $\,x\,$ leaves\/ $\,v(x)$, $\,u(x)\,$ and\/
$\,\Cal H_x$ invariant and acts transitively on the unit sphere in\/
$\,\Cal H_x$, then\/ $\,\Cal H=\Cal V^\perp$, and any\/ $\,G$-invariant
symmetric twice-co\-var\-i\-ant tensor field\/ $\,\bz\,$ on\/ $\,M\,$
satisfies condition\/ \f{\sm.6.b} at every point\/ $\,x\in M$.
\endproclaim
In fact, for such $\,\bz\,$ and $\,x\,$ the functions $\,\bz(v(x),w)$,
$\,\bz(u(x),w)\,$ and $\,\bz(w,w)\,$ of $\,w\in\Cal H_x$ are constant on the
unit sphere; the first two are also linear, so they must be identically zero.
This yields \f{\sm.6.b} and, applied to $\,\bz=g(x)$, gives
$\,\Cal H=\Cal V^\perp$.{\hfill\qd}
\example{Example \a\xa.4}Let our complex manifold be any $\,G$-invariant
nonempty open connected set $\,\,U\,$ in a complex vector space $\,V$ of
complex dimension $\,m\ge2$, carrying a fixed Hermitian inner product
$\,\langle\,,\rangle$, where $\,G\approx\hs\text{\rm U}\hs(m)\,$ is the group
of all unitary automorphisms $\,V\to V\nh$. For any $\,G$-invariant
K\"ahler metric $\,g\,$ on $\,\,U\nh$, a special \krp\ $\,\vp\,$ on
$\,(\hs U,g)\,$ (cf. \f{\ir.1}) can be obtained as follows.

Formula $\,u(x)=aix$, with any fixed real constant $\,a\ne0$, defines a
$\,G$-invariant holomorphic Killing field $\,u\,$ on $\,(\hs U,g)$. (Thus,
$\,u\,$ is an infinitesimal generator of the center subgroup of $\,G$.) We now
choose $\,\vp:M\to\bbR\,$ to be a \kip\ with $\,u=Jv\,$ for $\,v=\navp$, where
$\,\nabla\,$ is the $\,g$-gra\-di\-ent; by Lemma \a\kg.3, such $\,\vp\,$
exists and is unique up to an additive constant.

Since our $\,v,u\,$ are the same as in Remark \a\tb.1, Lemma \a\xa.3 now shows
that $\,\Cal V,\Cal H\,$ defined in Remark \a\tb.1 are $\,g$-or\-tho\-gon\-al
to each other, and the tensors $\,\bz=\,\ri\,$ and $\,\bz=\nabla d\vp\,$
satisfy \f{\sm.6.b} at every $\,x\in\,U\smallsetminus\{0\}$, which yields
\f{\ir.1}.
\endexample
\example{Example \a\xa.5}For an integer $\,m\ge2\,$ and a fixed point
$\,y\in\bbCP^m$, let $\,g\,$ be any $\,G$-invariant K\"ahler metric on
$\,\bbCP^m$, where $\,G\,$ is the
group of all biholomorphisms $\,\bbCP^m\to\bbCP^m$ that keep $\,y\,$ fixed and
preserve the standard (Fubini-Study) metric. Then $\,(\bbCP^m,g)\,$ admits a
function $\,\vp\,$ with \f{\ir.1}, obtained as follows.

As $\,G\approx\hs\text{\rm U}\hs(m)\,$ (see below), its center, isomorphic to
$\,\,\text{\rm U}\hs(1)$, is generated by a nontrivial holomorphic Killing
field $\,u\,$ on $\,(\bbCP^m,g)$, unique up to a factor; we choose
a \kip\ $\,\vp:\bbCP^m\to\bbR\,$ with $\,u=J(\navp)$, where $\,\nabla\,$ is
the $\,g$-gra\-di\-ent. Such $\,\vp\,$ exists and is unique up to an
additive constant (Lemma \a\kg.3).

The usual identification of $\,\bbC^m$ with an open dense set in $\,\bbCP^m$,
such that $\,y=0\in\bbC^m\subset\bbCP^m\nh$, makes $\,G\,$ act
on $\,\bbC^m$ as the matrix group $\,\,\text{\rm U}\hs(m)$. Restricting
$\,g\,$ to $\,\bbC^m$ we then obtain a $\,\,\text{\rm U}\hs(m)$-invariant
K\"ahler metric on $\,\bbC^m$. That $\,\vp\,$ is a special \krp\ on
$\,(\bbCP^m\nh,g)\,$ is now immediate from Example \a\xa.4: since $\,\bbC^m$
is dense in $\,\bbCP^m$, the eigenvector clause of \f{\ir.1} holds on
$\,\bbCP^m$ as well.
\endexample
\rmk{Remark \a\xa.6}More precisely, on
$\,\,\bbCP^{m-1}=\bbCP^m\smallsetminus\bbC^m$ (the hyperplane at infinity) the
eigenvector clause of \f{\ir.1} holds vacuously, since $\,u\,$ (and hence
$\,d\vp$) vanishes there; namely, the flow of $\,u\,$ in $\,\bbC^m$ consists
of multiples of the identity, all of which leave every line through $\,0\,$
invariant, i.e., keep every point at infinity fixed. Since $\,0\,$ is the only
point they keep fixed in $\,\bbC^m\nh$, it follows that the set of critical
points of $\,\vp$, i.e., zeros of $\,u$, is the union of two critical
manifolds (cf. Remark \a\cz.3(ii)): the one-point set $\,\{y\}\,$ and the
hyperplane at infinity.
\endrmk

\head\S\am. A local model\endhead
All special \krp s $\,\vp\,$ on K\"ahler manifolds arise (locally, at points
with $\,d\vp\ne0$) from the construction described below; see also \S\tl.

Given a positive $\,C^\infty$ function $\,Q\,$ of a real variable $\,\vp$,
defined on an open interval $\,\iyp\nh$, and a real constant $\,a\ne0$, a
$\,C^\infty$ function $\,r\hs$ of the variable $\,\vp\in\iyp$ with
$$dr/d\vp\,=\,ar/Q\hskip10pt\text{\rm and}\hskip10ptr\,>\,0\hskip10pt\text{\rm
on}\hskip8pt\iyp\ff\am.1$$
exists and is unique up to a constant factor, constituting a diffeomorphism
$$\iyp\ni\vp\,\mapsto\,r\in(r_-,r_+)\quad\text{\rm with}\quad
0\,\le\,r_-<\,r_+\le\,\infty\hs.\ff\am.2$$
Let there be given the following set of data:
$$\iyp\!,\vp,Q,r\hs;\qquad a,\ve,\y\hs;\qquad m,N,h\hs;\qquad
\Cal L,\Cal H,\langle\,,\rangle\hs.\hskip10pt\ff\am.3$$
Here $\,\iyp\subset\bbR\,$ is an open interval, $\,\vp\in\iyp$ is a real
variable, $\,Q\,$ is a positive $\,C^\infty$ function $\,Q\,$ of the variable
$\,\vp\in\iyp$, while $\,a\ne0\,$ is a real constant, $\,r\hs$ is a fixed
function of $\,\vp\in\iyp$ satisfying \f{\am.1}, and $\,\ve,\y\,$ are
constants such that either $\,\ve=0$, or $\,\y\notin\iyp$ and
$\,\ve=\,\text{\rm sgn}\,(\vp-\y)=\pm\hs1\,$ for all $\,\vp\in\iyp$. (When
$\,\ve=0$, we leave $\,\y\,$ undefined.) Next, $\,m\ge2,$ is an integer,
$\,(N,h)\,$ is a K\"ahler manifold of complex dimension $\,m-1\,$ (which we
assume to be Einstein unless $\,m=2$), while $\,\Cal L\,$ is a $\,C^\infty$
complex line bundle over $\,N\nh$, with a Hermitian fibre metric
$\,\langle\,,\rangle$, and $\,\Cal H\,$ is the horizontal distribution of a
connection in $\,\Cal L\,$ making $\,\langle\,,\rangle\,$ parallel, whose
curvature form (Remark \a\cn.1) equals $\,-\hs2\ve a\,$ times the K\"ahler
form of $\,(N,h)\,$ (cf. \f{\pn.5}).

Using the diffeomorphism \f{\am.2}, we treat functions of $\,\vp\in\iyp$ as
functions of $\,r\in(r_-,r_+)$. This includes $\,\vp\,$ itself, $\,Q\,$ and
$\,\fe\,$ defined by
$$\text{\rm$\fe=1\,$ \ (when $\,\ve=0$), \ or \ $\,\fe=2\hs|\vp-\y\hs|\,$ \
(when $\,\ve=\pm\hs1$).}\hskip5pt\ff\am.4$$
Let $\,\,U\,$ now be the open subset of $\,\Cal L\smallsetminus N\,$ given by
$\,r_-<r<r_+$, where $\,r\hs$ this time denotes the norm function of
$\,\langle\,,\rangle\,$ (see Remark \a\cn.2); thus, $\,r,\vp,Q,\fe\,$ can also
be regarded as $\,C^\infty$ functions $\,\,U\to\bbR\hs$. We define a metric
$\,g\,$ on $\,\,U\,$ by
$$g=\fe\hs\proj^*\nh h\hskip5pt\text{\rm on}\hskip5pt\Cal H\hs,\hskip12pt
g=(ar)^{-2}Q\,\text{\rm Re}\hskip1pt\langle\,,\rangle\hskip5pt\text{\rm on}
\hskip5pt\Cal V\hs,\hskip12ptg(\Cal H,\Cal V)=\{0\}\hs,\ff\am.5$$
where $\,\proj:\Cal L\to N\,$ is the bundle projection, $\,\Cal V\hs$ denotes
the vertical distribution in $\,\Cal L\hs$, and
$\,\,\text{\rm Re}\hskip1pt\langle\,,\rangle\,$ is the standard Euclidean
metric on each fibre of $\,\Cal L\hs$. (The last relation in \f{\am.5} means
that $\,\Cal H\,$ is $\,g$-or\-thog\-o\-nal to $\,\Cal V\nh$.)

Moreover (see Remark \a\cn.4), $\,\Cal L\,$ has a unique structure of a
holomorphic line bundle over $\,N\hs$ such that $\,\Cal H\,$ is
$\,J$-invariant. This turns the open submanifold $\,\,U\,$ of $\,\Cal L\,$
into a complex manifold of complex dimension $\,m$. According to \cite{\dml},
beginning of \S\dx\ and Remark \a\dx.1, for $\,\si,\ta,\la,\my,\fe\,$ as in
\f{\sm.3} and \f{\am.4},
\widestnumber\item{(e)}\roster
\item"(a)"$\,g\,$ is a K\"ahler metric on $\,\,U\nh$,
\item"(b)"$\,\vp\,$ is a special \krp\ on $\,(\hs U\nh,g)$, in the sense of
\f{\ir.1},
\item"(c)"$\,Q\,$ treated as a function $\,\,U\to\bbR\,$ is given by
$\,Q=g(\navp,\navp)$.
\item"(d)"If $\,\ve=0\,$ then $\,\si=0\,$ identically, while if
$\,\ve=\pm\hs1\,$ then $\,\si\ne0\,$ everywhere and $\,\y\,$ in \f{\am.3} is
the same as in Lemma \a\sm.2, i.e., $\,2\si=Q/(\vp-\y)$.
\item"(e)"Let $\,\kx:N\to\bbR\,$ be a function such that $\,h\,$ has the Ricci
tensor $\,\,\rih=\kx\hs h$. Identifying $\,\kx\,$ with the composite
$\,\kx\circ\proj$, we may treat it as a function on $\,\,U\nh$, constant
unless $\,m=2$. Then $\la=(\kx-\ve Y)/\!\fe$, with
$\,Y=\Delta\vp=\langle\hs g,\nabla d\vp\rangle$.
\endroster
\rmk{Remark \a\am.1}Let the data \f{\am.3} be chosen so that $\,h\,$ is an
Einstein metric, $\,\ve=\pm\hs1$, and $\,Q\,$ satisfies the differential
equation $\,p\hs Q'-Q+(m-1)p\hs Q/(\vp-\y)=\ve\hs p\kx-2\fy(\vp-\y)$, where
$\,Q'=\hs dQ/d\vp\,$ and $\,\kx\,$ is the constant in (e), while $\,p,\fy\,$
are constants with $\,p\ne0$. Also, let $\,\nabla,\,\ri\,\,$ be the
Levi-Civita connection and Ricci tensor of the metric $\,g\,$ with \f{\am.5}
on the manifold $\,\,U\,$ constructed above.

These $\,(\hs U\nh,g)\,$ and $\,\vp:U\to\bbR\,$ then satisfy the
condition $\,\nabla d\vp+p\,\ri\,=\fy g$.

In fact, $\,2\ta=\hs dQ/d\vp=Q'$ and $\,2\si=Q/(\vp-\y)\,$ by \f{\sm.5.i} and
(d), so that \f{\sm.5.ii} gives $\,Y=Q'+(m-1)Q/(\vp-\y)$, i.e., our assumption
yields $\,pY-Q=\ve\hs p\kx-2\fy(\vp-\y)$. Applying $\,d/d\vp\,$ we get
$\,-\hs2(\ta+p\my)=-\hs2\fy$, as $\,dY/d\vp=-\hs2\my\,$ (see the line
preceding Remark \a\sm.1) and $\,Q'=2\ta$, while
$\,-\hs2(\si+p\la)=(pY-Q-\ve\hs p\kx)/(\vp-\y)=-\hs2\fy\,$ since
$\,2\si=Q/(\vp-\y)\,$ and $\,2\la=(\ve\kx-Y)/(\vp-\y)$, by (e), as
$\,\fe=2\ve(\vp-\y)\,$ (cf. \f{\am.4}, with
$\,\ve=\,\text{\rm sgn}\,(\vp-\y)$). Our claim now follows from \f{\sm.3} with
$\,\ta+p\my=\fy\,$ and $\,\si+p\la=\fy$.

Note that, whenever such a triple $\,\,U\nh,g,\vp\,$ can be ``compactified''
as in \S\mw, it becomes an example of a K\"ahler-Ricci soliton (\cite{\ptv},
\cite{\tiz}, \cite{\cao}). See \cite{\dms} for details.
\endrmk

\head\S\tm. The case of open spherical shells\endhead
Suppose that we are given the data
$$\iyp\!,\vp,Q,r\hs;\qquad a,\ve,\y\hs;\qquad m,V,\langle\,,\rangle\hskip10pt
\ff\tm.1$$
which consist of an open interval $\,\iyp\subset\bbR\hs$, a positive
$\,C^\infty$ function $\,Q\,$ of the variable $\,\vp\in\iyp$, real constants
$\,a,\y\,$ and $\,\ve=\pm\hs1$ such that $\,\ve a>0\,$ and $\,\ve(\vp-\y)>0\,$
for all $\,\vp\in\iyp\nh$, a fixed function $\,r\hs$ of $\,\vp\in\iyp$
satisfying \f{\am.1}, a complex vector space $\,V$ of complex
dimension $\,m\ge2$, and a Hermitian inner product $\,\langle\,,\rangle\,$ in
$\,V\nh$. The symbol $\,r\hs$ also denotes the norm function
$\,V\to[\hs0,\infty)$, with $\,\zx\mapsto\langle\zx,\zx\rangle^{1/2}$.

Let $\,\,U\,$ be the open spherical shell in $\,V$ lying between spheres of
radii $\,r_+,r_-$ centered at $\,0$, i.e., given by $\,r_-<r<r_+$, with
$\,r_\pm$ as in \f{\am.2}. Our $\,\vp\,$ and $\,Q$, being $\,C^\infty$
functions of $\,r\in(r_-,r_+)$, thus become
$\,C^\infty$ functions $\,\,U\to\bbR\hs$. We now define a Riemannian metric
$\,g\,$ on $\,\,U\,$ by
$$|a|r^2g=2\hs|\vp-\y\hs|\,\text{\rm Re}\hskip1pt\langle\,,\rangle
\hskip5pt\text{\rm on}\hskip5pt\Cal H,\hskip8pt
a^2r^2g=Q\,\text{\rm Re}\hskip1pt\langle\,,\rangle\hskip5pt
\text{\rm on}\hskip5pt\Cal V,\hskip9ptg(\Cal H,\Cal V)=\{0\}\ff\tm.2$$
for $\,\Cal H,\Cal V\,$ as in Remark \a\tb.1,
$\,\,\text{\rm Re}\hskip1pt\langle\,,\rangle\,$ being the standard Euclidean
metric on $\,V\nh$.
\proclaim{Lemma \a\tm.1}The above construction of\/ $\,\,U,g\,$ and\/
$\,\vp:\hs U\to\bbR\,$ is a special case of that in\/ {\rm\S\am}. Namely,
our\/ $\,\,U,g,\vp\,$ are the same as those obtained from the data\/ \f{\am.3}
with\/ $\iyp\!,\vp,Q,r,a,\ve,\y\hs,m\,$ as in\/ \f{\tm.1},
$\,N,\Cal L,\Cal H,\langle\,,\rangle\,$ defined as in\/
{\rm\S\tb}\hskip4ptfor\/ $\,m,V,\langle\,,\rangle\,$ in\/ \f{\tm.1}, and\/
$\,h\,$ equal to\/ $\,1/|a|\,$ times the Fubini-Study metric on the projective
space\/ $\,N\nh$. The latter data\/ \f{\am.3} satisfy all the conditions
listed in\/ the paragraph following\/ \f{\am.3}, while equality between\/
$\,\,U\subset V\smallsetminus\{0\}\,$ and\/
$\,\,U\subset\Cal L\smallsetminus N\,$ makes sense due to the biholomorphic
identification\/ $\,V\smallsetminus\{0\}=\Cal L\smallsetminus N\,$ of\/
{\rm\S\tb}.

In particular, $\,g\,$ is a K\"ahler metric on\/ $\,\,U\,$ with the complex
structure of an open submanifold of\/ $\,V\nh$, and\/ $\,\vp\,$ is a special
\krp\ on\/ $\,(\hs U,g)$.
\endproclaim
In fact, in \S\tb\ we verified that the data \f{\am.3} described in the lemma
satisfy the assumptions made in \S\am. (Note that $\,\varOmega=\,-\hs\omfs$
equals $\,-\hs2\ve a\,$ times the K\"ahler form of $\,(N,h)$, as
$\,\ve a=|a|$.) The last two relations of \f{\am.5}, for our $\,g$, are
immediate from \f{\tm.2}. Finally, the first equality in \f{\am.5} follows
from \f{\am.4} and \f{\tm.2} with
$\,|a|r^2\hs\proj^*\nh h=\,\text{\rm Re}\hskip1pt\langle\,,\rangle\,$ on
$\,\Cal H\,$ (see Remark \a\tb.2) and with the identification
$\,\Cal L\smallsetminus N=V\smallsetminus\{0\}$, which also proves our claim
about the complex structure.{\hfill\qd}
\medskip
The triples $\,\,U,g,\vp\,$ constructed here also form a special case of those
described in Example \a\xa.4. In fact,  $\,g\,$ is a K\"ahler metric,
$\,u=J(\navp)\,$ is a holomorphic Killing field (Lemmas \a\tm.1 and \a\kg.3),
and $\,g,\vp\,$ are both invariant under the unitary group
$\,G\approx\hs\text{\rm U}\hs(m)$, since so are the Euclidean metric
$\,\,\text{\rm Re}\hskip1pt\langle\,,\rangle$, its norm function $\,r$, and
$\,\vp\,$ (which is a function of $\,r$). Since $\,u\,$ is $\,G$-invariant, it
generates the center of $\,G$.

\head\S\mm. Duality and metrics on annuli\endhead
Let $\,V^*$ be the dual space of a complex vector space $\,V$ of complex
dimension $\,1$. We define the {\it inversion biholomorphism\/}
$\,V\smallsetminus\{0\}\to V^*\smallsetminus\{0\}\,$ to be the assignment
$\,\zx\mapsto\zx^{-1}$, where $\,\zx^{-1}\in V^*$ is the $\,\bbC$-linear
functional $\,V\to\bbC\,$ sending $\,\zx\,$ to $\,1$. Any fixed Hermitian
inner product $\,\langle\,,\rangle\,$ in $\,V$ gives rise to a Hermitian
inner product $\,\langle\,,\rangle\nh^*$ in $\,V^*$ such that
$\,\langle\zx^{-1},\zx^{-1}\rangle\nh^*=\langle\zx,\zx\rangle^{-1}$ whenever
$\,\zx\in V\smallsetminus\{0\}$.

Given a septuple $\,\iyp\!,\vp,Q,r,a,V,\langle\,,\rangle\,$ consisting of
\vskip5pt
\settabs\+\noindent&\f{\mm.1}\hskip7pt&\cr
\+&&an open interval $\hs\iyp\subset\bbR\hs$, a positive $\hs C^\infty$
function $\hs\iyp\ni\vp\mapsto Q$, a real con-\cr
\+&\f{\mm.1}&stant $\,a\ne0$, a function $\,\hs r\hs\hs$ of
$\,\hs\vp\in\iyp$ satisfying \f{\am.1}, and a complex vec-\cr
\+&&tor space $\,V\nh$ of complex dimension $1$ with a Hermitian inner product
$\hs\langle\,,\rangle$,\cr
\vskip4pt
\noindent
let $\,\,U\,$ be the open annulus in $\,V$ given by $\,r_-<r<r_+$. Here
$\,r\hs$ also denotes the norm function of $\,\langle\,,\rangle$, which allows
us to treat $\,r$, $\,\vp\,$ and $\,Q\,$ as $\,C^\infty$ functions
$\,\,U\to\bbR\hs$. Formula
$\,\gx=(ar)^{-2}Q\,\text{\rm Re}\hskip1pt\langle\,,\rangle\,$ now defines a
Riemannian metric $\,\gx\,$ on the annulus $\,\,U\nh$, conformal to the
standard Euclidean metric $\,\,\text{\rm Re}\hskip1pt\langle\,,\rangle$.

Next, let us replace $\,r,a,V,\langle\,,\rangle\,$ in these data by
$\,r^*,a^*$ with $\,r^*=1/r$, $\,a^*=-\hs a\,$ and
$\,V^*\nh,\langle\,,\rangle\nh^*$ as in the beginning of this section, but
keep the same $\,\iyp\!,\vp,Q\,$ and $\,m=1$. Since \f{\am.1} implies that
$\,dr^*/d\vp=a^*r^*/Q$, the new data satisfy the same assumptions, and may be
used as above to define a metric $\,\gx^*$ on an open annulus
$\,\,U^*\subset V^*\smallsetminus\{0\}$. Then, with $\,r^*$ also standing for
the norm function of $\,\langle\,,\rangle\nh^*$,
\widestnumber\item{(b)}\roster
\item"(a)"The assignment $\hs\zx\mapsto\zx^{-1}$ is an isometry
$\hs(\hs U,\gx)\to(\hs U^*\nh,\gx^*)$, under which $\,r\hs$ treated as a
function on $\,\,U\,$ corresponds to the function $\,1/r^*$ on $\,\,U^*\nh$.
\item"(b)"$\vp\,$ as a function of $\,r\hs$ with \f{\am.1} is related to
$\,\vp\,$ viewed, similarly, as a function of $\,r^*$, in such a way that the
assignment $\,r\mapsto r^*=1/r\,$ leaves the corresponding value of $\,\vp\,$
unchanged.
\endroster
In fact, regarded as functions of $\,\vp\,$ inverse to those in (b), our
$\,r,r^*$ satisfy $\,r^*=1/r$, which implies (b). To verify (a), let us use
the multiplicative notation $\,\zeta\zx\in\bbC\,$ for evaluating functionals
$\,\zeta\in V^*$ on vectors $\,\zx\in V\nh$. Thus,
$\,|\zeta\zx|^2=\langle\zx,\zx\rangle\hs\langle\zeta,\zeta\rangle\nh^*$. (To
see this, assume that $\,\zx\ne0\,$ and write $\,\zeta\,$ as a scalar times
$\,\zx^{-1}$, cf. the beginning of this section.) Setting
$\,\zeta(\sa)=[\zx(\sa)]^{-1}$ for any $\,C^1$ curve
$\,\sa\mapsto\zx(\sa)\in\,U$, and differentiating the resulting relation
$\,\zeta(\sa)\zx(\sa)=1$, we now see that the differential of
$\,\zx\mapsto\zx^{-1}$ at any point $\,\zx\in\,U\,$ sends each tangent vector
$\,\dot\zx\in T_\zx\hs U=V\,$ to the vector
$\,\dot\zeta=-\hs(\zeta\dot\zx)\hs\zeta\in V^*$ tangent to $\,\,U^*$ at the
point $\,\zeta=\zx^{-1}$. As
$\,\langle\zeta,\zeta\rangle\nh^*=\langle\zx,\zx\rangle^{-1}$ (see above), we
thus have $\,\langle\dot\zeta,\dot\zeta\rangle\nh^*
=\langle\dot\zx,\dot\zx\rangle/\langle\zx,\zx\rangle^2
=\langle\dot\zx,\dot\zx\rangle/r^4$, i.e., the pull\-back under
$\,\zx\mapsto\zx^{-1}$ of the Euclidean metric
$\,\,\text{\rm Re}\hskip1pt\langle\,,\rangle\nh^*$ on
$\,\,U^*$ is $\,1/r^4$ times the Euclidean metric on $\,\,U$. The pull\-back
of the function $\,r^*$ on $\,\,U^*$, that is, the composite of
$\,\zx\mapsto\zx^{-1}$ followed by $\,r^*$, clearly is $\,1/r\,$ on
$\,\,U$, as claimed in (a). (See the definition of
$\,\langle\,,\,\rangle\nh^*$.) Also, as $\,Q\,$ is a function of $\,\vp$, (b)
remains valid when $\,\vp\,$ is replaced by $\,Q$, so that the pull\-back of
$\,Q\,$ is $\,Q$. Consequently, the pull\-back of
$\,\gx^*=(a^*r^*)^{-2}Q\,\text{\rm Re}\hskip1pt\langle\,,\rangle\nh^*$ is
$\,\gx=(ar)^{-2}Q\,\text{\rm Re}\hskip1pt\langle\,,\rangle$, as required.

\head\S\ib. The inversion biholomorphism\endhead
Let $\,\Cal L^*$ be the dual of a holomorphic line bundle $\,\Cal L\,$ over a
complex manifold $\,N\nh$. As in \f{\cn.2}, the symbols $\,\Cal L,\Cal L^*$
also stand for their total spaces, and $\,N\hs$ is identified with the zero
sections $\,N\subset\Cal L\,$ and $\,N^*\subset\Cal L^*$. We now define
$\,M\,$ to be the complex manifold obtained from the disjoint union
$\,\Cal L\cup\Cal L^*$ by identifying the open subsets
$\,\Cal L\smallsetminus N\,$ and $\,\Cal L^*\smallsetminus N^*$ via the {\it
inversion biholomorphism\/}
$\,\Cal L\smallsetminus N\to\Cal L^*\smallsetminus N^*$ given by
$\,(y,\zx)\mapsto(y,\zx^{-1})$, in the notation of \f{\cn.2}, where
$\,\zx^{-1}\in\Cal L_y^*$ is the unique $\,\bbC$-linear functional
$\,\Cal L_y\to\bbC\,$ that sends $\,\zx\,$ to $\,1\,$ (cf. \S\mm). This makes
$\,M\,$ a holomorphic $\,\bbCP^1$ bundle over $\,N\nh$. If $\,N\hs$ is
compact, so is $\,M\nh$, and one then refers to $\,M\,$ as the {\it projective
compactification\/} of $\,\Cal L\hs$.

Equivalently, we could define $\,M\,$ to be the bundle associated with the
principal $\,\,\text{\rm GL}\hs(1,\bbC)$-bundle of $\,\Cal L\,$ via the
obvious multiplicative action of
$\,\,\text{\rm GL}\hs(1,\bbC)=\bbC\smallsetminus\{0\}\,$ on the Riemann sphere
$\,\bbC\cup\{\infty\}$.
\rmk{Remark \a\ib.1}For $\,\Cal L\,$ and $\,M\,$ as above, the inversion
biholomorphism clearly sends the horizontal distribution $\,\Cal H\,$ of any
$\,C^\infty$ linear connection in the line bundle $\,\Cal L\,$ onto the
horizontal distribution, also denoted $\,\Cal H$, of its dual connection in
$\,\Cal L^*$. Thus, $\,\Cal H\,$ has an extension from $\,\Cal L\subset M\,$
to a $\,C^\infty$ distribution on $\,M$.

Also, any Hermitian fibre metric $\,\langle\,,\rangle\,$ in $\,\Cal L\,$ gives
rise to a Hermitian fibre metric $\,\langle\,,\rangle\nh^*$ in $\,\Cal L^*$
obtained, in each fibre, as in \S\mm.
\endrmk
\proclaim{Lemma \a\ib.2}Let the data\/ \f{\am.3} have the properties listed
in the paragraph following\/ \f{\am.3}. For\/ $\,\,U\nh,g\,$ and\/
$\,\vp:U\to\bbR\,$ determined by them as in\/ {\rm\S\am},
\widestnumber\item{(ii)}\roster
\item"(i)"The same assumptions hold for the new set of data obtained if one
leaves\/ $\,\iyp\!,\vp,Q,\ve,\y,m,N,h\,$ unchanged and replaces the function\/
$\,r\hs$ of the variable\/ $\,\vp\,$ by\/ $\,r^*=1/r$, the constant\/ $\,a\,$
by\/ $\,a^*=-\hs a$, and\/ $\,\Cal L,\langle\,,\rangle,\Cal H\,$ by\/
$\,\Cal L^*,\langle\,,\rangle\nh^*,\Cal H\,$ described in\/ {\rm Remark
\a\ib.1}.
\item"(ii)"The construction of\/ {\rm\S\am} applied to the new data in\/
{\rm(i)} leads to analogous objects\/ $\,\,U^*\!,g^*$ and\/
$\,\vp^*:U^*\to\bbR\,$ such that the inversion biholomorphism\/
$\,\Cal L\smallsetminus N\to\Cal L^*\smallsetminus N^*$ sends\/
$\,\,U\nh,g,\vp\,$ onto\/ $\,\,U^*\!,g^*,\vp^*$.
\endroster
\endproclaim
\demo{Proof}The curvature forms of a given connection and its dual differ by
sign, since so do their connection forms (see Remark \a\cn.1) relative to
local sections without zeros having the form $\,w\,$ and $\,w^{-1}$. This
implies (i). Let $\,r,r^*$ now also stand for specific
functions $\,\Cal L\to\bbR\,$ and $\,\Cal L^*\to\bbR\hs$, namely, the norm
functions of our fibre metrics (Remark \a\cn.2). That the inversion
biholomorphism sends $\,\,U\,$ onto $\,\,U^*$ is clear as $\,\,U\nh,\hs U^*$
are given by $\,r_-<r<r_+$ and $\,r^*_-<r^*<r^*_+$, with $\,r^*_\pm=1/r_\pm$,
and that it makes $\,\vp,r\,$ and the restriction of $\,g\,$ to the vertical
distribution $\,\Cal V\hs$ in $\,\,U\,$ correspond to their counterparts in
$\,\,U^*$ is immediate from (a), (b) in \S\mm. By \f{\am.5} and Remark
\a\ib.1, the same holds for $\,g\,$ and the horizontal distribution
$\,\Cal H\,$ in $\,\,U$. Finally, again by \f{\am.5}, $\,\Cal H,\Cal V\hs$ are
orthogonal to each other both in $\,\Cal L\,$ and in $\,\Cal L^*$. Combined
with Remark \a\ib.1, this completes the proof.{\hfill\qd}
\enddemo

\head\S\os. A one-sided boundary condition\endhead
Let $\,Q,\vp,\iy,\iyp\!,\vp_0,a,r\,$ have the following properties:
\vskip4pt{\rm
\hbox{\hskip-.8pt
\vbox{\hbox{\f{\os.1}}\vskip15pt}
\hskip9pt
\vbox{
\hbox{$Q\hs$ is a $\hs C^\infty\nh$ function of the real variable $\hs\vp$,
defined on a half-open}
\hbox{interval $\,\hs\hs\iy$,\hskip4.5ptpositive on its interior
$\,\,\hs\iyp\nh$,\hskip4.5ptand such that at the only}
\hbox{endpoint $\,\hs\vp_0\hs$ of $\,\hs\iy\hs$ we have $\,\,Q\,=\,0\,\,$ and
$\,\,dQ/d\vp\hs=\hs2a\hs\ne\hs0$, while}
\hbox{$r\,$ is a positive $\,C^\infty$ function of $\,\hs\vp\in\iyp\hs$
satisfying equation \f{\am.1},}}}}
\vskip4pt
\noindent
and let $\,r_\pm$ be as in \f{\am.2}. Then
\widestnumber\item{(b)}\roster
\item"(a)"$r_-=0$, i.e., $\,r\to0\,$ as $\,\vp\to\vp_0$, and $\,Q/r^2$ has
a positive limit as $\,\vp\to\vp_0$.
\item"(b)"$\vp\,$ and $\,Q/r^2$ are $\,C^\infty$ functions of
$\,r^2\in[\hs0,r_+^2)\,$ with $\,Q/r^2>0\,$ at $\,r=0$.
\endroster
In fact, $\,Q/(\vp-\vp_0)\,$ is a $\,C^\infty$ function of $\,\vp\in\iy\,$
equal to $\,2a\,$ at $\,\vp=\vp_0$ (Remark \a\pn.1), and so \f{\am.1} implies
that $\,2\hs d\hs[\log r]/d\vp=2a/Q\,$ equals $\,1/(\vp-\vp_0)\,$ plus a
$\,C^\infty$ function of $\,\vp$, i.e., $\,\log r^2$ equals
$\,\log|\vp-\vp_0|\,$ plus a $\,C^\infty$ function of $\,\vp\in\iy$. Hence
$\,r^2/(\vp-\vp_0)\,$ {\it is a\/} $\,C^\infty$ {\it function of\/}
$\,\vp\in\iy\,$ {\it with a nonzero value at\/} $\,\vp_0$. Now, as
$\,Q/|\vp-\vp_0|\,$ and $\,|\vp-\vp_0|/r^2$ both have positive limits as
$\,\vp\to\vp_0$ (the former limit being $\,2|a|$), the same follows for
$\,Q/r^2$, which proves (a). In view of \f{\os.1}, \f{\am.1} and the phrase
italicized above, the assignment $\,\vp\mapsto r^2$ is a $\,C^\infty$
diffeomorphism of $\,\iy\,$ onto $\,[\hs0,r_+^2)$, sending the endpoint
$\,\vp_0$ to $\,0$, and so (a) implies (b).

\head\S\md. Metrics on disk bundles\endhead
Let a set \f{\am.3} of data have the properties listed
in the paragraph following \f{\am.3}, and let \f{\os.1} hold for
$\,Q,\vp,\iy,\iyp\!,\vp_0,a,r\,$ consisting of the same $\,\iyp\!,\vp,Q,a\,$
as in \f{\am.3}, a fixed finite endpoint $\,\vp_0$ of $\,\iyp\nh$, and
$\,\iy=\iyp\cup\hs\{\vp_0\}$. Thus, we require $\,Q\,$ to have a $\,C^\infty$
extension to $\,\iy$. Let us also assume that, in \f{\am.3}, either $\,\ve=0$,
or $\,\ve=\pm\hs1\,$ and $\,\vp_0\ne\y$.

The construction of \S\am\ now yields a triple $\,\,U\nh,g,\vp\,$ such that
$\,\,U\subset\Cal L\,$ is the open set given by $\,0<r<r_+$, for $\,r_+$ as in
\f{\am.2}, where $\,r\hs$ is the norm function of the fibre metric
$\,\langle\,,\rangle\,$ in $\,\Cal L\,$ (Remark \a\cn.2), while $\,g\,$ is a
K\"ahler metric on $\,\,U\,$ and $\,\vp\,$ is a special \krp\ on
$\,(\hs U\nh,g)$, as in \f{\ir.1}.

These $\,g\,$ and $\,\vp\,$ also have $\,C^\infty$ extensions to a metric and
a function, still denoted $\,g,\vp$, on the open set $\,\,\uo\subset\Cal L\,$
given by $\,0\le r<r_+$, that is, on the bundle of open disks of radius
$\,r_+$ in $\,\Cal L\hs$. The resulting triple $\,\,\uo\nh,g,\vp\,$ satisfies
\f{\ir.1} as well.

In fact, by \f{\am.5} our $\,g\,$ is a real fibre metric on
$\,TU=\Cal H\oplus\Cal V\,$ obtained as the direct sum of
$\,\fe\hs\proj^*\nh h\,$ in $\,\Cal H\,$ and
$\,\theta\,\text{\rm Re}\hskip1pt\langle\,,\rangle\,$ in $\,\Cal V$, with
$\,\theta=Q/(ar)^2$ and $\,\fe\,$ as in \f{\am.4}. The required extensions
exist since the distributions $\,\Cal H,\Cal V\hs$ and the fibre metrics
$\,\proj^*\nh h\,$ and $\,\,\text{\rm Re}\hskip1pt\langle\,,\rangle\,$ on them
are defined and of class $\,C^\infty$ everywhere in $\,\Cal L\,$ (cf.
Remark \a\cn.2), while, by (b) in \S\os, the functions $\,\vp,\theta,\fe\,$
have $\,C^\infty$ extensions to $\,\,\uo\nh$, which are positive in the case
of $\,\theta\,$ and $\,\fe\,$ (the latter due to our assumption that
$\,\vp_0\ne\y\,$ unless $\,\ve=0$). Now \f{\ir.1} for $\,\,\uo\nh,g,\vp\,$
follows since $\,\,U\,$ is dense in $\,\,\uo\nh$.

\head\S\ob. The case of an open ball\endhead
\proclaim{Lemma \a\ob.1}Let a set\/ \f{\tm.1} of data satisfying the
conditions listed in the paragraph following\/ \f{\tm.1}, with\/ $\,m\ge2$,
also have the property that\/ $\,\hs\y\,$ is an endpoint of\/ $\,\iyp\nh$,
while the function\/ $\,Q\,$ of the variable\/ $\,\vp\in\iyp$ has a\/
$\,C^\infty$ extension to the half-open interval\/ $\,\iy=\iyp\nh\cup\{\y\}\,$
with\/ $\,Q=0\,$ and\/ $\,dQ/d\vp=2a\ne0\,$ at\/ $\,\vp=\y\hs$.

The open spherical shell\/ $\,\,U\,$ defined in\/ {\rm\S\tm} then is a
punctured ball, i.e., has the inner radius\/ $\,\hs r_-=0$, while\/ $\,g\,$
with\/ \f{\tm.2} and\/ $\,\vp:\hs U\to\bbR\,$ admit\/ $\,C^\infty$ extensions
to a metric/function on the solid ball\/ $\,\,U\cup\{0\}$.
\endproclaim
\demo{Proof}Let the vector fields $\,v,u\,$ on $\,\,U\cup\{0\}\,$ be as in
Remark \a\tb.1 (for our $\,a$), and let $\,\xi,\xi'$ be the $\,1$-forms on
$\,\,U\cup\{0\}\,$ with
$\,\xi=\,\text{\rm Re}\hskip1pt\langle v,\,\cdot\,\rangle\,$ and
$\,\xi'=\,\text{\rm Re}\hskip1pt\langle u,\,\cdot\,\rangle$. Then $\,g\,$ is,
on $\,\,U\nh$, a combination of $\,\xi\otimes\xi+\xi'\otimes\xi'$ and
$\,\,\text{\rm Re}\hskip1pt\langle\,,\rangle\,$ with the coefficients
$\,[Q-2a(\vp-\y)]/(ar)^4$ and $\,2(\vp-\y)/(ar^2)$. In fact,
$\,(\vp-\y)/a>0\,$ due to the assumptions on $\,\ve\,$ in the lines following
\f{\tm.1}, so that, using \f{\pn.4} and \f{\tm.2} we obtain equal values when
both tensors are evaluated on two vectors, one of which is in $\,\Cal H\,$ and
the other in $\,\Cal H\,$ or $\,\Cal V$, and it is also immediate for the
vector fields $\,v,v$, or $\,v,u$, or $\,u,u$, as
$\,\langle v,v\rangle=\langle u,u\rangle=a^2r^2$ and
$\,\,\text{\rm Re}\hskip1pt\langle v,u\rangle=0\,$ (see Remark \a\cn.2).

Both coefficients are $\,C^\infty$ functions of the variable
$\,r^2\in[\hs0,r_+^2)$. In fact, for $\,2(\vp-\y)/(ar^2)\,$ this is clear from
Remark \a\pn.1 (with $\,\sa=r^2$) and \S\os\ with $\,\vp_0=\y\hs$. Now,
$\,Q/r^2$ and $\,2a(\vp-\y)/r^2$, treated as $\,C^\infty$ functions of
$\,r^2\in[\hs0,r_+^2)\,$ (see (b) in \S\os), have the same positive value
at $\,r^2=0$, since \f{\am.1} gives $\,Q/r^2=2a\hs d\vp/d(r^2)$. Thus, their
difference divided by $\,r^2$ is a $\,C^\infty$ function of
$\,r^2\in[\hs0,r_+^2)\,$ (Remark \a\pn.1 for $\,\sa=r^2$).
Positivity of $\,2a(\vp-\y)/r^2$ at $\,r^2=0\,$ also shows that the limit of
$\,g\,$ at $\,0\in V\,$ is positive definite, completing the proof.{\hfill\qd}
\enddemo
\rmk{Remark \a\ob.2}Conditions $\,\ve a>0\,$ and $\,\ve(\vp-\y)>0\,$ for
all $\,\vp\in\iyp\nh$, required in \f{\tm.1}, follow from each other under the
remaining hypotheses of Lemma \a\ob.1. In fact, since $\,\y\,$ is an endpoint
of $\,\iyp$ and $\,Q>0\,$ on $\,\iyp\nh$, while $\,Q=0\,$ at $\,\y\hs$, the
sign of $\,dQ/d\vp\,$ at $\,\y\,$ must be the same as that of $\,\vp-\y\,$ for
$\,\vp\in\iyp\nh$.
\endrmk
\rmk{Remark \a\ob.3}For reasons mentioned at the end of \S\tm, the triples
$\,\,U,g,\vp\,$ constructed above form another special subset of those in
Example \a\xa.4. Conversely, every triple $\,\,U,g,\vp\,$ of Example \a\xa.4,
in which $\,\,U\,$ is an open ball, can also be obtained as described in this
section. This fact, not needed for our argument, follows if one applies the
local-structure theorem established in \S\tl\ to $\,M=\,U\,$ and $\,y=0$.
\endrmk

\head\S\mw. Special K\"ahler-Ricci potentials on \ $\,\bbCP^1$ \
bundles\endhead
Examples of functions $\,\vp\,$ with \f{\ir.1} on compact K\"ahler manifolds,
in all complex dimensions $\,m\ge2$, can be constructed as follows. Suppose
that
\vskip4pt
\hbox{\hskip-.8pt
\vbox{\hbox{\f{\mw.1}}\vskip15pt}
\hskip9pt
\vbox{
\hbox{$[\vp_{\text{\rm min}},\vp_{\text{\rm max}}]\,$ is
a nontrivial closed interval of the variable $\,\vp\,$ with a $\,C^\infty$}
\hbox{function
$\hs[\vp_{\text{\rm min}},\vp_{\text{\rm max}}]\ni\vp\mapsto Q\in\bbR\hs$,
which is positive on the open inter-}
\hbox{val\hskip8pt$(\vp_{\text{\rm min}},\vp_{\text{\rm max}})$\hskip7ptand
vanishes at the
endpoints\hskip7pt$\vp_{\text{\rm min}},
\hs\vp_{\text{\rm max}}$,\hskip4.5ptwhile the}
\hbox{values of $\,\,dQ/d\vp\,\,$ at the endpoints are mutually opposite and
nonzero.}
}}
\vskip4pt
\noindent
Next, let some data \f{\am.3} satisfy all conditions listed in the paragraph
following \f{\am.3}, and have $\,N,\iyp\!,Q,a,\ve,\y\,$ such that $\,N\hs$ is
compact, $\,Q\,$ is the restriction of $\,Q\,$ in \f{\mw.1} to
$\,\iyp=(\vp_{\text{\rm min}},\vp_{\text{\rm max}})$, while $\,dQ/d\vp=2a\,$
at a fixed endpoint
$\,\vp_0$ of $\,[\vp_{\text{\rm min}},\vp_{\text{\rm max}}]$, and either
$\,\ve=0$, or $\,\ve=\pm\hs1\,$ and
$\,\y\notin[\vp_{\text{\rm min}},\vp_{\text{\rm max}}]$.

These assumptions then also hold for a new set of data, analogous to
\f{\am.3}, which is obtained as in Lemma \a\ib.2(i), so that
$\,r,a,\Cal L,\langle\,,\rangle\,$ are replaced by some specific
$\,r^*,a^*,\Cal L^*,\langle\,,\rangle\nh^*$. 
Using \S\md\ with $\,r_+=\infty\,$ (cf. Remark \a\mw.1 below), we see
that the construction of \S\am\ applied to our original data (or, the new
data) leads to a K\"ahler metric $\,g\,$ on $\,\Cal L\,$ (or, $\,g^*$ on
$\,\Cal L^*$) along with a special \krp\ $\,\vp\,$ on $\,(\Cal L,g)\,$ (or,
$\,\vp^*$ on $\,(\Cal L^*,g^*)$).

>From Lemma \a\ib.2(ii), with $\,\,U=\Cal L\smallsetminus N\,$ and
$\,\,U^*=\Cal L^*\smallsetminus N^*$ it now follows that $\,g,\vp\,$ and
$\,g^*,\vp^*$ together form a K\"ahler metric and a special \krp, again
denoted $\,g,\vp$, on the projective compactification $\,M\,$ of
$\,\Cal L\,$ (cf. \S\ib).
\rmk{Remark \a\mw.1}Assuming \f{\mw.1} and choosing $\,r,r_+,r_-$ with
\f{\am.1} - \f{\am.2} for
$\,\iyp=(\vp_{\text{\rm min}},\vp_{\text{\rm max}})\,$ and $\,a\,$ such that
$\,\pm\hs2a\,$ are the values of $\,dQ/d\vp\,$ at the endpoints of
$\,\iyp\nh$, we have $\,r_-=0\,$ and $\,r_+=\infty$. In fact, (a) in \S\os\
for these $\,r,a\,$ and one endpoint gives $\,r_-=0$, while, when applied to
$\,r^*=1/r$, $\,a^*=-\hs a\,$ and the other endpoint, it yields $\,r^*_-=0$,
i.e., $\,r_+=\infty$.
\endrmk
\rmk{Remark \a\mw.2}The manifold $\,M\,$ constructed above contains two
holomorphically embedded copies of $\,N\nh$, namely the zero sections
$\,N\subset\Cal L\,$ and $\,N^*\subset\Cal L^*$. By (c) in \S\am, $\,d\vp\,$
vanishes precisely at the zeros of the function $\,Q\,$ prescribed in
\f{\mw.1} and treated as a function on $\,\,U=\Cal L\smallsetminus N\,$ (or
$\,\,U^*=\Cal L^*\smallsetminus N^*$), via the dependence of $\,\vp\,$ on the
norm function $\,r\hs$ (or, $\,r^*$). Since \f{\mw.1} and (a) in \S\os\ show
that $\,Q\,$ has in $\,\,U\,$ (or, $\,\,U^*$) the same zeros as $\,r\hs$ (or,
$\,r^*$), $\,d\vp\,$ vanishes in $\,M\,$ along the zero sections only, i.e.,
$\,\vp:M\to\bbR\,$ has two critical manifolds: $\,N\hs$ and $\,N^*$.
\endrmk

\head\S\cp. Special K\"ahler-Ricci potentials on \ $\,\bbCP^m$\endhead
The following construction is a more explicit version of Example \a\xa.5; see
Remark \a\cp.1 below.

Let us assume \f{\mw.1}, and let the conditions listed in the paragraph
following \f{\tm.1} hold for some data \f{\tm.1}. Furthermore, let
$\,\iyp\!,Q,a,\y\,$ in \f{\tm.1} be such that
$\,\iyp=(\vp_{\text{\rm min}},\vp_{\text{\rm max}})\,$ and $\,Q\,$ is the
restriction of $\,Q\,$ in \f{\mw.1} to $\,\iyp\nh$, while
$\,\y\,$ is an endpoint of $\,\iyp$ and $\,2a\,$ is the value of $\,dQ/d\vp\,$
at $\,\vp=\y\hs$.

According to Lemmas \a\ob.1 and \a\tm.1, the construction of \S\tm\ applied to
these data is a special case of that in \S\am, and yields a K\"ahler metric
$\,g\,$ on the complex vector space $\,V$ along with a special \krp\
$\,\vp\,$ on $\,(V\nh,g)$. (The solid ball in Lemma \a\ob.1 is $\,V$ itself,
since $\,r_+=\infty\,$ by Remark \a\mw.1.) The objects \f{\am.3} leading, as
in \S\am, to $\,g\,$ and $\,\vp\,$ on
$\,M'=V\smallsetminus\{0\}=\Cal L\smallsetminus N\,$ are our
$\iyp\!,\vp,Q,r,a,\ve,\y\hs,m\,$ in \f{\tm.1}, the projective space
$\,N\hs$ of $\,V\nh$, the metric $\,h\,$ such that $\,|a|h\,$ is the
Fubini-Study metric, and the
tautological bundle $\,\Cal L\,$ over $\,N\hs$ with the standard fibre metric
and the horizontal distribution $\,\Cal H\,$ of the canonical connection
(\S\tb).

Let $\,N^*$ stand for $\,N\hs$ treated as the zero section
$\,N^*\subset\Cal L^*$ in the dual bundle $\,\Cal L^*\nh$, cf. \f{\cn.2}. By
Lemma \a\ib.2(ii), the biholomorphism $\,V\smallsetminus\{0\}
=\Cal L\smallsetminus N\to\Cal L^*\smallsetminus N^*$ identifies $\,g,\vp\,$
on $\,V\smallsetminus\{0\}\,$ with $\,g^*,\vp^*$ on
$\,\Cal L^*\smallsetminus N^*$ which are obtained as in \S\am\ from the new
data $\,\iyp\!,\vp,Q,1/r\hs,-\hs a,\ve,\y\hs,m,N,h,
\Cal L^*,\Cal H,\langle\,,\rangle\nh^*$ analogous to \f{\am.3}. Therefore,
$\,g,\vp\,$ give rise to a K\"ahler metric and a special \krp, still denoted
$\,g,\vp$, on the complex manifold $\,M\,$ obtained from the disjoint union
$\,V\cup\Cal L^*$ by using the above biholomorphism to identify the open sets
$\,V\smallsetminus\{0\}\subset V\,$ and
$\,\Cal L^*\smallsetminus N^*\subset\Cal L^*$. (Note that, applying \S\md\
to the new data and the endpoint of $\,\iyp$ other than $\,\y\hs$, we
can extend $\,g^*,\vp^*$ from $\,\Cal L^*\smallsetminus N^*$ to $\,\Cal L^*$.)

This $\,M\,$ is clearly biholomorphic to the projective space
$\,\bbP\approx\bbCP^m$ formed by all complex lines through $\,0\,$ in
$\,V\times\bbC\,$ (cf. \S\tb). Namely, we have holomorphic embeddings
$\,V\smallsetminus\{0\}\to\hs\bbP\,$ and
$\,\Cal L^*\smallsetminus N^*\to\hs\bbP\,$ sending any
$\,\zx\in V\smallsetminus\{0\}\,$ to the line spanned by
$\,(\zx,1)\in V\times\bbC\hs$, and any
$\,(y,\zeta)\in\Cal L^*\smallsetminus N^*$ (see \f{\cn.2}) to the graph of
the linear functional $\,\zeta\,$ on the line $\,y\subset V\,$ (cf. \S\tb);
the graph is itself a line through zero in
$\,y\times\bbC\subset V\times\bbC\hs$, that is, an element of $\,\bbP$. The
resulting transition mapping, obtained from the former embedding followed by
the inverse of the latter, is precisely the biholomorphism
$\,V\smallsetminus\{0\}\to\Cal L^*\smallsetminus N^*$ we just used, as it
takes $\,\zx\in V\smallsetminus\{0\}\,$ to
$\,(y,\zeta)\,$ such that $\,y\subset V\,$ is a line through zero and
the graph of the functional $\,\zeta:y\to\bbC\,$ is spanned by (i.e.,
contains) $\,(\zx,1)$, which means that $\,\zx\,$ spans $\,y\,$ and
$\,\zeta\,$ sends $\,\zx\,$ to $\,1$, i.e., $\,\zeta=\zx^{-1}$ (notation of
\S\ib).
\rmk{Remark \a\cp.1}The construction just described is a special case of
Example \a\xa.5, as one sees restricting $\,g,\vp\,$ to $\,V\subset M\,$ and
using Remark \a\ob.3. Thus, by Remark \a\xa.6, $\,\vp\,$ obtained here has two
critical manifolds, which realize case 2) of (ii) in \S\bo\ below.

Conversely, every triple $\,\,\bbCP^m,g,\vp\,$ of Example \a\xa.5 can also be
obtained as described in this section. This fact, which will not be used, is
an immediate consequence of Theorem \a\gc.2 in \S\gc.
\endrmk

\head\S\dc. Dimensions of critical manifolds\endhead
Suppose that $\,\vp\,$ is a special \krp\ on a K\"ahler manifold $\,(M,g)\,$
of complex dimension $\,m\,$ (see \f{\ir.1}) and $\,N\hs$ is a critical
manifold of $\,\vp$, cf. Remark \a\cz.3(ii), while $\,Q=g(\navp,\navp)$,
$\,Y=\Delta\vp$, and $\,\si,\ta:M'\to\bbR\,$ are as in \f{\sm.3}, $\,M'$
denoting the open set on which $\,d\vp\ne0\,$ (i.e., $\,Q>0$). Thus, either
$\,\si=0\,$ identically on $\,M'$ or $\,\si\ne0\,$ everywhere in $\,M'$ (see
Lemma \a\sm.2). To refer to the former case, we will just write $\,\si=0$,
while the latter one will be tacitly assumed whenever we mention the constant
$\,\y\,$ defined (only when $\,\si\ne0$) in Lemma \a\sm.2.

We have $\,dY=0\,$ wherever $\,d\vp=0$, as
$\,dY=-\hs2\,\ri\hs(\navp,\,\cdot\,)\,$ (cf. the two lines  preceding Remark
\a\sm.1), so that $\,Y=\Delta\vp=\langle\hs g,\nabla d\vp\rangle\,$ is
constant on $\,N\nh$. Letting $\,\vp_0$ be the constant value of $\,\vp\,$ on
$\,N\nh$, we define a real constant $\,a$, depending on $\,N\nh$, by
$$\text{\rm$2a=Y$\hskip5pton\hskip4pt$N$\hskip6ptif\hskip5pt$\si=0$
\hskip4ptor\hskip6pt$\vp_0\ne\y\hs$,\hskip6ptand\hskip6pt$2a=Y/m$\hskip5pton
\hskip4pt$N$\hskip4ptif\hskip6pt$\vp_0=\y\hs$.}\ff\dc.1$$
Then $\,\ta(x)\to a\,$ as $\,x\to y\in N\nh$, where $\,x\,$ is a variable
point of $\,M'\nh$. Also,
$$\text{\rm$a\ne0$\hskip6ptand\hskip6pt$\nabla_{\!w}v=aw$\hskip6ptfor\ every\
vector\hskip5pt$w$\hskip5ptnormal\ to\hskip4pt$N$\hskip3ptat\ any\ point,}
\hskip12pt\ff\dc.2$$
where $\,v=\navp$. Furthermore, one of the following two cases must occur:
$$\alignedat2
&\text{\rm a)}\quad&&
\text{\rm$N$\hskip4ptis\ a\ complex\ submanifold\ of\ complex\ codimension
\hskip3pt$1$\hskip5ptin\hskip4pt$M$,\hskip5ptor\hskip13pt}\\
&\text{\rm b)}\quad&&
\text{\rm$N$\hskip4ptconsists\ of\ a\ single\ point.}
\endalignedat\ff\dc.3$$
Finally,
\widestnumber\item{(iii)}\roster
\item"(i)"In case a) of \f{\dc.3}, $\,\si=0\,$ on $\,M'$ or $\,\vp_0\ne\y\hs$.
\ In case b) of \f{\dc.3}, $\,\vp_0=\y\hs$.
\item"(ii)"If $\,\si=0\,$ and $\,m\ge2$, then no critical manifold of
$\,\vp\,$ is a one-point set.
\item"(iii)"In the case where $\,\si\,$ is not identically zero, a point
$\,y\in M\,$ has $\,\vp(y)=\y\,$ if and only if $\,\{y\}\,$ is a critical
manifold of $\,\vp$.
\item"(iv)"In case a) of \f{\dc.3} with $\,m\ge2$, the limit relation
$\,\Cal H_x\to T_yN\,$ as $\,x\to y$, for any $\,y\in N\nh$, holds in an
appropriate Grassmannian bundle, with $\,\Cal H\,$ as in \f{\sm.1}, $\,x\,$
being a variable point of $\,M'\nh$.
\endroster
In fact, let us fix $\,y\in N\nh$. Since $\,M'$ is dense in $\,M\,$ (Remark
\a\kg.4), we may choose a sequence of points in $\,M'$ converging to $\,y\,$
and, at each point $\,x\,$ of the sequence, an orthonormal basis of
$\,T_xM\,$ formed by eigenvectors of $\,(\nabla d\vp)(x)$, the last two of
which correspond to the eigenvalue $\,\ta(x)$, and the others to $\,\si(x)$,
cf. \f{\sm.3}. A subsequence of this sequence of bases converges, in a
suitable frame bundle, to an orthonormal basis of $\,T_xM\,$ that has all the
properties just listed for $\,x=y$, with some eigenvalues $\,\ta_0,\si_0$ that
are limits of the $\,\ta(x)\,$ and $\,\si(x)$.

If $\,\si=0\,$ or $\,\vp(y)\ne\y\hs$, then $\,\si_0=0$, which is obvious when
$\,\si=0\,$ and, if $\,\vp(y)\ne\y\hs$, follows if we let $\,x\to y\,$ in
$\,Q=2(\vp-\y)\si\,$ (see Lemma \a\sm.2).

If $\,\vp(y)=\y\hs$, we must have $\,\ta_0=\si_0$. In fact, let us choose a
curve $\,\sa\mapsto x(\sa)\,$ as in Remark \a\sm.4(ii), so that
$\,\dot\vp\ne0\,$ for all $\,\sa\ne0\,$ close to $\,0$. Also,
$\,\ta(x(\sa))\to\ta_0$ as $\,\sa\to0$, and similarly for $\,\si$, since
$\,\ta_0,\si_0$ are the limits of all convergent sequences of such values
$\,\ta(x(\sa)),\hs\si(x(\sa))\,$ (to see this, consider two separate cases,
$\,\ta_0=\si_0$ and $\,\ta_0\ne\si_0$). Hence, by l'Hospital's
rule, $\,Q/(\vp-\y)\,$ evaluated at $\,x(\sa)\,$ tends, as $\,\sa\to0$, to the
limit of $\,\dot Q/(\vp-\y)\dot{\,}=2\ta\dot\vp/\dot\vp\,$ (see Remark
\a\sm.4(i)), that is, to $\,2\ta_0$, while, by Lemma \a\sm.2,
$\,Q/(\vp-\y)\to2\si_0$ as $\,\sa\to0$.

Our $\,\ta_0,\si_0$ are the eigenvalues of $\,(\nabla d\vp)(y)$, with the
multiplicities $\,2(m-1)\,$ and $\,2$. Hence the constant value of
$\,Y=\Delta\vp\,$ on $\,N\hs$ equals $\,2\ta_0+2(m-1)\hs\si_0$, that is, its
value at $\,y\in N\nh$. As $\,\si_0=0$ in one case discussed above and
$\,\si_0=\ta_0$ in the other, we obtain $\,\ta_0=a\,$ for $\,a\,$ given by
\f{\dc.1}. Also, $\,a\ne0$, or else we would have $\,\si_0=\ta_0=0\,$ in both
cases, contradicting the relation $\,(\nabla d\vp)(y)\ne0\,$ in Remark
\a\kg.4. According to \f{\pn.3}, the complex space $\,T_yN\,$ thus is the
orthogonal direct sum of two subspaces, of which one is the eigenspace of
$\,(\nabla v)(y)\,$ for the unique nonzero eigenvalue $\,\ta_0=a$, and the
other is the kernel of $\,(\nabla v)(y)\,$ (trivial when $\,\si_0=\ta_0$,
of complex codimension one when $\,\si_0=0$), so that Remark \a\cz.3(iii)
implies \f{\dc.2}, \f{\dc.3} and (i). Next, (ii) follows since its premise
precludes the case $\,\si_0=\ta_0\ne0$, and, as we saw, the remaining case
$\,\si_0=0\,$ yields \f{\dc.3.a} for {\it every\/} critical manifold $\,N\nh$, and, similarly, (iii) is obvious from (i) and the fact that, by Lemma \a\sm.2,
we have $\,\vp\ne\y\,$ at all points with $\,d\vp\ne0$. Finally, the limit
relations
$\,\ta(x)\to a\,$ and (iv) follow since the convergence involving the
$\,\ta(x)\,$ and $\,\ta_0$, as well as that for orhonormal bases, was
established for some subsequence of {\it any\/} given sequence of points
$\,x\in M'$ tending to $\,y\in N\nh$.
\proclaim{Lemma \a\dc.1}Suppose that\/ $\,N\hs$ is a critical manifold of a
function\/ $\,\vp\,$ satisfying\/ \f{\ir.1} on a K\"ahler manifold\/
$\,(M,g)$, cf. {\rm Remark \a\cz.3(ii)}, and\/
$\,[\hs0,\ell)\ni\sa\mapsto x(\sa)\in M\,$ is a unit-speed geodesic with\/
$\,\dot x(0)\,$ normal to\/ $\,N\hs$ at\/ $\,x(0)=y\in N\nh$,
where\/ $\,\dot x=\hs dx/d\sa$, and such that\/ $\,d\vp\ne0\,$ at\/
$\,x(\sa)\,$ for all\/ $\,\sa\in(0,\ell)$. If we set\/ $\,v=\navp\,$ and let\/
$\,\,\sgn\,a=\pm\hs1\,$ stand for the sign of\/ $\,a\,$ in\/ \f{\dc.1} --
\f{\dc.2}, then\/ $\,\dot x=(\sgn\,a)\hs v/|v|\,$ at\/ $\,x(\sa)$, for every\/
$\,\sa\in(0,\ell)$, and, for\/ $\,\sa\in[\hs0,\ell)$,
$$d\vp/d\sa\,=\,(\hs\sgn\,a)\hs\sqrt{Q\,},\hskip10pt\text{\rm with\ the\
initial\ value}\hskip6pt\vp\,=\,\vp_0\hskip6pt\text{\rm at}\hskip6pt\sa=0\hs,
\ff\dc.4$$
where\/ $\,d\vp/d\sa=d\hs[\vp(x(\sa))]/d\sa\,$ and\/ $\,\vp_0$ is the constant
value of\/ $\,\vp\,$ on\/ $\,N\nh$.
\endproclaim
In fact, $\,\nabla_{\!v}v=\ta\hskip.4ptv\,$ by \f{\sm.4} with $\,w=v$, and
$\,\nabla_{\!w}v=aw\,$ for $\,w=\dot x(0)\,$ by \f{\dc.2}, so that Lemma
\a\gd.4(b) applied to the Levi-Civita connection $\,\nabla\,$ of $\,g$, the
geodesic segment $\,X\subset M\,$ which is the image of $\,\sa\mapsto x(\sa)$,
and our $\,y,v,a\,$ yields $\,\dot x=\pm\hs v/|v|\,$ for some sign $\,\pm\,$
and all $\,\sa\in(0,\ell)$. Remark \a\sm.4(i) now
gives $\,\dot\vp=g(v,\dot x)\,$ and $\,\dot Q=2\ta\dot\vp$. At any $\,\sa>0\,$
close to $\,0\,$ we thus have $\,\pm\hs\dot\vp>0\,$ and $\,a\ta>0\,$ (due to
the relation $\,\ta(x)\to a\,$ preceding \f{\dc.2}), so that
$\,\pm\hs a\dot Q>0$, which implies $\,\pm\hs a>0$, since $\,\dot Q>0$. (Note
that $\,Q(x(\sa))>0\,$ for such $\,\sa$, as $\,Q=g(\navp,\navp)$, while
$\,Q(x(0))=0$.) Finally, \f{\dc.4} now follows since
$\,d\vp/d\sa=\hs d_{\dot x}\vp\,$ and, by \f{\sm.2},
$\,|\navp|=\sqrt Q$.{\hfill\qd}

\head\S\cg. A consequence of Gauss's Lemma\endhead
The {\it normal exponential mapping\/} of a submanifold $\,N\hs$ of a
Riemannian manifold $\,(M,g)\,$ is the restriction of
$\,\text{\rm Exp}:\hs\uexp\to M\,$ to the set $\,\,\uexp\cap\Cal L\hs$, with
$\,\Cal L\,$ denoting the total space of the normal bundle of $\,N\hs$ (see
\f{\cn.2}) and $\,\uexp\subset TM\,$ defined as in Remark \a\gd.3 for the
Levi-Civita connection $\,\nabla\,$ of $\,(M,g)$.

For $\,M,g,N\nh,\Cal L\,$ as above and $\,y\in N\nh$, let $\,\sa\,$ be the
norm function (Remark \a\cn.2) of the real fibre metric in $\,\Cal L\,$
obtained by restricting $\,g\,$ to $\,\Cal L\hs$. The inverse mapping theorem
allows us to choose a connected neighborhood $\,N'$ of $\,y\,$ in $\,N\hs$ and
a number $\,\ell\in(0,\infty)\,$ such that, for the open subset $\,\,U'$ of
$\,\Cal L'$ given by $\,0\le\sa<\ell$, where $\,\Cal L'$ is the portion of
$\,\Cal L\,$ lying over $\,N'\nh$, we have $\,\,U'\subset\hs\uexp$ and the
normal exponential mapping sends $\,\,U'$ diffeomorphically onto an open set
in $\,M$.

The following fact is well-known (and also immediate from (d) in \S\cv):
\proclaim{Gauss's Lemma}Under these assumptions, all half-open geodesic
segments of length\/ $\,\ell$, emanating from\/ $\,N'$ in directions normal
to\/ $\,N\nh$, intersect orthogonally the\/ $\,\text{\rm Exp}$-images of all
level sets of the norm function restricted to\/ $\,\,U'\nh$.{\hfill\qd}
\endproclaim
One of its consequences is
\proclaim{Lemma \a\cg.1}For a function\/ $\,\vp\,$ satisfying\/ \f{\ir.1} on a
K\"ahler manifold\/ $\,(M,g)$, let\/ $\,Q=g(\navp,\navp)\,$ and let\/
$\,\ta:M'\to\bbR\,$ be characterized by\/ \f{\sm.3}, with\/ $\,M'$ standing
for the open set on which\/ $\,d\vp\ne0$. Then
\widestnumber\item{(c)}\roster
\item"(a)"$\ta\,$ has a unique extension to a\/ $\,C^\infty$ function\/
$\,M\to\bbR\hs$, also denoted\/ $\,\ta$.
\item"(b)"Every point of\/ $\,M\,$ has a neighborhood\/ $\,\,U\,$ on which\/
$\,Q\,$ is a\/ $\,C^\infty$ function of\/ $\,\vp$, that is, a composite
consisting of\/ $\,\vp\,$ followed by a\/ $\,C^\infty$ function\/
$\,\vp\mapsto Q\,$ defined on a suitable interval of the variable\/ $\,\vp\,$
and such that\/ $\,dQ/d\vp=2\ta\,$ for\/ $\,dQ/d\vp\,$ and\/ $\,\ta\,$ treated
as functions on\/ $\,\,U\nh$.
\endroster
\endproclaim
\demo{Proof}At points with $\,d\vp\ne0$, (b) is obvious from Remark
\a\sm.5(a). Suppose now that $\,y\in M\,$ is a point at which $\,d\vp=0$, and
let $\,N\hs$ be the critical manifold of $\,\vp\,$ containing $\,y\,$ (cf.
Remark \a\cz.3(ii)). We may choose $\,N'\nh,\ell,\hs U'$ as in the second
paragraph of this section and, making $\,N'$ and $\,\ell\,$ smaller if
necessary, also require that $\,d\vp\ne0\,$ at every point of
$\,\text{\rm Exp}\hs(\hs U'\smallsetminus N')$. (Cf. Lemma \a\cz.2(a) for
$\,u=J(\navp)$.)

The gradients $\,v=\navp\,$ and
$\,\nabla Q=2\ta\hskip.4ptv\,$ (see \f{\sm.5.i}), which, by Lemma \a\dc.1, are
tangent to the geodesic segments mentioned in Gauss's Lemma, must therefore be
normal to the $\,\text{\rm Exp}$-images of all level sets of the norm function
restricted to $\,\,U'\nh$. Any such level set is a bundle of positive\diml\
spheres over $\,N'$ (cf. the inequality in Lemma \a\cz.2(c) for
$\,u=J(\navp)$), unless it is the zero section $\,N'\nh$, i.e., the zero
level; therefore, it is connected, and so $\,\vp,Q\,$ must both be constant
along its $\,\text{\rm Exp}$-image. Thus, both $\,\vp\,$ and
$\,Q$, restricted to $\,\text{\rm Exp}\hs(\hs U')\,$ and then pulled back to
$\,\,U'$ via $\,\text{\rm Exp}\hs$, are functions of the norm function.

Let $\,(-\ell,\ell)\ni\sa\mapsto x(\sa)\in M\,$ be any unit-speed geodesic
such that $\,x(0)\in N'$ and $\,\dot x(0)\,$ is normal to $\,N\hs$ at
$\,x(0)$, where $\,\dot x=\hs dx/d\sa$. As
$\,x(\sa)=\,\text{\rm Exp}\hs(x(0),\hs\sa\hs\dot x(0))\,$ and the value of the
norm function at $\,(x(0),\hs\sa\dot x(0))\,$ is $\,|\sa|$, it follows that
$\,\vp,Q\,$ treated as $\,C^\infty$ functions of the variable
$\,\sa\in(-\ell,\ell)\,$ (via the substitutions $\,\vp(x(\sa)),\hs Q(x(\sa))$)
depend just on $\,|\sa|$, i.e., are even. Their restrictions to
$\,[\hs0,\ell)\,$ describe how their $\,\text{\rm Exp}$-pull\-backs depend on
the norm function (also denoted $\,\sa$). By \f{\dc.4}, the dependence of
$\,\vp\,$ on $\,\sa\,$ is homeomorphic, i.e., $\,Q\,$ restricted to
$\,\text{\rm Exp}\hs(\hs U')\,$ is also a function of $\,\vp$. Finally,
$\,\dsq\hs\vp/d\sa^2\ne0\,$ at $\,\sa=0\,$ in view of Remark \a\sm.4(ii),
since, by \f{\dc.2} and \f{\pn.3}, $\,\dot x(0)\,$ is an eigenvector of
$\,(\nabla d\vp)(y)\,$ for the eigenvalue $\,a\ne0$. Assertion (b) for the
point $\,y\,$ is therefore immediate from Remark \a\pn.2. Finally, relation
$\,dQ/d\vp=2\ta$, valid locally in $\,M'$ (see Remark \a\sm.5(a)) can now be
used to define a $\,C^\infty$ extension of $\,\ta\,$ to a suitable
neighborhood of any given point in $\,M\,$ and, as every such extension is
unique (due to denseness of $\,M'\nh$, cf. Remark \a\kg.4), all such
extensions together form a function $\,\ta:M\to\bbR\hs$. This completes the
proof.{\hfill\qd}
\enddemo
For $\,M,g,\vp\,$ as in Lemma \a\cg.1 and $\,\ta,\si:M'\to\bbR\,$ given by
\f{\sm.3}, the unique $\,C^\infty$ extension of $\,\ta\,$ to $\,M\,$ provided
by  Lemma \a\cg.1(a) leads to a similar extension of $\,\si$.
In fact, $\,\si:M\to\bbR\,$ then is defined either by \f{\sm.5.ii} with
$\,m\ge2$, or by $\,\si=0\,$ when $\,M\,$ is of complex dimension $\,1$. Both
extensions are constant on every critical manifold $\,N\hs$ of $\,\vp$.

Furthermore, $\,\ta=a\,$ on $\,N\nh$, for the constant $\,a\ne0\,$ depending
on $\,N\hs$ as in \f{\dc.1} -- \f{\dc.2} (due to the relation
$\,\ta(x)\to a\,$ preceding \f{\dc.2}). Therefore, $\,\si=a\,$ on $\,N\hs$
when $\,N\hs$ consists of a single point, and $\,\si=0\,$ on $\,N\hs$
otherwise; this is clear from \f{\sm.5.ii} restricted to $\,N\hs$ (so that
$\,\ta=a$) along with \f{\dc.1} and (i) in \S\dc.

\head\S\pf. Isometric actions of the circle\endhead
For a $\,C^2$ function $\,\vp\,$ on a Riemannian
manifold $\,(M,g)$, let $\,\,\text{\rm Crit}^1(\vp)\,$ be the set of those
critical points $\,y\,$ of $\,\vp\,$ at which the Hessian
$\,\,\text{\rm Hess\hskip.3pt}_y\vp\,$ has exactly one nonzero eigenvalue (of
any multiplicity). Thus,
$$\text{\rm Hess\hskip.3pt}_y\vp\quad\text{\rm is\hskip6ptsemidefinite\hskip6pt
for\hskip6ptevery}\quad y\in\,\text{\rm Crit}^1\hs(\vp)\hs.\ff\pf.1$$
If $\,y\in\,\text{\rm Crit}^1(\vp)\,$ and $\,a\,$ is the nonzero eigenvalue of
$\,\,\text{\rm Hess\hskip.3pt}_y\vp$, while $\,N\hs$ is the critical manifold
of $\,\vp\,$ containing $\,y$, and $\,u=Jv\,$ with $\,v=\navp$, then, for any
$\,\zx\in T_yM$,
$$\text{\rm$\nabla_{\!\zx}u=0\,$ \ if \ $\,\zx\in T_yN\nh$, \ and \
$\,\nabla_{\!\zx}u=aJ\zx\,$ \ if \ $\,\zx\in(T_yN)^\perp\nh$.}\ff\pf.2$$
In fact, as $\,\nabla v\,$ commutes with
$\,J\,$ (see Lemma \a\kg.2(ii) and \f{\kg.1.a}), so does
$\,\nabla u=J\circ(\nabla v)$, by \f{\kg.1.b}. Thus, $\,(\nabla u)(y)\,$ is
complex-linear has the same eigenvectors as $\,(\nabla v)(y)\,$ (or
$\,(\nabla d\vp)(y)$, cf. \f{\pn.3}), its eigenvalues being $\,i\,$ times
those of $\,(\nabla d\vp)(y)$ Now \f{\pf.2} is obvious from Remark
\a\cz.3(iii) and \f{\dc.2}. Next, for the set $\,\,\text{\rm Crit}\hs(\vp)\,$
of all critical points of $\,\vp$, Remark
\a\cz.3(iii), \f{\dc.2} and \f{\pn.3} give
$$\text{\rm Crit}\hs(\vp)\,=\,\,\text{\rm Crit}^1(\vp)\quad\text{\rm if}
\hskip7pt\vp\hskip7pt\text{\rm satisfies\ \hs\f{\ir.1}\hs\ on\ a\ K\"ahler\
manifold.}\ff\pf.3$$
\proclaim{Lemma \a\pf.1}Let\/ $\,\vp\,$ be a \kip\ on a K\"ahler manifold\/
$\,(M,g)$, cf. {\rm\S\pn}, and let a point\/ $\,y\in M\,$ lie in the set\/
$\,\,\text{\rm Crit}^1(\vp)\,$ defined above, so that\/ $\,u=J(\navp)\,$ is a
Killing field and\/ $\,u(y)=0$. If\/ $\,\,U,\,U'$ are chosen as in\/ {\rm
Lemma \a\cz.1}, for these\/ $\,u\,$ and\/ $\,y$, then the flow of\/ $\,u\,$
restricted to\/ $\,\,U\,$ is periodic, i.e., represents an isometric action
on\/ $\,\,U\,$ of the circle group $\,S^1$. The minimum period of the flow
of\/ $\,u\,$ on\/ $\,\,U\,$ equals\/ $\,2\pi/|a|$, where\/ $\,a\,$ is the
nonzero eigenvalue of\/ $\hs\,\nabla d\vp\,$ at\/ $\,y$.
\endproclaim
In fact, according to Lemma \a\gd.5, $\,u\,$ restricted to $\,\,U\,$ is the
$\,\,\e_{\hskip.4pty}$-image of the linear vector field on $\,\,U'$ given by
the skew-adjoint (and hence diagonalizable) operator
$\,\zx\mapsto\nabla_{\!\zx}u\,$ with the eigenvalues $\,ai\,$ and $\,0$,
or just $\,ai\,$ (see \f{\pf.2}).{\hfill\qd}
\proclaim{Corollary \a\pf.2}Let\/ $\,\vp:M\to\bbR\,$ be a \kip\ on a complete
K\"ahler manifold\/ $\,(M,g)\,$ such that the set\/
$\,\,\text{\rm Crit}^1(\vp)\,$ defined above is nonempty. Then
\widestnumber\item{(ii)}\roster
\item"(i)"The flow of the Killing vector field\/ $\,u=J(\navp)\,$ is periodic,
i.e., constitutes an isometric\/ $\,S^1$ action on\/ $\,(M,g)$.
\item"(ii)"The absolute value of the nonzero eigenvalue of\/ $\,\nabla d\vp\,$
is the same at all points of\/ $\,\,\text{\rm Crit}^1(\vp)$.
\endroster
\endproclaim
In fact, (i), (ii) are both obvious from Lemma \a\pf.1 and the unique
continuation property for isometries (Remark \a\kg.1): in (ii),
$\,2\pi/|\ta(y)|\,$ is, by Lemma \a\pf.1, the minimum period of the flow of
$\,u$, and so it is the same for all
$\,y\in\,\text{\rm Crit}^1(\vp)$.{\hfill\qd}
\medskip
\proclaim{Corollary \a\pf.3}Let\/ $\,\vp\,$ satisfy\/ \f{\ir.1} on a complete
K\"ahler manifold\/ $\,(M,g)$, and let\/ $\,\ta:M\to\bbR\,$ be the continuous
extension to\/ $\,M\nh$, described in\/ {\rm Lemma \a\cg.1}, of the eigenvalue
function\/ $\,\ta\,$ in\/ \f{\sm.3}. Then the restriction of the function\/
$\,|\ta|\,$ to the set\/ $\,\,\text{\rm Crit}\hs(\vp)\,$ of critical points
of\/ $\,\vp\,$ is constant and positive.
\endproclaim
This is clear from Corollary \a\pf.2(ii) and \f{\pf.3}. (Constancy of
$\,\ta\,$ on each connected component $\,N\hs$ of
$\,\,\text{\rm Crit}\hs(\vp)\,$ has also been shown at the end of
\S\cg.){\hfill\qd}

\head\S\bo. Boundary conditions\endhead
Let $\,\vp\,$ be a special \krp\ on a compact K\"ahler manifold $\,(M,g)\,$ of
complex dimension $\,m\ge1$, cf. \f{\ir.1}. Then
\widestnumber\item{(iv)}\roster
\item"(i)"$\vp\,$ has exactly two critical manifolds, defined as in
Remark \a\cz.3(ii), and they are the $\,\vp$-pre\-im\-ages of its extremum
values $\,\vp_{\text{\rm max}}$ and $\,\vp_{\text{\rm min}}$.
\item"(ii)"One of the following two cases must occur:
\itemitem{\rm 1)}Both critical manifolds of $\,\vp\,$ are of complex
codimension one;
\itemitem{\rm 2)}One critical manifold of $\,\vp\,$ is of complex
codimension $\,1$, and the other consists of a single point.
\item"(iii)"$Q=g(\navp,\navp)\,$ is a $\,C^\infty$ function of $\,\vp$,
that is, a composite consisting of $\,\vp\,$ followed by a $\,C^\infty$
function
$\,[\vp_{\text{\rm min}},\vp_{\text{\rm max}}]\ni\vp\mapsto Q\in\bbR\hs$.
\item"(iv)"The function
$\,[\vp_{\text{\rm min}},\vp_{\text{\rm max}}]\ni\vp\mapsto Q\in\bbR\,$ in
(iii) satisfies \f{\mw.1}.
\endroster
In fact, by Example \a\mb.1, \f{\pf.3}, \f{\pf.1} and the inequality in
Lemma \a\cz.2(c) with $\,u=J(\navp)$, our $\,M\,$
and $\,\vp\,$ satisfy the assumptions, and hence the conclusions, of
Corollary \a\mb.5. This gives (i). Now (ii) easily follows from \f{\dc.3}.
In fact, unless $\,m=1$, the critical manifolds of $\,\vp\,$ cannot both
consist of single points, for if they did, we would have
$\,\vp_{\text{\rm max}}=\vp_{\text{\rm min}}=\y\,$ (by (ii), (iii) in \S\dc),
contradicting nonconstancy of $\,\vp\,$ in \f{\ir.1}. (Also, if both critical
manifolds were single points, with $\,m\ge2$, Reeb's theorem \cite{\mor} would
imply that $\,M\,$ is a topological $\,n$-sphere, $\,n\ge4$, admitting no
K\"ahler metric.)

By (i), the open set $\,M'\subset M\,$ on which $\,d\vp\ne0\,$ is the disjoint
union of the $\,\vp$-pre\-im\-ages of all values in
$\,(\vp_{\text{\rm min}},\hs\vp_{\text{\rm max}})$. Each of those
$\,\vp$-pre\-im\-ages is connected, due to the assertion of Corollary \a\mb.5
(which, as we saw, hold in our case); therefore, $\,Q=g(\navp,\navp)\,$ is
constant on it (Remark \a\sm.5(b)). This gives (iii), the
$\,C^\infty$-differentiability property of the assignment
$\,[\vp_{\text{\rm min}},\vp_{\text{\rm max}}]\to\bbR\,$ being now obvious
from the analogous local conclusion in Lemma \a\cg.1.

Lemma \a\cg.1 also gives $\,2\ta=\hs dQ/d\vp\,$ on the interval
$\,[\vp_{\text{\rm min}},\hs\vp_{\text{\rm max}}]$. The assertion about
$\,|\ta|\,$ in Corollary \a\pf.3 thus shows that
$\,|\hs dQ/d\vp|\,$ has the same positive value at both endpoints
$\,\vp_{\text{\rm min}},\hs\vp_{\text{\rm max}}$. Finally, since
$\,Q=g(\navp,\navp)$, the function $\,\vp\mapsto Q\,$ is positive on the open
interval $\,(\vp_{\text{\rm min}},\hs\vp_{\text{\rm max}})\,$ (formed by
non-critical values of $\,\vp$, cf. (i)), and vanishes at its endpoints
$\,\vp_{\text{\rm min}},\hs\vp_{\text{\rm max}}$. Hence $\,dQ/d\vp>0\,$ at
$\,\vp_{\text{\rm min}}$ and $\,dQ/d\vp<0\,$ at $\,\vp_{\text{\rm max}}$,
which yields (iv).

\head\S\db. The distance between the critical manifolds\endhead
For $\,[\vp_{\text{\rm min}},\vp_{\text{\rm max}}]\,$ with a function
$\,\vp\mapsto Q\,$ satisfying \f{\mw.1}, let us set
$$L\,=\,\int_{\vpsu_{\text{\rm min}}}^{\vpsu_{\text{\rm max}}}
{\,d\vp\over\sqrt{Q\,}}\,\in\,(0,\infty)\hs.\ff\db.1$$
To see that $\,L<\infty$, use $\,Q\,$ as the variable of integration near
either endpoint.
\proclaim{Lemma \a\db.1}Let\/ $\,N,N^*$ be the two critical manifolds of a
function\/ $\,\vp\,$ satisfying\/ \f{\ir.1} on a compact K\"ahler
manifold\/ $\,(M,g)$, and let\/ $\,L\hs$ be the invariant given by \f{\db.1},
with\/ $\,Q=g(\navp,\navp)\,$ treated as a function of\/ $\,\vp$, cf. {\rm
\S\bo}.
\widestnumber\item{(c)}\roster
\item"(a)"$L\,$ is the minimum distance between\/ $\,N\hs$ and any given
point\/ $\,y'\in N^*$.
\item"(b)"Every point\/ $\,x\in M\,$ at which\/ $\,d\vp\ne0\,$ can be joined
to\/ $\,N\hs$ by a geodesic, normal to\/ $\,N\nh$, of some length\/
$\,\ell\in(0,L)$.
\item"(c)"If\/ $\,X\subset M\,$ is a geodesic of length\/ $\,L\hs$ with
endpoints\/ $\,y,y'$ such that\/ $\,y\in N\,$ and\/ $\,X\hs$ is normal to\/
$\,N\hs$ at\/ $\,y$, then\/ $\,y'\in N^*$ and\/ $\,Q>0\,$ on\/
$\,X\smallsetminus\{y,y'\}$.
\endroster
\endproclaim
\demo{Proof}For $\,X,y,y'$ as in (c), let $\,X'$ be the maximal half-open
geodesic segment containing $\,y\,$ as an endpoint along with all points of
$\,X\hs$ sufficiently close to
$\,y\,$ and such that $\,d\vp\ne0\,$ everywhere in $\,X'\smallsetminus\{y\}$,
and let $\,[\hs0,\ell)\ni \sa\mapsto x(\sa)\,$ be an arc-length
parameterization of $\,X'$. By \f{\dc.4}, $\,d\sa=\pm\hs Q^{-1/2}\hs d\vp\,$
on $\,X'$, so that
$\,\ell=\int_0^\ell d\sa=|\int_{\vpsu_0}^{\vpsu'}Q^{-1/2}\hs d\vp\hs|$, where
$\,\vp_0=\vp(y)\in\{\vp_{\text{\rm min}},\vp_{\text{\rm max}}\}\,$ is the
value of $\,\vp\,$ on $\,N\nh$, and $\,\vp'=\vp(x(\ell))\,$ with
$\,x(\ell)=\,\lim_{\,\sa\to\ell}\,x(\sa)$. (Note that $\,\ell\le L<\infty\,$
by \f{\db.1}, and $\,\ell<L\,$ unless
$\,\vp'\in\{\vp_{\text{\rm min}},\vp_{\text{\rm max}}\}$.) However, maximality
of $\,X'$ now gives $\,(d\vp)(x(\ell))=0$, and so, as $\,\vp(x(\ell))=\vp'$,
(i) in \S\bo\ shows that
$\,\{\vp_0,\vp'\}=\{\vp_{\text{\rm min}},\vp_{\text{\rm max}}\}$, i.e.,
$\,\ell=L$. Consequently, $\,X'=X\hs$ and (c) follows.

Given $\,y'\in N^*$, let $\,y\,$ be the point of $\,N\hs$ nearest to $\,y'$,
and let $\,X'$ be a minimizing geodesic segment of some length $\,L'$, joining
$\,y'$ to $\,y$. As (c) implies that every point in a given critical manifold
lies at the distance $\,L\,$ from some point in the other critical manifold,
we have $\,L'\le L$. On the other hand, $\,L'\ge L$. In fact, if we had
$\,L'<L$, by extending $\,X'$ beyond $\,y'$ so as to obtain a geodesic segment
$\,X\hs$ of length $\,L\,$ we would conclude, from the final clause of (c),
that $\,y'$ is not a critical point of $\,\vp$. (Note that $\,X'$ is normal to
$\,N\hs$ at $\,y\,$ due to our distance-minimizing choice of $\,y\,$ and
$\,X'$.) Hence $\,L'=L$, which gives (a).

To prove (b), let us connect any $\,x\in M'=M\smallsetminus(N\cup N^*)\,$ with
the point $\,y\,$ nearest to it in $\,N\cup N^*$ by a minimizing geodesic
segment $\,X'$ of some length $\,\ell>0$. Thus, $\,\ell<L$, or else some point
of $\,X'$ would lie at the distance $\,L\,$ from $\,y$, and so, by (c), it
would be a point of $\,N\cup N^*$, closer to $\,x\,$ than $\,y\,$ is.
Extending $\,X'$ beyond $\,x$, we obtain a geodesic segment $\,X\hs$ of length
$\,L\,$ and, by (c), one of the endpoints of $\,X\hs$ lies in $\,N\nh$.
Moreover, $\,X\hs$ must be normal to $\,N\hs$ at that endpoint, since, by (a),
$\,X\hs$ is a minimum-length curve joining $\,N\hs$ to $\,N^*$. This completes
the proof.{\hfill\qd}
\enddemo

\head\S\nd. The normal horizontal distribution\endhead
The normal bundle $\,\Cal L\,$ of any submanifold $\,N\hs$ of a Riemannian
manifold $\,(M,g)\,$ carries the usual {\it normal connection\/}
$\,\nabla\nrm$, characterized by
$\,\nabla\nrm_{\!\dot y}w=[\nabla_{\!\dot y}w]\nrm$ whenever
$\,t\mapsto w(t)\in T_{y(t)}M\,$ is a $\,C^1$ vector field normal to $\,N\,$
along a $\,C^1$ curve $\,t\mapsto y(t)\in N\nh$. Here $\,\dot y=\hs dy/dt$,
while $\,\nabla\,$ is the Levi-Civita connection of $\,(M,g)$, and
$\,\,...\nrm\,$ stands for the component normal to $\,N\nh$.

Let $\,\Cal L\,$ now denote both the normal bundle of $\,N\nh$, and the total
space thereof (see \f{\cn.2}), where $\,N\hs$ is a critical manifold of a
function $\,\vp\,$ with \f{\ir.1} on a K\"ahler manifold $\,(M,g)$, cf. Remark
\a\cz.3(ii). By \f{\dc.3}, two cases are possible:
\widestnumber\item{(b)}\roster
\item"(a)"$N\,$ is a complex submanifold of complex codimension $\,1\,$ in
$\,M$, so that $\,\Cal L\,$ is a complex line bundle over $\,N\nh$, or
\item"(b)"$N=\{y\}\,$ for some point $\,y\in M$, and so
$\,\Cal L=\{y\}\times T_yM$.
\endroster
The {\it normal horizontal distribution\/} of $\,N\,$ is a distribution
$\,\hn$ on $\,\Cal L\smallsetminus N\nh$, defined as follows. In case
(a), $\,\hn$ is the restriction to $\,\Cal L\smallsetminus N\,$ of
the horizontal distribution of the normal connection in $\,\Cal L\,$
(see above), while, in case (b), $\,\hn$ coincides with the distribution
$\,\Cal H\,$ of Remark \a\tb.1 for $\,V=T_yM\,$ with the Hermitian inner
product $\,\langle\,,\rangle\,$ whose real part is $\,g(y)$, provided that one
identifies $\,\Cal L=\{y\}\times T_yM\,$ with $\,T_yM\,$ (as we will always
do). Note that $\,\hn$ is not only a real vector subbundle of the tangent
bundle $\,T(\Cal L\smallsetminus N)$, but also a {\it complex\/} vector
bundle, with the complex structure in each fibre $\,\hnyz=\hn_{(y,\zx)}$
inherited, in case (b), from the ambient space $\,T_yN\,$ (in which $\,\hnyz$
is contained as a complex subspace), or provided, in case (a), by requiring
the differential at $\,(y,\zx)\,$ of the bundle projection $\,\Cal L\to N\,$
to be a complex linear operator $\,\hnyz\to T_yN\nh$.

We also define vector fields $\,\vn,\un$ on $\,\Cal L\,$ to be $\,v,u\,$ in
\f{\cn.3} (in case (a)), or $\,v,u\,$ in Remark \a\tb.1 (in case (b), with
$\,V=T_yM\,$ and $\,\langle\,,\rangle\,$ as above), where $\,a\,$ is the
constant determined by $\,N\hs$ via \f{\dc.1}.
\rmk{Remark \a\nd.1}As noted in \S\tb, $\,\hn\hs$ is, also in case (b), the
horizontal distribution of a connection. Therefore, every vector in $\,\hn$ at
any given point of $\,\Cal L\smallsetminus N\,$ is tangent to a curve in
$\,\Cal L\smallsetminus N\,$ which is {\it horizontal}\hs, i.e.,
tangent to $\,\hn$ at every point.
\endrmk
\rmk{Remark \a\nd.2}Given a totally geodesic submanifold $\,N\hs$ of a
Riemannian manifold $\,(M,g)$, a point $\,y\in N\nh$, and vectors
$\,w,w'\in T_yN\nh$, let $\,\nabla,R\,$ be the Levi-Civita connection and
curvature tensor of $\,(M,g)$, and let a Riemannian metric $\,h\,$ on $\,N\hs$
be a constant multiple of its submanifold metric. Then, for any vector
$\,\xi\,$ tangent (or, normal) to $\,N\hs$ at $\,y$, the value
$\,R(w,w')\xi\,$ coincides with the one obtained by replacing $\,R\,$ with the
curvature tensor of $\,(N,h)\,$ (or, respectively, of the normal connection in
the normal bundle $\,\Cal L\,$ of $\,N$).

In fact, extending $\,w,w',\xi\,$ to $\,C^\infty$ vector fields on a
neighborhood $\,\,U\,$ of $\,y\,$ in $\,M\,$ tangent/normal to $\,N\hs$ along
$\,N\cap\hs U\nh$, we see that $\,\nabla_{\!w}\xi$, restricted to
$\,N\cap\hs U\nh$, is the covariant derivative relative to the Levi-Civita
connection of $\,(N,h)\,$ (or, respectively, the normal connection in
$\,\Cal L$), and our claim is obvious from \f{\cn.1}.
\endrmk

\head\S\cu. Critical manifolds and curvature\endhead
\proclaim{Lemma \a\cu.1}Let\/ $\,\vp\,$ be a special \krp\ on a K\"ahler
manifold\/ $\,(M,g)$, cf. \f{\ir.1}, and let\/
$\,v,u,\Cal V,\Cal H,Q,\si,\ta\,$ be given by\/ \f{\sm.1} -- \f{\sm.3}, so
that\/ $\,\si,\ta\,$ are\/ $\,C^\infty$ functions on the open set\/ $\,M'$ on
which\/ $\,d\vp\ne0$. For any two\/ $\,C^\infty$ vector fields\/ $\,w,w'$
defined on an open subset of\/ $\,M'$ and orthogonal to\/ $\,v\,$ and\/
$\,u\,$ at every point, we then have, with\/ $\,R\,$ denoting the curvature
tensor,
\widestnumber\item{(ii)}\roster
\item"(i)"$\hs Q\hs R(w,w')v\,=\,2(\si-\ta)\si\hs g(Jw,w')\hs u\hs.$
\item"(ii)"Either of\/ $\,g(R(w,v)w'\nh,v)\,$ and\/ $\,g(R(w,u)w'\nh,u)\,$
equals\/ $\,g(w,w')\,$ times a function which does not depend on the choice
of\/ $\,w\,$ and\/ $\,w'\nh$.
\endroster
\endproclaim
\demo{Proof}Let $\,...\vrt,\,\hs...\hrz$ be the $\,\Cal V\hs$ and $\,\Cal H\,$
components relative to the decomposition $\,TM'=\Cal H\oplus\Cal V$. Then
$\,\,Q\hs[\nabla_{\!w}w']\vrt=-\hs\si\hs[g(w,w')v+g(Jw,w')u]\,\,$ and
$$Q\hs[\hs w\hs,w'\hs]\vrt\,=\,-\hs2\si\,g(Jw,w')\hs u\hskip7pt\text{\rm for\
any\ local}\hskip4ptC^\infty\hskip3.5pt\text{\rm sections}\hskip4ptw,w'
\hskip3pt\text{\rm of}\hskip4pt\Cal H\hs.\ff\cu.1$$
(See \cite{\dml}, formula \f{\rr.1}.) Since, by \f{\sm.5.i}, $\,\si\,$ is
constant in the direction of $\,w$, \f{\sm.4} gives
$\,\nabla_{\!w}\nabla_{\!w'}v=\si\hs\nabla_{\!w}w'\nh$. Also, from \f{\sm.4},
$\,\nabla_{\![w,w']}v=\si\hs[w,w']\hrz+\hs\ta\hs[w,w']\vrt
=(\ta-\si)\hs[w,w']\vrt+\si\hs[w,w']$. As $\,\nabla\,$ is torsion-free, (i)
now easily follows from \f{\cn.1} and \f{\cu.1}. Next, writing
$\,\langle\,,\rangle\,$ for $\,g(\hskip3pt,\hskip2.2pt)\,$ we have
$\,\langle\nabla_{\!v}\nabla_{\!w}w',v\rangle
=d_v\langle\nabla_{\!w}w',v\rangle-\langle\nabla_{\!w}w',\nabla_{\!v}v\rangle
=(\si\ta-\hs d_v\si)\langle w,w'\rangle-\si\hs d_v\langle w,w'\rangle$, since
$\,\nabla_{\!v}v=\ta\hskip.4ptv\,$ (by \f{\sm.4} for $\,w=v$) and
$\,\langle\nabla_{\!w}w',v\rangle
=\langle[\nabla_{\!w}w']\vrt,v\rangle=-\hs\si\langle w,w'\rangle\,$ (cf. the
above formula for $\,\,Q\hs[\nabla_{\!w}w']\vrt$ and \f{\sm.2}), while
$\,\langle\nabla_{\!v}w',v\rangle=-\hs\langle w',\nabla_{\!v}v\rangle=0\,$
(as $\,\nabla_{\!v}v=\ta\hskip.4ptv$), and so
$\,\langle\nabla_{\!w}\nabla_{\!v}w',v\rangle
=-\hs\langle\nabla_{\!v}w',\nabla_{\!w}v\rangle
=-\hs\si\langle w,\nabla_{\!v}w'\rangle\,$ (since \f{\sm.4} gives
$\,\nabla_{\!w}v=\si\hskip.4ptw$). Next, the local flow of $\,v\,$ leaves
$\,\Cal H\,$ invariant (see \cite{\dml}, Remark \a\lc.3, discussion of
condition (a)), so that $\,[w,v]\,$ is a section of $\,\Cal H\,$ and
\f{\sm.4} gives $\,\nabla_{\!w}v=\si\hskip.4ptw$,
$\,\nabla_{\![w,v]}v=\si\hskip.4pt[w,v]$. Hence, as $\,\nabla\,$ is
torsion-free, $\,\langle\nabla_{\![w,v]}w',v\rangle
=-\hs\langle w',\nabla_{\![w,v]}v\rangle
=\si\langle w',[v,w]\rangle
=\si\langle w',\nabla_{\!v}w\rangle
-\si^2\langle w,w'\rangle$. Now \f{\cn.1} with
$\,d_v\langle w,w'\rangle=\langle\nabla_{\!v}w,w'\rangle
+\langle w,\nabla_{\!v}w'\rangle\,$ yields assertion (ii) for
$\,\langle R(w,v)w'\nh,v\rangle$. However,
$\,\langle R(w,u)w'\nh,u\rangle=\langle R(Jw,v)Jw'\nh,v\rangle\,$ (and so (ii)
for $\,\langle R(w,u)w'\nh,u\rangle\,$ follows from
$\,\langle Jw,Jw'\rangle=\langle w,w'\rangle$). Namely, in each tangent space,
the operator $\,w'\mapsto R(w,v)w'$ is complex-linear, i.e., commutes with
$\,J$, since relation $\,\nabla J=0\,$ means that $\,\nabla\,$
is a connection in $\,TM\,$ treated as a {\it complex\/} vector bundle. Due to
skew-adjoint\-ness of $\,J$, this gives
$\,\langle R(w,u)w'\nh,u\rangle=\langle R(w,u)w'\nh,Jv\rangle
=-\hs\langle R(w,u)Jw'\nh,v\rangle$, which in turn equals
$\,-\hs\langle R(Jw'\nh,v)w,u\rangle=-\hs\langle R(Jw'\nh,v)w,Jv\rangle
=\langle R(Jw'\nh,v)Jw,v\rangle=\langle R(Jw,v)Jw'\nh,v\rangle$. This
completes the proof.{\hfill\qd}
\enddemo
Suppose that $\,\vp\,$ satisfies \f{\ir.1} on a K\"ahler manifold $\,(M,g)\,$
and a given critical manifold $\,N\hs$ of $\,\vp\,$ is of complex codimension
$\,1\,$ in $\,M$, cf. \f{\dc.3}, while $\,\ve,a\,$ are the constants
described in Remark \a\sm.3 and \f{\dc.1}, and $\,\fe:M\to\bbR\,$ is given by
\f{\am.4} with $\,\y\,$ as in Lemma \a\sm.2. (The cases $\,\si=0\,$ and
$\,\si\ne0\,$ in Lemma \a\sm.2 correspond to $\,\ve=0\,$ and $\,\ve=\pm\hs1$.)

We then define a K\"ahler metric $\,h\,$ on the complex manifold $\,N\hs$ to
be $\,1/\!\fe_0$ times the restriction of $\,g\,$ to $\,TN\nh$, where
$\,\fe_0$ is the constant value of $\,\fe\,$ on $\,N\nh$. Note that
$\,\fe_0>0\,$ as $\,\vp\ne\y\,$ on $\,N\hs$ when $\,\ve=\pm\hs1$, cf. (i) in
\S\dc.
\proclaim{Lemma \a\cu.2}Let\/ $\,N\hs$ be a critical manifold of a special
\krp\/ $\,\vp\,$ on a K\"ahler manifold\/ $\,(M,g)\,$ such that\/
$\,\dimc M=m\ge2\,$ and\/ $\,\dimc N=m-1$, and let\/ $\,\Cal L\,$ be the
normal bundle of\/ $\,N\hs$. Then, for\/ $\,\ve,a,h\,$ described above, and
with the normal connection and curvature form as in\/ {\rm\S\nd} and\/ {\rm
Remark \a\cn.1},
\widestnumber\item{(b)}\roster
\item"(a)"The K\"ahler manifold\/ $\,(N,h)\,$ is Einstein unless\/ $\,m=2$.
\item"(b)"The curvature form of the normal connection in\/ $\,\Cal L\,$
equals\/ $\,-\hs2\ve a\,$ times the K\"ahler form of\/ $\,(N,h)$, defined as
in\/ \f{\pn.5}.
\endroster
\endproclaim
\demo{Proof}Let $\,v,u,\Cal H\,$ be as in \f{\sm.1}. For any $\,x\in M'\nh$,
i.e., a point $\,x\in M\,$ at which $\,v(x)\ne0$, we define symmetric bilinear
forms $\,h(x)\,$ and $\,\rih\nh(x)\,$ on the space $\,\Cal H_x$ by declaring
$\,h(x)\,$ to be the restriction of $\,g(x)/\!\fe(x)\,$ to $\,\Cal H_x$, with
$\,\fe\,$ as in \f{\am.4} (so that $\,\fe>0\,$ on $\,M'\nh$, as $\,\vp\ne\y\,$
on $\,M'$ by Lemma \a\sm.2), and letting $\,\rih\nh(x)\,$ assign to vectors
$\,w,w'\in\Cal H_x$ the value
$\,\rih(w,w')\,$ equal to $\,\sum_jg(R(w,e_j)w',e_j)$, where $\,R\,$ is the
curvature tensor of $\,(M,g)\,$ and the $\,e_j$ run through any
$\,g(x)$-or\-tho\-nor\-mal basis of $\,\Cal H_x$. Since such $\,e_j$ along
with $\,Q^{-1/2}v\,$ and $\,Q^{-1/2}u\,$ (at $\,x$) then form a
$\,g(x)$-or\-tho\-nor\-mal basis of $\,T_xM$, cf. \f{\sm.2}, our
$\,\rih\nh(x)\,$ and the Ricci tensor $\,\ri\hs(x)\,$ of $\,g\,$ at $\,x\,$
are related by  $\,\rih(w,w')=\,\ri\hs(w,w')-[g(R(w,v)w',v)+g(R(w,u)w',u)]/Q$,
for $\,w,w'\in\Cal H_x$, with $\,v,u,Q\,$ standing for their values at $\,x$.
As $\,\ri\,=\la\hs g\,$ on $\,\Cal H\,$ (see \f{\sm.3}), Lemma \a\cu.1(ii)
shows that $\,\rih\nh(x)\,$ equals a scalar times $\,h(x)$.

As $\,M'$ is dense in $\,M\,$ (Remark \a\kg.4), choosing a sequence of points
$\,x\in M'$ converging to any given $\,y\in N\,$ and using the limit relation
$\,\Cal H_x\to T_yN\,$ (see (iv) in \S\dc) along with Remark \a\nd.2, we see
that the the Ricci tensor of $\,h\,$ at $\,y\,$ is a scalar multiple of
$\,h(y)$, which proves (a).

Lemma \a\cu.1(i) and \f{\sm.1} give $\,R(w,w')\xi=i\hs\varOmega(w,w')\xi\,$
for any point $\,x\in M'$ and vectors $\,w,w'\in\Cal H_x$,
$\,\xi\in\Cal H_x^\perp$, where $\,\varOmega(w,w')\in\bbR\,$ equals
$\,2(\si-\ta)\si/Q\,$ at $\,x\,$ times
$\,g(Jw,w')$. Let a variable point $\,x\in M'$ now tend to any fixed
$\,y\in N\nh$. Relation $\,\Cal H_x\to T_yN\,$ (cf. (iv) in \S\dc) then gives
the same for $\,y\,$ and $\,T_yN\,$ instead of $\,x\,$ and $\,\Cal H_x$,
while $\,\si(x)\to0$, $\,\ta(x)\to a\,$ (see end of \S\cg) and, unless
$\,\si=0\,$ identically,  Lemma \a\sm.2 yields
$\,Q(x)/\si(x)\to2\hs[\vp(y)-\y\hs]=\ve\fe(y)$. Now (b) is immediate from
Remark \a\nd.2 along with the definitions of the curvature form and K\"ahler
form (see Remark \a\cn.1 and \f{\pn.5}), which completes the proof.{\hfill\qd}
\enddemo

\head\S\cv. Variations and partial covariant derivatives\endhead
Let $\,(\sa,t)\mapsto x(\sa,t)\in M\,$ be a fixed $\,C^\infty$ {\it
variation of curves\/} in a manifold $\,M$, that is, a $\,C^\infty$ mapping
with real variables $\,\sa,t\,$ ranging independently over some intervals. By
$\,(\sa,t)${\it-de\-pen\-dent\/} functions $\,\varphi\,$ or vector fields
$\,w\,$ we then mean assignments sending each $\,(\sa,t)\,$ to
$\,\varphi(\sa,t)\in\bbR\,$ or $\,w(\sa,t)\in T_{x(\sa,t)}M$.
Differentiability of such objects is also well-defined, since they clearly are
just sections of appropriate pull\-back bundles. In particular, the velocities
of the curves $\,\sa\mapsto x(\sa,t)\,$ and $\,t\mapsto x(\sa,t)$, with
$\,t\,$ or $\,\sa\,$ fixed, form $\,(\sa,t)$-de\-pen\-dent vector fields, here
denoted $\,x_\sa$ and $\,x_t$, and having, in any local coordinates, the
components $\,x_\sa^j=\,\partial x^j/\partial\sa\,$ and
$\,x_t^j=\,\partial x^j/\partial t$, where $\,x^j(\sa,t)\,$ are the components
of $\,x(\sa,t)$. Ordinary vector fields $\,u\,$ on $\,M\,$ or functions
$\,f:M\to\bbR\,$ are treated as $\,(\sa,t)$-de\-pen\-dent ones that assign
$\,u(x(\sa,t))\,$ or $\,f(x(\sa,t))\,$ to any $\,(\sa,t)$.

We use the subscript notation $\,\varphi_\sa,\hs\varphi_t$ for the
partial derivatives of $\,(\sa,t)$-de\-pen\-dent $\,C^1$ functions
$\,\varphi$, including ordinary $\,C^1$ functions on $\,M$. If, in addition,
there is a fixed connection $\,\nabla\,$ in the tangent bundle $\,TM$, we may
differentiate $\,(\sa,t)$-de\-pen\-dent $\,C^1$ vector fields $\,w\,$
covariantly with respect to either parameter $\,\sa\,$ or $\,t\,$ (i.e., along
the curves mentioned above), obtaining $\,(\sa,t)$-de\-pen\-dent fields
$\,w_\sa,\hs w_t$ equal to $\,\nabla_{\!\dot x}w\,$ for $\,\dot x=x_\sa$ (or,
$\,\dot x=x_t$), with the local-coordinate expressions
$\,w_\sa^j=\,\partial w^j/\partial \sa\,+\,\varGamma_{\!kl}^jx_\sa^kw^l\,$ and
$\,w_t^j=\,\partial w^j/\partial t\,+\,\varGamma_{\!kl}^jx_t^kw^l$. Here
$\,\varGamma_{\!kl}^j$ are the component functions of $\,\nabla\nh$, evaluated
at $\,x(\sa,t)$.

Applied to $\,x_\sa$ and $\,x_t$, this leads to the $\,(\sa,t)$-de\-pen\-dent
fields $\,x_{\sa\sa}=(x_\sa){}_\sa$, $\,x_{\sa t}=(x_\sa){}_t$, etc. Thus,
$\,x_{\sa\sa}=0\,$ identically if and only if all the curves
$\,\sa\mapsto x(\sa,t)\,$ are uniform-parameter geodesics. If $\,\nabla\,$ is
torsion-free, then $\,\varGamma_{\!kl}^j=\hs\varGamma_{\!lk}^j$, and so
\widestnumber\item{($*$)}\roster
\item"($*$)"$x_{\sa t}=x_{t\sa}$.
\endroster
Let us now assume that $\,\nabla\,$ is the Levi-Civita connection of a fixed
Riemannian metric $\,g\,$ on $\,M$, while $\,N\hs$ is a submanifold of $\,M\,$
and $\,t\mapsto\zy(t)\,$ is a $\,C^\infty$ unit vector field normal to $\,N\,$
along some given $\,C^\infty$ curve $\,t\mapsto y(t)\in N\nh$, where $\,t\,$
ranges over some interval. Let us set
$\,x(\sa,t)=\,\text{\rm Exp}\hs(y(t),\sa\hs\zy(t))\,$ for all $\,\sa\,$ in
some interval of the form $\,[\hs0,\ell\hs]\,$ with $\,\ell>0$, where
$\,\,\text{\rm Exp}:\hs\uexp\to M\,$ is defined as in Remark \a\gd.3. (Such
$\,\ell\,$ exists, i.e., $\,(y(t),\sa\hs\zy(t))\in\hs\uexp$ for all $\,\sa,t$,
provided that one replaces the original interval of $\,t\,$ with a suitable
subinterval.) Then
\widestnumber\item{(a)}\roster
\item"(a)"$\,|x_\sa|=1\,$ and $\,x_{\sa\sa}=0\,$ for all $\,\sa,t$,
\item"(b)"$\,x_\sa(0,t)=\zy(t)\,$ is unit and normal to $\,N\nh$, and
$\,x_t(0,t)=\dot y(t)\,$ is tangent to $\,N\nh$, with $\,\dot y=\hs dy/dt$,
\item"(c)"$x_{\sa t}=\nabla_{\!\dot y}\zy\,$ at $\,\sa=0\,$ and any $\,t$,
\item"(d)"$\langle x_\sa,x_t\rangle=0\,$ for all $\,\sa,t$, which is known as
{\it Gauss's Lemma\/} (cf. \S\cg),
\endroster
with $\,\langle\,,\rangle\,$ standing for $\,g(\hskip3pt,\hskip2.2pt)$. In
fact, the formula for $\,x(\sa,t)\,$ implies (a) -- (c). The Leibniz rule for
$\,(\sa,t)$-de\-pen\-dent functions such as $\,\langle x_\sa,x_t\rangle\,$
yields $\,\langle x_\sa,x_t\rangle\nh_\sa=\langle x_{\sa\sa},x_t\rangle
+\langle x_\sa,x_{t\sa}\rangle$, and, from ($*$),
$\,2\langle x_\sa,x_{t\sa}\rangle=2\langle x_\sa,x_{\sa t}\rangle
=\langle x_\sa,x_\sa\rangle_t$. Hence (a) gives
$\,\langle x_\sa,x_t\rangle\nh_\sa=0$, i.e., $\,\langle x_\sa,x_t\rangle\,$
does not depend on $\,\sa$, and (d) follows since, by (b),
$\,\langle x_\sa,x_t\rangle=0\,$ when $\,\sa=0$.

Still making all the assumptions listed in the paragraph following ($*$), let
us also suppose that $\,(M,g)\,$ is a K\"ahler manifold with a special \krp\
$\,\vp\,$ (cf. \f{\ir.1}), while $\,N\hs$ is a critical manifold of $\,\vp\,$
(Remark \a\cz.3(ii)), $\,v,u,Q,\si,\ta\,$ are given by \f{\sm.1} -- \f{\sm.3}
(so that $\,\si,\ta\,$ are $\,C^\infty$ functions on the open set $\,M'$
defined by $\,d\vp\ne0$), and $\,x(\sa,t)\in M'$ for all $\,\sa>0\,$ and all
$\,t$. Then, for all $\,\sa,t\,$ with $\,\sa>0$,
\widestnumber\item{(a)}\roster
\item"(e)"$|v|=|u|=Q^{1/2}$,
\item"(f)"$v=\pm\hs Q^{1/2}x_\sa$,
\item"(g)"$\langle u,x_t\rangle\nh_\sa
=\pm\hs2\langle u,x_t\rangle\hs\ta\hs Q^{-1/2}=2\hs\langle u,x_{\sa t}\rangle$,
\item"(h)"$Q_\sa=\pm\hs2\ta\hs Q^{1/2}\,$ and $\,\si_\sa
=\pm\hs2(\ta-\si)\hs\si\hs Q^{-1/2}$,
\endroster
with a specific fixed sign $\,\pm\hs$, namely, the sign of $\,a\,$ in
\f{\dc.1} - \f{\dc.2}. In fact, \f{\sm.2} implies (e), while Lemma \a\dc.1 and
(a) give $\,v=\pm\hs|v|x_\sa\,$ with the required sign $\,\pm\hs$, so that
(f) follows from (e). Also, $\,f_\sa=\langle x_\sa,\nabla\!f\rangle\,$ for
$\,f=Q\,$ and $\,f=\si$, and so \f{\sm.5.i} combined with (f), (a) yields (h).
Next, $\,\langle u,x_{t\sa}\rangle
=-\hs\langle u_t,x_\sa\rangle=\langle u_\sa,x_t\rangle$. Namely, the first
relation follows from the Leibniz rule and ($*$), as
$\,\langle u,x_\sa\rangle=0\,$ (by (f), since $\,\langle u,v\rangle=0$, cf.
\f{\sm.2}), while the second is clear from \sky\ of $\,\nabla u\,$ (\S\kg), as
$\,u_\sa=(\nabla u)x_\sa\,$ (see above). The Leibniz rule now yields
$\,\langle u,x_t\rangle\nh_\sa=\langle u_\sa,x_t\rangle
+\langle u,x_{t\sa}\rangle=2\langle u_\sa,x_t\rangle$. This in turn implies
both $\,\langle u,x_t\rangle\nh_\sa=2\langle u,x_{t\sa}\rangle\,$ and
$\,\langle u,x_t\rangle\nh_\sa=\pm\hs2\langle u,x_t\rangle\hs\ta\hs Q^{-1/2}$
(since \f{\sm.4} gives
$\,\nabla_{\!v}u=\nabla_{\!v}(Jv)=J\nabla_{\!v}v=\ta\hskip.4ptu$, so that, by
(e) -- (f), $\,u_\sa=\pm\hs\ta\hs Q^{-1/2}u$), which proves (g).

\head\S\fv. The normal exponential mapping\endhead
Suppose that $\,N\hs$ is a critical manifold of a special \krp\ on a K\"ahler
manifold $\,(M,g)\,$ of complex dimension $\,m\ge2\,$ (see \f{\ir.1} and
Remark \a\cz.3(ii)), while $\,\hn,\vn,\un\,$ (or, $\,\Cal H,v,u$) are the
complex vector bundles and vector fields defined in \S\nd\ (or, \f{\sm.1}).
Letting $\,\Cal L\,$ stand, again, for the total space of the normal bundle of
$\,N\nh$, we also assume that the image of the set
$\,\,U'=(\hs\uexp\cap\Cal L)\smallsetminus N\,$ under the normal exponential
mapping of $\,N\,$ (see \S\cg) is contained in the open set $\,M'\subset M\,$
on which $\,d\vp\ne0$. Recall that, by Lemma \a\sm.2, $\,\si\,$ in \f{\sm.3}
is either identically zero, or nonzero everywhere in $\,M'$. If $\,\si\ne0\,$
on $\,M'$, there is a constant $\,\y\,$ such that $\,Q/\si=2(\vp-\y)\,$ on
$\,M'$ and, by (i) in \S\dc, in case (a) of \S\nd\ we have $\,\vp\ne\y\,$
everywhere in $\,N\nh$, while, in case (b) of \S\nd, $\,N=\{y\}\,$ and
$\,\vp(y)=\y$.

In addition, let $\,\varTheta:T_{(y,\zx)}\Cal L\to T_xM\,$ denote the
differential of $\,\text{\rm Exp}\,\hs$ at some fixed $\,(y,\zx)\in\hs U'\nh$,
with $\,x=\hs\text{\rm Exp}\hs(y,\zx)$, and let the symbol $\,\hnyz$ stand for
$\,\hn_{(y,\zx)}$.
\proclaim{Lemma \a\fv.1}Under these assumptions,
$\,\varTheta(\hnyz)\,\subset\,\Cal H_x$ and\/ $\,\varTheta:\hnyz\to\Cal H_x$
is complex-linear. Also, letting\/ $\,w,w'\in T_yN\,$ be the images of any
given\/ $\,\xi,\xi'\in\hnyz\,$ under the differential at\/ $\,(y,\zx)\,$ of
the bundle projection\/ $\,\Cal L\to N\nh$, we have
\widestnumber\item{(iii)}\roster
\item"(i)"$g(\varTheta\xi,\varTheta\xi')=g(w,w')\,$ if\/ $\,\si=0\,$
identically on\/ $\,M'\nh$.
\item"(ii)"$\lx\vp(y)-\y\hs\rx\,g(\varTheta\xi,\varTheta\xi')
=\,\lx\vp(x)-\y\hs\rx\,g(w,w')\,$ if\/ $\,\si\ne0\,$ on\/ $\,M'$ and case\/
{\rm(a)} of\/ {\rm\S\nd} occurs.
\item"(iii)"$a\hs g(\zx,\zx)\,g(\varTheta\xi,\varTheta\xi')
=2\,\lx\vp(x)-\y\hs\rx\,g(\xi,\xi')\,$
with\/ $\,a\,$ as in\/ \f{\dc.1}, if\/ $\,\si\ne0\,$ on\/ $\,M'$ and case\/
{\rm(b)} of\/ {\rm\S\nd} occurs.
\endroster
Finally, the\/ $\,\varTheta$-images of\/ $\,\vn(y,\zx)\,$ and\/
$\,\un(y,\zx)\,$ are\/ $\,|a\zx|\hs v(x)/|v(x)|\,$ and, respectively,
$\,u(x)$, with\/ $\,a\,$ as in\/ \f{\dc.1}.
\endproclaim
\demo{Proof}Let $\,y(t),\zy(t),x(\sa,t)\,$ be as in the paragraph following
($*$) in \S\cv, and in addition such that, in case (a) of \S\nd,
the unit vector field $\,t\mapsto\zy(t)\,$ normal to $\,N\,$ along the curve
$\,t\mapsto y(t)\in N\,$ is parallel relative to the Levi-Civita connection of
$\,(M,g)$, while, in case (b) of \S\nd, $\,y(t)=y\,$ for all $\,t\,$ and
$\,\dot\zy=\hs d\zy/dt\in T_yM\,$ is $\,g(y)$-or\-thog\-o\-nal to $\,\zy(t)\,$
and $\,J\zy(t)\,$ for every $\,t$. These assumptions mean that, for any fixed
$\,\sa$, the curve $\,t\mapsto(y(t),s\zy(t))\,$ in $\,\,U'$ is horizontal in
the sense of Remark \a\nd.1. (In case
(a) we assume that $\,\nabla_{\!\dot y}\zy=0$, rather than just
$\,[\nabla_{\!\dot y}\zy]\nrm=0\,$ as required by the definition of the normal
connection in \S\nd, since $\,N\hs$ is totally geodesic, cf. Remark
\a\cz.3(ii), and so $\,\nabla_{\!\dot y}\zy\,$ is normal to $\,N\hs$ whenever
$\,\zy\,$ is.)

Writing $\,\langle\,,\rangle\,$ for $\,g(\hskip3pt,\hskip2.2pt)\,$ we have
$\,\langle v,x_t\rangle=\langle u,x_t\rangle=0\,$ for all $\,\sa,t\,$
(notation of \S\cv). First, $\,\langle v,x_t\rangle=0\,$ by (f), (d) in \S\cv.
Next, (g) -- (h) in \S\cv\ yield $\,[\langle u,x_t\rangle/Q]_\sa=0$,
i.e., $\,\langle u,x_t\rangle/Q\,$ is constant as a function of $\,\sa\,$
with fixed $\,t$. To see that its constant value is $\,0\,$ we may evaluate
its limit as $\,\sa\to0\,$ using l'Hospital's rule and noting that, by (g),
(h) in \S\cv, $\,\langle u,x_t\rangle\nh_\sa/Q_\sa
=\pm\hs\langle u,x_{\sa t}\rangle\hs Q^{-1/2}/\ta$. In case (a),
$\,\langle u,x_t\rangle\nh_\sa/Q_\sa
=\pm\hs\langle u/|u|,x_{\sa t}\rangle/\ta\,\to\,0\,$ as $\,\sa=0$, by (c), (e)
in \S\cv, since $\,\nabla_{\!\dot y}\zy=0$, while $\,\ta=a\ne0\,$ on $\,N\hs$
(see end of \S\cg). In case (b),
$\,\langle u,x_t\rangle\nh_\sa/Q_\sa=\langle Jx_\sa,x_{\sa t}\rangle/\ta
\to\langle J\zy(t),\dot\zy\rangle/a=0\,$ as $\,\sa\to0\,$ by (f), (b) in
\S\cv\ with $\,u=Jv$, $\,\ta=a\ne0\,$ on $\,N\nh$, and our orthogonality
assumption for case (b). Thus, $\,\langle u,x_t\rangle=0$.

Relations $\,\langle v,x_t\rangle=\langle u,x_t\rangle=0\,$ state that
$\,v,u\,$ are $\,g$-normal to the $\,\text{\rm Exp}\hs$-image of every
horizontal curve in $\,\,U'\nh$. Hence $\,\varTheta(\hnyz)\subset\Cal H_x$
(cf. Remark \a\nd.1 and \f{\sm.1}).

To prove (i) -- (iii) we may assume, due to symmetry of $\,g$, that
$\,\xi=\xi'$. By \f{\pn.2} and \f{\sm.3}, $\,v_t=\si\hskip.4ptx_t$ for every
$\,\sa,t$, since $\,x_t(\sa,t)\in\Cal H_{x(\sa,t)}$. (Note that $\,v_t$ equals
$\,(\nabla v)x_t$, the covariant derivative of $\,v\,$ in the direction of
$\,x_t$, cf. \S\cv\ and \f{\pn.1}.) Also,
$\,Q_t=\langle x_t,\nabla Q\rangle=0$, as $\,\nabla Q=2\ta\hskip.4ptv\,$
(see \f{\sm.5.i}), while $\,\langle v,x_t\rangle=0$, and so, since
$\,x_\sa=\pm\hs Q^{-1/2}v\,$ by (f) in \S\cv, we have
$\,x_{\sa t}=\pm\hs Q^{-1/2}\si\hskip.4ptx_t$. The Leibniz rule and ($*$) in
\S\cv\ now imply
$\,\langle x_t,x_t\rangle\nh_\sa=2\langle x_t,x_{\sa t}\rangle
=\pm\hs2\langle x_t,x_t\rangle\hs\si\hskip.4ptQ^{-1/2}$. This, along with (b)
in \S\cv, yields $\,\langle x_t,x_t\rangle=g(\dot y,\dot y)\,$ when $\,\si=0$,
thus proving (i); at the same time, combined with (h) in \S\cv, it gives
$\,[\langle x_t,x_t\rangle\hs\si/Q]_\sa=0$, so that
$\,\langle x_t,x_t\rangle\hs\si/Q\,$ is constant as a function of $\,\sa$.
When $\,\si\ne0\,$ on $\,M'\nh$, we have $\,Q/\si=2(\vp-\y)\,$ (by Lemma
\a\sm.2), so that $\,\langle x_t,x_t\rangle/(\vp-\y)\,$ is constant as a
function of $\,\sa$, and we can find its constant value by evaluating its
limit as $\sa\to0$. Specifically, in case (a) of \S\nd, $\,\vp(y)\ne\y\,$
(see (i) in \S\dc), and (ii) follows as
$\,\langle x_t,x_t\rangle/(\vp-\y)\,$ at $\,x=x(\sa,t)\,$ equals
$\,g(\dot y,\dot y)/(\vp-\y)\,$ at $\,y=x(0,t)$, cf. (b) in \S\cv. In case (b)
of \S\nd, however, we find the limit by applying l'Hospital's rule {\it
twice}\hs, as $\,\vp(y)=\y\,$ (cf. (i) in \S\dc) and $\,d\vp=0\,$ at $\,y$,
while $\,x_t=0\,$ at $\,\sa=0\,$ by (b) in \S\cv; this gives
$\,2\langle x_{t\sa},x_{t\sa}\rangle\,$ in the numerator (at $\,\sa=0$) and
$\,\vp_{\sa\sa}$ in the denominator. By Remark \a\sm.4(ii) and
($*$), (c) in \S\cv, $\,\vp_{\sa\sa}=\ta(y)=a\,$ and $\,x_{t\sa}$ equals
$\,\dot\zy=\hs d\zy/dt\,$ at $\,\sa=0$. (Note that $\,\ta=a\,$ on $\,N\nh$,
cf. end of \S\cg, while $\,\nabla_{\!\dot y}\zy=\dot\zy\,$ as $\,y(t)=y\,$ is
constant.) Now (iii) follows: $\,\langle x_t,x_t\rangle/(\vp-\y)\,$ at any
$\,x=x(\sa,t)\,$ is the same as at $\,y=x(0,t)$, i.e.,
$\,\langle x_t,x_t\rangle/(\vp-\y)=2g(\dot\zy,\dot\zy)/a
=2g(\xi,\xi)/(a\sa^2)=2g(\xi,\xi)/\lx a\hs g(\zx,\zx)\rx$, where
$\,\zx=\sa\hs\zy(t)\,$ and $\,\xi=\sa\hs\dot\zy(t)$.

The equality $\,x_{\sa t}=\pm\hs Q^{-1/2}\si\hskip.4ptx_t$ established in the
last paragraph means that
$$\nabla_{\!\dot x}w\,=\,\pm\hs Q^{-1/2}\si\hskip.4ptw\,,\qquad\text{\rm with}
\quad\dot x=x_\sa\,,\ff\fv.1$$
where $\,w=x_t$ stands for the vector field $\,\sa\mapsto w(\sa)=x_t(\sa,t)\,$
along the geodesic $\,\sa\mapsto x(\sa,t)$, for any fixed $\,t$. As
$\,\nabla J=0$, \f{\fv.1} holds for $\,\tilde w=Jw\,$ whenever it does
for $\,w$. In case (a) of \S\nd, the $\,\text{\rm Exp}\hs$-pre\-im\-age
of $\,w\,$ is a vector field along the curve
$\,\sa\mapsto(y(t),\sa\hs\zy(t))\in\Cal L\,$ which arises as the horizontal
lift of $\,w(0)$. (In fact, at any $\,\sa,t$, the pre\-im\-age is the velocity
vector of the curve $\,t\mapsto(y(t),\sa\hs\zy(t))$, which we chose to be
horizontal, and which has the projection image $\,t\mapsto y(t)\,$ with the
velocity $\,w(0)$, cf. (b) in \S\cv.) Replacing $\,w(0)\,$ by $\,Jw(0)\,$
causes such a horizontal-lift field to become multiplied by $\,i\,$ in the
complex vector bundle $\,\hn$ (see \S\nd), and at the same time results in
replacing $\,w\,$ with \f{\fv.1} by $\,\tilde w=Jw$, since a solution $\,w\,$
to \f{\fv.1} is determined by the initial value $\,w(0)$. This proves our
complex-linearity assertion in case (a) of \S\nd. In the remaining case (b),
under the identification $\,\Cal L=T_yM\,$ (\S\nd), the
$\,\text{\rm Exp}\hs$-pre\-im\-age of $\,w\,$ is the vector field
$\,\sa\mapsto\sa\hs\dot\zy(t)\,$ along the line segment
$\,\sa\mapsto\sa\hs\zy(t)\in T_yM\,$ (with fixed $\,t$). Therefore,
$\,w(0)=0\,$ and $\,w(\sa)/\sa\,$ has a limit as $\,\sa\to0$, equal (in view
of the local-coordinate formula for $\,\nabla_{\!\dot x}w$) to the value of
$\,\nabla_{\!\dot x}w\,$ at $\,\sa=0$. As the
$\,\text{\rm Exp}\hs$-pre\-im\-age of the limit is $\,\dot\zy(t)$, it follows
that, in case (b) of \S\nd, a solution $\,w\,$ to \f{\fv.1} is uniquely
determined by $\,(\nabla_{\!\dot x}w)(0)=\dot\zy(t)\,$ (with fixed $\,t$).
Thus, replacing $\,w\,$ by $\,\tilde w=Jw\,$ amounts to using
$\,J\dot\zy(t)\,$ instead of $\,\dot\zy(t)$, i.e., to multiplying the
$\,\text{\rm Exp}\hs$-pre\-im\-age of $\,w\,$ by $\,i\,$ in the complex vector
bundle $\,\hn$ (cf. \S\nd), which establishes our complex-linearity claim also
in case (b).

Finally, relation $\,\dot x=(\sgn\,a)\hs v/|v|\,$ in Lemma \a\dc.1, for
$\,x(\sa)=\,\text{\rm Exp}\hs(y,\sa\hs\zx/|\zx|)\,$ at $\,\sa=|\zx|$, shows
that $\,\varTheta\,$ sends $\,\zx/|\zx|$, treated as a vertical vector in
$\,T_{(y,\zx)}\Cal L\,$ (cf. \S\cn), onto $\,(\sgn\,a)\hs v(x)/|v(x)|$.
Multiplying both vectors by $\,a|\zx|$, we obtain our assertion about the
$\,\varTheta$-image of $\,\vn(y,\zx)$. Also, since $\,\zx\in (T_yN)^\perp$,
\f{\dc.2} and \f{\pf.2} -- \f{\pf.3} give $\,\nabla_{\!\zx}u=i\hs a\zx$. Now
$\,\varTheta(\un(y,\zx))=u(x)\,$ by Lemma \a\gd.5, since on the normal space
$\,\Cal L_y\subset T_yM\,$ (identified, as usual, with
$\,\{y\}\times\Cal L_y\subset\Cal L$), the normal exponential mapping of
$\,N\hs$ coincides with $\,\,\e_{\hskip.4pty}$. This completes the
proof.{\hfill\qd}
\enddemo

\head\S\nx. A normal exponential diffeomorphism\endhead
The normal exponential mapping of $\,N\nh$, used below, was introduced in
\S\cg; this time it is defined on the whole total space of the normal bundle
of $\,N\nh$, since, due to compactness of $\,M$, we have $\,\,\uexp=TM\,$
(notation of Remark \a\gd.3).
\proclaim{Lemma \a\nx.1}Suppose that\/ $\,L\hs$ is given by\/ \f{\db.1} for\/
$\,\vp_{\text{\rm min}},\vp_{\text{\rm max}}$ and\/ $\,\vp\mapsto Q\,$
determined as in\/ {\rm(iii)} of\/ {\rm\S\bo} by a function\/ $\,\vp\,$
satisfying\/ \f{\ir.1} on a compact K\"ahler manifold\/ $\,(M,g)$, and\/
$\,N\hs$ is a critical manifold of\/ $\,\vp$, cf. {\rm Remark \a\cz.3(ii)}.
If\/ $\,\Cal L\,$ denotes the total space of the normal bundle of\/ $\,N\nh$,
while\/ $\,M'\subset M\,$ and\/ $\,\Cal L'\subset\Cal L\,$ are the open sets
defined by\/ $\,d\vp\ne0\,$ and, respectively, $\,0<\sa<L$, where\/ $\,\sa\,$
is the norm function\/ $\,\Cal L\to[\hs0,\infty)\,$ corresponding as in\/ {\rm
Remark \a\cn.2} to the fibre metric obtained by restricting\/ $\,g\,$ to\/
$\,\Cal L\hs$, then the restriction to\/ $\,\Cal L'$ of the normal exponential
mapping of\/ $\,N\hs$ is a\/ $\,C^\infty$ diffeomorphism\/
$\,\,\text{\rm Exp}:\Cal L'\to M'$.
\endproclaim
\demo{Proof}The $\,\,\text{\rm Exp}\hs$-image of any open line segment of
length $\,L\,$ emanating from $\,0\,$ in any fibre $\,\Cal L'_y$ of the
punctured-disk bundle $\,\Cal L'$ has the form $\,X\smallsetminus\{y,y'\}$,
with $\,X\hs$ and $\,y,y'$ as in Lemma \a\db.1(c), so that, by Lemma
\a\db.1(c), $\,X\smallsetminus\{y,y'\}\subset M'$. Hence
$\,\,\text{\rm Exp}\,\,$ actually sends $\,\Cal L'$ into $\,M'\nh$.

Surjectivity of $\,\,\text{\rm Exp}:\Cal L'\to M'$ is obvious from Lemma
\a\db.1(b). To prove its injectivity, suppose that $\,(y,\zx)\in\Cal L'$ and
$\,x=\hs\text{\rm Exp}\hs(y,\zx)\in M'$. Since $\,0<|\zx|<L$, we can express
$\,(y,\zx)\,$ in terms of $\,x\,$ by travelling backwards along the unit-speed
geodesic $\,t\mapsto x(t)=\,\text{\rm Exp}\hs(y,t\hs\zx/|\zx|)$, which has
$\,x(0)=y$, $\,\dot x(0)=\zx/|\zx|$, $\,x(\sa)=x\,$ (where
$\,\dot x=\hs dx/dt\,$ and $\,\sa=|\zx|\in(0,L)$) and, by Lemma \a\dc.1,
$\,\dot x(\sa)=w(x)\,$ for the vector field $\,w=(\hs\sgn\,a)\hs v/|v|\,$ on
$\,M'$ (with $\,v,a\,$ defined as in Lemma \a\dc.1). In fact, the
re-parameterized geodesic
$\,t\mapsto y(t)=\,\text{\rm Exp}\hs(x,-\hs tw(x))\,$ clearly has $\,y(0)=x$,
$\,\dot y(0)=-\hs w(x)$, $\,y(\sa)=y\,$ and $\,\dot y(\sa)=-\hs\zx/|\zx|$, so
that $\,(y,\zx)=(y(\sa),-\hs\sa\hs\dot y(\sa))$. Moreover, $\,\sa\,$ is
uniquely determined by $\,x\,$ and depends $\,C^\infty$-differentiably on
$\,x\,$ (via $\,\vp(x)$), since the assignment $\,\sa\mapsto\vp\,$ with
\f{\dc.4} is a $\,C^\infty$ diffeomorphism
$\,(0,L)\to(\vp_{\text{\rm min}},\vp_{\text{\rm max}})\,$ (cf. \f{\db.1}). The
last formula for $\,(y,\zx)\,$ thus shows that $\,(y,\zx)\,$ is determined by
$\,x$, i.e., $\,\,\text{\rm Exp}:\Cal L'\to M'$ is injective, and its inverse
$\,M'\to\Cal L'$ is of class $\,C^\infty$. This completes the proof.{\hfill\qd}
\enddemo
\proclaim{Lemma \a\nx.2}Let\/ $\,\Cal L\,$ be the total space, with\/
\f{\cn.2}, of a real/com\-plex vector bundle with a
Riem\-ann\-i\-an/Her\-mit\-i\-an fibre metric\/ $\,\langle\,,\rangle\,$ over a
manifold\/ $\,N\nh$, and let\/ $\,\Cal L'\subset\Cal L\,$ denote the open set
given by\/ $\,0<\sa<L$, for some\/ $\,L\in(0,\infty)$, where\/ $\,\sa\,$ is
the norm function\/ $\,\Cal L\to[\hs0,\infty)$, cf. {\rm Remark \a\cn.2}.
Then, setting\/ $\,\varPhi(y,\zx)=(y,\sa\zx/|\zx|)\,$ for\/ $\,y\in N\hs$
and\/ $\,\zx\in\Cal L_y\smallsetminus\{0\}$, where\/ $\,\sa\in(0,L)\,$ depends
on\/ $\,r=|\zx|\in(0,\infty)\,$ via a fixed\/ $\,C^\infty$ diffeomorphism\/
$\,(0,\infty)\to(0,L)$, we obtain a\/ $\,C^\infty$ diffeomorphism\/
$\,\varPhi:\Cal L\smallsetminus N\to\Cal L'$.
\endproclaim
In fact, if $\,(y,w)=\varPhi(y,\zx)\,$ and $\,r=|\zx|>0$, then
$\,w=\sa\zx/r\hs$ and $\,\zx=rw/\sa\,$ with $\,r\hs$ obtained from
$\,\sa=|w|\,$ via the inverse diffeomorphism
$\,(0,L)\to(0,\infty)$.{\hfill\qd}
\rmk{Remark \a\nx.3}With the notations and under the assumptions of Lemma
\a\nx.2, let $\,\Cal H\,$ be the horizontal distribution of a fixed connection
in $\,\Cal L\,$ making $\,\langle\,,\rangle\,$ parallel, and let
$\,\varXi:T_{(y,\zx)}\Cal L\to T_{\varPhi(y,\zx)}\Cal L\,$ denote the
differential of $\,\varPhi\,$ at any given point
$\,(y,\zx)\in\Cal L\smallsetminus N\nh$. Then, for $\,\Cal H'=\Cal H_{(y,\zx)}$
and $\,\Cal H''=\Cal H_{\varPhi(y,\zx)}$,
\widestnumber\item{(c)}\roster
\item"(a)"$\varXi\,$ maps $\,\Cal H'$ onto $\,\Cal H''$.
\item"(b)"In the case where $\,\Cal L\,$ is a {\it complex\/} vector bundle,
$\,\varXi\,$ sends the vectors $\,v(y,\zx)\,$ and $\,u(y,\zx)$, defined by
\f{\cn.3} for any fixed $\,a\in\bbR\smallsetminus\{0\}$, onto $\,r/\sa\,$
times $\,d\sa/dr\,$ times $\,v(\varPhi(y,\zx))\,$ and, respectively, onto
$\,u(\varPhi(y,\zx))$, with $\,\sa\,$ and $\,d\sa/dr\,$ evaluated at
$\,r=|\zx|$.
\endroster
In fact, (a) is immediate: the norm function $\,r\hs$ is constant along any
horizontal curve in $\,\Cal L\hs$, and so $\,\varPhi\,$ multiplies such a curve
by a constant factor. Also, (b) for $\,u\,$ is obvious as the flow of $\,u\,$
consists of multiplications by scalars of modulus one, each of which commutes
with $\,\varPhi\,$ (treated as a mapping
$\,\Cal L\smallsetminus N\to\Cal L\smallsetminus N$). Finally, (c) for $\,v\,$
follows since an integral curve $\,t\mapsto(y,e^{at}\zx)\,$ of $\,v$, for a
fixed $\,\zx\in\Cal L_y$ with $\,|\zx|=1$, has the $\,\varPhi$-image
$\,t\mapsto(y,w(t))\,$ such that $\,w(t)=\sa\zx$, with $\,\sa\,$ depending on
$\,t\,$ via the given diffeomorphism $\,r\mapsto\sa$, where $\,r=e^{at}$.
Hence $\,dw/dt\,$ equals $\,r/\sa\,$ times $\,d\sa/dr\,$ times $\,aw(t)$, as
required.
\endrmk

\head\S\gc. A global classification of special K\"ahler-Ricci
potentials\endhead
We will now show that every triple $\,M,g,\vp\,$ formed by a special \krp\
$\,\vp\,$ on a compact K\"ahler manifold $\,(M,g)\,$ is biholomorphically
isometric to one of the examples constructed in \S\mw\ and \S\cp. Since we
already know that, conversely, each of those examples constitutes a compact
K\"ahler manifold with a a special \krp, the result of this section amounts to
a complete classification theorem for such triples $\,M,g,\vp$.

Note that, in contrast with its use elsewhere in the text, the symbol
$\,(\srf,\gx)\,$ stands in this section for a Riemannian (or K\"ahler)
manifold {\it of any dimension}\hs.
\proclaim{Lemma \a\gc.1}Let\/ $\,(\srf,\gx)$, $\,(M,g)\,$ be complete
Riemannian manifolds with open subsets\/ $\,\srf'\subset\srf$,
$\,M'\subset M\,$ such that both\/ $\,\srf\smallsetminus\srf'$ and\/
$\,M\smallsetminus M'$ are unions of finitely many compact submanifolds of
codimensions greater than one. Any isometry\/ $\,\varPsi\,$
of\/ $\,(\srf'\!,\gx)\,$ onto\/ $\,(M'\!,g)\,$ then can be uniquely extended
to an isometry of\/ $\,(\srf,\gx)\,$ onto\/ $\,(M,g)$. If, in addition,
$\,(\srf,\gx)\,$ and\/ $\,(M,g)\,$ are K\"ahler manifolds and the isometry\/
$\,\varPsi:\srf'\to M'$ is a biholomorphism, then so is the extension\/
$\,\srf\to M$.
\endproclaim
In fact, by the codimension hypothesis $\,\srf'$ (or, $\,M'$) is connected and
dense in $\,\srf\,$ (or, in $\,M$), and the inclusion mappings
$\,\srf'\to\srf$, $\,M'\to M\,$ are distance-preserving. Thus, as metric
spaces, $\,\srf\,$ and $\,M\,$ are the completions of $\,\srf'$ and $\,M'$.
Our assertion now follows since distance-preserving mappings are isometries
\cite{\kno}, p. 169, the K\"ahler case being obvious from a continuity
argument.{\hfill\qd}
\proclaim{Theorem \a\gc.2}Let\/ $\,\vp\,$ be a special \krp\ on a compact
K\"ahler manifold\/ $\,(M,g)\,$ of complex dimension\/ $\,m\ge2$, cf.
\f{\ir.1}. Then either
\widestnumber\item{(iii)}\roster
\item"(i)"$M,g,\vp\,$ are, up to a biholomorphic isometry, obtained as in\/
{\rm\S\mw;}\hskip4ptor,
\item"(ii)"$M\,$ can be biholomorphically identified with\/ $\,\bbCP^m$ in
such a way that\/ $\,g,\vp\,$ become the objects constructed in\/ {\rm\S\cp}.
\endroster
\endproclaim
\demo{Proof}We denote $\,N,N^*$ the two critical manifolds of $\,\vp$, ordered
so that either
\widestnumber\item{(b)}\roster
\item"(a)"both $N,N^*$ are of complex codimension one, or
\item"(b)"$N=\{y\}\,$ for some $\,y\in M\,$ (cf. (i), (ii) in \S\bo).
\endroster
In case (a) (or, (b)), we define the data \f{\am.3} (or, \f{\tm.1}) as
follows. First, in both cases, $\,m\,$ is the complex dimension of $\,M$,
while $\,\iyp=(\vp_{\text{\rm min}},\vp_{\text{\rm max}})$, the variable
$\,\vp\in\iyp$ and the function $\,Q\,$ of $\,\vp\,$ are chosen to be the ones
appearing in the assignment
$\,[\vp_{\text{\rm min}},\vp_{\text{\rm max}}]\ni\vp\mapsto Q\,$ determined by
our $\,M,g,\vp\,$ as in (iii) of \S\bo; $\,r\hs$ is a fixed
function with \f{\am.1}, $\,\ve\,$ is the invariant defined in Remark \a\sm.3,
$\,a\in\bbR\,$ satisfies \f{\dc.1} with our $\,N\nh$, and $\,\y\,$ is the
constant with $\,Q/\si=2(\vp-\y)\,$ on $\,M'$ (notation of Lemma \a\sm.2).
Thus, $\,\y\,$ remains undefined when $\,\si=0\,$ identically on $\,M'$, cf.
Lemma \a\sm.2, which, by (ii) in \S\dc, may happen only in case (a).

Next, in case (a), $\,N\hs$ is the critical manifold chosen above,
$\,\Cal L\,$ is its normal bundle, $\,\Cal H\,$ in \f{\am.3} is the horizontal
distribution of the normal connection in $\,\Cal L\,$ (see \S\nd),  while
$\,\langle\,,\rangle\,$ is the Hermitian fibre metric in $\,\Cal L\,$ whose
real part is the restriction of $\,g\,$ to $\,\Cal L\hs$, and $\,h\,$ is the
metric on $\,N\hs$ defined in the paragraph preceding Lemma \a\cu.2. On the
other hand, in case (b), $\,V=T_yM\,$ and $\,\langle\,,\rangle\,$ is the
Hermitian inner product with the real part $\,g(y)$.

The data \f{\am.3} (or, \f{\tm.1}) just defined in case (a) (or, (b)) satisfy
the conditions listed in the paragraph following \f{\am.3} (or, \f{\tm.1}; see
Lemma \a\cu.2 and Remark \a\ob.2). They also fulfill the additional
requirements in \S\mw\ (or, \S\cp): namely, (iv) in \S\bo\ gives \f{\mw.1} in
both cases, while, by (iii) in \S\dc,
$\,\y\notin[\vp_{\text{\rm min}},\vp_{\text{\rm max}}]\,$ in case (a) and
$\,\vp(y)=\y\,$ in case (b). The construction of \S\mw\ (or, \S\cp) now yields
a compact K\"ahler manifold of complex dimension $\,m$, which we will denote
$\,(\srf,\gx)$, rather than $\,(M,g)$, and a special \krp\ on $\,(\srf,\gx)$,
still denoted $\,\vp$.

We define the set $\,\srf'$ to be $\,\Cal L\smallsetminus N\,$ in case (a) and
$\,T_yM\smallsetminus\{0\}\,$ in case (b), so that $\,\srf'$ may be treated as
an open subset of $\,\srf\,$ (cf. \S\mw\ or \S\cp). Then, in both cases,
$\,\srf'=\Cal L\smallsetminus N\nh$, since $\,\Cal L\,$ is, also in case (b),
the total space of the normal bundle of $\,N\nh$, provided that we identify
$\,\{y\}\times T_yM\,$ with $\,T_yM\,$ as in \S\nd.

The inverse of our function $\,r\hs$ of the variable $\,\vp$, satisfying
\f{\am.1}, is a diffeomorphism $\,(0,\infty)\ni r\mapsto\vp\in\iyp\nh$.
Another diffeomorphism is  $\,\iyp\nh\ni\vp\mapsto\sa\in(0,L)$, for $\,L\,$ as
in \f{\db.1}, with the inverse characterized by \f{\dc.4}. The composite
assignment $\,r\mapsto\vp\mapsto\sa\,$ is a $\,C^\infty$ diffeomorphism
$\,(0,\infty)\to(0,L)$, leading to a diffeomorphism
$\,\varPhi:\Cal L\smallsetminus N\to\Cal L'$ defined as in Lemma \a\nx.2.
Then, by Lemmas \a\nx.1 and \a\nx.2, the composite
$\,\varPsi=\,\text{\rm Exp}\hs\circ\varPhi\,$ is a $\,C^\infty$ diffeomorphism
$\,\Cal L\smallsetminus N\to M'\nh$, that is, $\,\srf'\to M'\nh$.

Just as we did for $\,\vp$, let us use the symbol $\,Q\,$ both for a function
on $\,M'$ (with $\,Q=g(\navp,\navp)$) and a function on $\,\srf'$, obtained
from $\,\vp\,$ on $\,\srf'$ via our assignment $\,\vp\mapsto Q\,$ (cf. (iii)
in \S\bo). The diffeomorphism $\,\varPsi:\srf'\to M'$ then makes $\,\vp\,$
and $\,Q\,$ on $\,\srf'$ correspond to $\,\vp\,$ and $\,Q\,$ on $\,M'\nh$. In
fact, $\,\vp\,$ becomes a function on $\,\srf'=\Cal L\smallsetminus N\,$ by
being treated as the composite of the diffeomorphism $\,r\mapsto\vp\,$ (see
the last paragraph) with $\,r\hs$ which now stands for the norm function
$\,\Cal L\to\bbR\,$ (cf. \S\mw, \S\cp). On the other hand, the norm function
restricted to $\,\Cal L'$, which we denote $\,\sa$, is mapped by
$\,\,\text{\rm Exp}\,\,$ onto the arc-length parameter (with the initial value
$\,0$) for normal geodesics emanating from $\,N\nh$, so that our claim for
$\,\vp\,$ follows from relation \f{\dc.4}, established in \S\db\ for
$\,\vp:M\to\bbR\,$ and the arc length $\,\sa$, but also used in the present
proof to define $\,\vp\,$ as a function on $\,\srf'$. The claim for $\,Q\,$
now is obvious since $\,Q\,$ on $\,\srf'$ is the same function of $\,\vp\,$ as
$\,Q\,$ on $\,M'\nh$.

The diffeomorphism $\,\varPsi:\srf'\to M'$ sends the objects $\,\hn,\vn,\un$
in $\,\srf'=\Cal L\smallsetminus N\nh$, introduced in \S\nd, onto
$\,\Cal H,v,u\,$ in $\,M'$, defined by \f{\sm.1}. This is clear from the
initial and final clauses of Lemma \a\fv.1 combined with Remark \a\nx.3,
where, in the case of $\,\vn$ and $\,v$, we use the fact that
$\,|a|\hs r\hs d\sa/dr=Q^{1/2}$ by \f{\am.1} and \f{\dc.4}, while the factor
$\,|a\zx|/|v(x)|\,$ in Lemma \a\fv.1 equals $\,|a|s\hs Q^{-1/2}$, cf.
\f{\sm.2}.

Furthermore, $\,\varPsi\,$ is a biholomorphism, i.e., has a complex-linear
differential at every point. Namely, its complex-linearity holds separately on
$\,\hn$ and on the distribution in $\,\Cal L\smallsetminus N\,$ spanned by
$\,\vn$ and $\,\un$, the former conclusion being immediate from Lemma \a\fv.1
and Remark \a\cn.3, the latter obvious as $\,\un=J\vn$ in
$\,\Cal L\smallsetminus N\,$ (\f{\cn.3}, or Remark \a\tb.1) and $\,u=Jv\,$ in
$\,M\,$ (see \f{\sm.1}).

Next, $\,\varPsi\,$ is an isometry of $\,(\srf'\!,\gx)\,$ onto $\,(M'\!,g)$.
In fact, by (i) -- (iii) in Lemma \a\fv.1 combined with Remark \a\cn.3, the
differential of $\,\varPsi\,$ at any point is isometric when restricted to
$\,\hn$. Also, as we just saw, $\,\varPsi\,$ sends
$\,\hn\,$ and the vector fields $\,\vn,\un$ (which span the
$\,\gx$-or\-thog\-o\-nal complement of $\,\hn$ in $\,\srf'$)
onto $\,\Cal H\,$ and $\,v,u\,$ (which span the $\,g$-or\-thog\-o\-nal
complement of $\,\Cal H\,$ in $\,M'$). Finally, due to the last line in Remark
\a\cn.2, relations \f{\sm.2} remain valid if one replaces $\,Q\,$ (in $\,M'$)
and $\,g,v,u\,$ with $\,Q\,$ (in $\,\srf'$) and $\,\gx,\vn\nh,\un\nh$.

Our $\,(\srf,\gx)$, $\,(M,g)\,$ and $\,\srf',M'$ clearly satisfy the
assumptions of Lemma \a\gc.1, as $\,M\smallsetminus M'=N\cup N^*$
((i) in \S\bo) and $\,\srf\smallsetminus\srf'=N\cup N^*$ (Remarks \a\mw.2 and
\a\cp.1). Combined with the two preceding paragraphs, this completes the
proof.{\hfill\qd}
\enddemo

\head\S\cf. Compact conformally-Einstein K\"ahler manifolds\endhead
The simplest examples of quadruples $\,\mgmt\,$ with \f{\ir.2} for which
$\,M\,$ is compact are provided by certain well-known Riemannian products
having $\,S^2$ as a factor; see \cite{\dml}, \S\tg. In this section we
describe their immediate generalization, in which $\,g\,$ is a locally
reducible metric on the total space $\,M\,$ of an $\,S^2$ bundle with a flat
connection.

Suppose that we are given an integer $\,m\ge2$, a real number $\,K>0$, a
compact K\"ahler-Einstein manifold $\,(N,h)\,$ of complex dimension $\,m-1\,$
with the Ricci tensor $\,\,\rih=(3-2m)K\hs h$, and a $\,C^\infty$ complex line
bundle $\,\Cal L\,$ over $\,N\hs$ with a Hermitian fibre metric and a fixed
flat connection making the metric parallel (i.e., a flat
$\,\,\text{\rm U}\hs(1)\,$ connection). The simplest choice of such
$\,\Cal L\,$ is the product bundle $\,\Cal L=N\times\bbC$.

Let $\,\Cal E=N\times\bbR\,$ now denote the product real-line bundle over
$\,N\nh$, with the obvious (``constant'') Riemannian fibre metric, and let
$\,M\,$ be the unit-sphere bundle of the direct sum $\,\Cal L\oplus\Cal E$.
Thus, $\,M\,$ is a $\,2$-sphere bundle over $\,N\nh$, with
$\,TM=\Cal H\oplus\Cal V$, where $\,\Cal V\hs$ is the vertical distribution
(tangent to the fibres), and $\,\Cal H\,$ is the restriction to $\,M\,$ of the
horizontal distribution of the direct-sum connection in
$\,\Cal L\oplus\Cal E$. Since the
latter connection is flat, the distributions $\,\Cal V,\Cal H\,$ are both
integrable. We now define a metric $\,g\,$ on $\,M\,$ by choosing $\,g\,$ on
$\,\Cal V\hs$ to be $\,1/K\,$ times the standard unit-sphere metric of each
fibre, declaring $\,\Cal V\hs$ and $\,\Cal H\,$ to be $\,g$-or\-thog\-o\-nal,
and letting $\,g\,$ on $\,\Cal H\,$ to be the pull\-back of $\,h\,$ under the
bundle projection $\,M\to N\nh$. Finally, we define $\,\vp:M\to\bbR\,$ to be
any nonzero constant times the restriction to $\,M\,$ of the composite
$\,\Cal L\oplus\hs\Cal E\to\Cal E\to\bbR\,$ consisting of the direct-sum
projection morphism $\,\Cal L\oplus\hs\Cal E\to\Cal E\,$ followed by the
Cartesian-product projection $\,\Cal E=N\times\bbR\to\bbR\hs$.

Let $\,\,U\,$ be any open submanifold of $\,N\,$ which admits a system of
parallel, orthonormal, trivializing sections for both $\,\Cal L\,$ and
$\,\Cal E$. We may use such sections to trivialize the portion $\,M_U$ of
$\,M\,$ lying over $\,\,U\nh$, i.e., identify it with $\,\,U\times S^2$, where
$\,S^2$ is the unit sphere centered at $\,0\,$ in a Euclidean $\,3$-space
$\,V\nh$. Since $\,\Cal H,\Cal V\hs$ are the factor distributions of such a
product decomposition, it is clear that $\,(M_U\nh,g)\,$ is a
Riemannian-product manifold with the factors $\,(U,h)\,$ and $\,S^2$, the
latter carrying a metric of constant Gaussian curvature $\,K\,$ (namely,
$\,1/K\,$ times the submanifold metric it inherits from $\,V$); in fact,
$\,g\,$ on $\,M_U$ is the sum of the pull\-backs of $\,h\,$ and the $\,S^2$
metric via the factor projections, which, for the $\,S^2$ factor, follows from
invariance of the fibre metric in $\,\Cal E\,$ under parallel transports. This
local-product structure makes $\,(M,g)\,$ a K\"ahler manifold: the bundle
$\,\Cal V\hs$ (and the $\,2$-sphere fibres of $\,M$) are naturally oriented,
since so are $\,\Cal L\,$ and $\,\Cal E$. Finally, in terms of such a
local-product decomposition, our $\,\vp\,$ is a function on $\,M\nh$, constant
in the direction of the $\,N\hs$ factor, while, as a function
$\,S^2\to\bbR\hs$, it is the restriction to $\,S^2$ of a nonzero linear
homogeneous function $\,V\to\bbR\hs$, with $\,V$ as above.

All quadruples $\,\mgmt\,$ obtained from the construction just described
satisfy \f{\ir.3}, along with condition (a) in Proposition \a\tf.1 below. The
last claim is immediate from \f{\sm.3}: namely, $\,\si=0$, i.e.,
$\,\nabla d\vp=0\,$ on $\,\Cal H$, since $\,\Cal V,\Cal H\,$ defined above are
obviously the same as in \f{\sm.1}, and $\,\vp\,$ is constant along the factor
distribution $\,\Cal H$. As for \f{\ir.3}, it is easily verified through
direct local calculations using the Riemannian-product description of $\,g\,$
in the preceding paragraph; such calculations can also be found in
\cite{\dml}, section \tg.

\head\S\rs. Some metrics on the Riemann sphere\endhead
The {\it Riemann sphere\/} $\,\srf\,$ of a complex
vector space $\,V$ of complex dimension $\,1\,$ is obtained, as usual, from
the disjoint union of $\,V$ and its dual $\,V^*$ by identifying the open
sets $\,V\smallsetminus\{0\}\,$ and $\,V^*\smallsetminus\{0\}\,$ via the
inversion biholomorphism $\,\zx\mapsto\zx^{-1}$ (see \S\mm). In this way
both $\,V,V^*$ become identified with open subsets of $\,\srf$. Namely,
$\,V=\srf\smallsetminus\{\infty\}\,$ and $\,V^*=\srf\smallsetminus\{0\}$,
where $\,0\in\srf\,$ stands for $\,0\in V\hs$ and $\,\infty\,$ is the zero
functional $\,0\in V^*$ treated as an element of $\,\srf$.

A given Hermitian inner product $\,\langle\,,\rangle\,$ in $\,V$ and a
function $\,\vp\mapsto Q\,$ satisfying \f{\mw.1} on an interval
$\,[\vp_{\text{\rm min}},\vp_{\text{\rm max}}]$, along with a selected
endpoint $\,\vp_0\in\{\vp_{\text{\rm min}},\vp_{\text{\rm max}}\}$, then give
rise to a  septuple $\,\iyp\!,\vp,Q,r,a,V,\langle\,,\rangle\,$ with \f{\mm.1},
which consisting, besides $\,V,\langle\,,\rangle$, of
$\,\iyp=(\vp_{\text{\rm min}},\vp_{\text{\rm max}})$, the variable
$\,\vp\in\iyp\nh$, our function $\,Q\,$ of $\,\vp$, the constant $\,a\,$ such
that $\,dQ/d\vp=2a\,$ at $\,\vp=\vp_0$, and a fixed solution $\,r\hs$ to
\f{\am.1}.

Formula $\,\gx=(ar)^{-2}Q\,\text{\rm Re}\hskip1pt\langle\,,\rangle\,$ now
defines, as in \S\mm, a Riemannian metric on $\,V\smallsetminus\{0\}$, the
annulus $\,\,U\,$ of \S\mm\ being $\,V\smallsetminus\{0\}\,$ since Remark
\a\mw.1 gives $\,(r_-,r_+)=(0,\infty)\,$ in \f{\am.2}. Here $\,r\hs$ stands,
as usual, also for the norm function of $\,\langle\,,\rangle$.

The metric $\,\gx\,$ then has a 
$\,C^\infty$ extension to a metric, 
also denoted $\,\gx$, on the Riemann sphere $\,\srf$, in which
$\,V\smallsetminus\{0\}\,$ is contained as the open subset
$\,\srf\smallsetminus\{0,\infty\}$.

In fact, $\,\gx\,$ has an extension from $\,V\smallsetminus\{0\}\,$ to a
$\,C^\infty$ metric on $\,V\nh$, since \S\os\ allows us to treat
$\,Q/r^2$ as a positive $\,C^\infty$ function of $\,r^2\in[\hs0,\infty)$. The
same applies to the metric $\,\gx^*$ on $\,V^*\smallsetminus\{0\}\,$ obtained
from the corresponding ``dual'' data as in \S\mm, while (a) in \S\mm\ states
precisely that the two metrics agree on the intersection
$\,\srf\smallsetminus\{0,\infty\}=V\cap\hs V^*\nh$.
\example{Example \a\rs.1}Formula $\,Q=K(\vp_0^{\hskip1pt2}-\vpsq)\,$ with
any constants $\,K>0\,$ and $\,\vp_0\ne0\,$ defines a function
$\,\vp\mapsto Q\,$ which satisfies \f{\mw.1} with
$\,\vp_{\text{\rm max}}=-\hs\vp_{\text{\rm min}}=|\vp_0|$.
\endexample
\rmk{Remark \a\rs.2}The above construction for
$\,\vp_{\text{\rm min}},\vp_{\text{\rm max}}$ and $\,\vp\mapsto Q\,$ chosen as
in Example \a\rs.1 with any given $\,K>0\,$ and $\,\vp_0\ne0\,$ yields a
metric $\,\gx\,$ {\it constant Gaussian curvature\/} $\,K\,$ on the Riemann
sphere $\,\srf$, and $\,\varphi:\srf\to\srf_1$ defined below is an isometry
between $\,(\srf,\gx)\,$ and the unit sphere $\,\srf_1$ about $\,(0,0)\,$ in
$\,V\nh\times\bbR\,$ with $\,1/K\,$ times its submanifold metric.

In fact, let $\,\varphi(\zx)=\chi(\zx)/|\chi(\zx)|=\chi(\zx)/|\rho_0|\,$ for
$\,\zx\in V\,$ and $\,\varphi(\infty)=(0,-\hs1)$, with $\,\chi:V\to\hsf\,$
given by $\,\chi(0)=(0,\rho_0)\,$ and
$\,\chi(\zx)=(\sqrt Q\hs\zx/|\zx|,\hs\sqrt K\hs\vp)\,$ for
$\,\zx\in V\smallsetminus\{0\}$, where $\,\rho_0=\sqrt K\hs\vp_0$ and
$\,\hsf\subset V\times\bbR\,$ is the
sphere, centered at $\,0$, of radius $\,|\rho_0|$. Here both $\,\vp\,$ and
$\,Q=K(\vp_0^{\hskip1pt2}-\vpsq)\,$ depend on $\,r=|\zx|\,$ via \f{\am.1}
with $\,\iyp=(-\hs|\vp_0|,\hs|\vp_0|)\,$ and $\,a=-\hs K\nh\vp_0$ (cf.
\f{\os.1}). Since \f{\am.2} is a diffeomorphism, $\,\chi\,$ is injective and its image misses just
the point $\,(0,-\rho_0)\in\hsf$.

The pull\-back under $\,\chi\,$ of the  submanifold metric of $\,\hsf\,$
equals $\,Q/r^2$ times the Euclidean metric
$\,\,\text{\rm Re}\hskip1pt\langle\,,\rangle\,$ on $\,V\nh$. Namely,
$\,\chi\,$ sends {\it lines\/} (through $\,0\,$ in $\,V$) and {\it circles\/}
(about $\,0$, in $\,V$) into {\it meridians\/} and, respectively, {\it
parallels\/} in $\,\hsf$, in the cartographic terminology based on the {\it
poles\/} $\,(0,\pm\hs\sqrt K\hs\vp_0)$. Our lines are orthogonal to circles,
and meridians to parallels; thus, all we need to show is that $\,\chi\,$
restricted to any line or circle deforms the metric by the conformal factor
$\,Q/r^2$.

For the circles, this is immediate: the obvious $\,S^1$ actions make
$\,\chi\,$ equivariant, while $\,\chi\,$ sends the circle of any radius
$\,r>0\,$ onto a parallel which is a circle of radius $\,\sqrt Q$, where
$\,Q=K(\vp_0^{\hskip1pt2}-\vpsq)$, with the required ratio
$\,\sqrt{Q/r^2\,}\,$ of the radii.

For the lines, let $\,\zx(r)=r\zx_0$ with $\,\langle\zx_0,\zx_0\rangle=1$.
Then $\,\chi(\zx(r))=(\sqrt Q\hs\zx_0,\hs\sqrt K\hs\vp)\,$ (where $\,\vp,Q\,$
depend on $\,r\in(0,\infty)\,$ as before), and so, as
$\,dQ/d\vp=-\hs2K\nh\vp$, we have
$\,|\hs d\hs[\chi(\zx(r))]/dr|^2=K(Q+K\nh\vpsq)(d\vp/dr)^2/Q$. This equals
$\,Q/r^2$, i.e., $\,Q/r^2$ times $\,|\hs d\hs[\zx(r)]/dr|^2$. (In fact,
$\,Q+K\nh\vpsq=K\nh\vp_0^{\hskip1pt2}$, $\,a=-\hs K\nh\vp_0$ and, by
\f{\am.1}, $\,d\vp/dr=Q/(ar)$.)

Obviously, $\,\varphi\,$ has an isometric extension $\,\srf\to\srf_1$ (cf.
Lemma \a\gc.1), with $\,\varphi(0)=(0,1)$, $\,\varphi(\infty)=(0,-\hs1)$.
\endrmk

\head\S\bc. Special cases of conditions \f{\mw.1}\endhead
\proclaim{Lemma \a\bc.1}Let\/
$\,f=(k-1)\beta^{\hs k+1}-(k+1)\beta^{\hs k}+(k+1)\beta-(k-1)\,$ for\/
$\,\beta\in\bbR\,$ and an integer\/ $\,k\ge2$. Then\/ $\,f\ne0\,$ unless\/
$\,\beta=1\,$ or\/ $\,\beta=(-1)^k$.
\endproclaim
In fact, 
$\,\dsq f/d\beta^2=k(k+1)(k-1)(\beta-1)\beta^{\hs k-2}$, and so
$\,f'=\hs\df/d\beta\,$ is strictly decreasing (or, increasing) on $\,(0,1)\,$
(or, on $\,(1,\infty)$), while $\,(-1)^kf'$ is strictly decreasing on
$\,(-\infty,0)$. Evaluating $\,f'$ at $\,1$, $\,0\,$ and $\,-\hs1$, we now
obtain $\,f'>0\,$ on $\,(0,1)\cup(1,\infty)$, and, if $\,k\,$ is even,
$\,f'>0\,$ on $\,(-\infty,0\hs]\,$ while, if $\,k\,$ is odd, $\,f'<0\,$ on
$\,(-\infty,\beta_0)\,$ and $\,f'>0\,$ on $\,(\beta_0,0)$, for some
$\,\beta_0\in(-1,0)$. Therefore, evaluating $\,f\,$ at $\,1$, $\,0\,$ and
$\,-\hs1$, we see that $\,f>0\,$ on $\,(-\infty,-1)\,$ and $\,f<0\,$ on
$\,(-1,1)\,$ (for odd $\,k$), $\,f<0\,$ on $\,(-\infty,1)\,$ (for even $\,k$),
and $\,f>0\,$ on $\,(1,\infty)\,$ (for all $\,k$).{\hfill\qd}
\medskip
Another, purely algebraic proof of Lemma \a\bc.1 can be obtained by noting
that $\,f\,$ equals $\,(\beta-1)^3\hs\varPi(\beta)\,$ with
$\,\varPi(\beta)=\sum_{j=1}^{k-1}j\hs(k-j)\hs\beta^{\hs j-1}$, while
$\,\varPi(\beta)/(\beta+1)\,$ or $\,\varPi(\beta)\,$ is a sum of squares, as
$\,\varPi(\beta)\,=\,2^{2-k}\sum_{1\,\le\,j\,\le\,k/2}\hskip3pt
j\hs{k+1\choose2j+1}(\beta+1)^{k-2j}(\beta-1)^{2j-2}$.
\proclaim{Lemma \a\bc.2}Let\/ \f{\mw.1} hold for\/ $\,\vp\mapsto Q\,$ given
by one of the equations
$$\alignedat2
&\text{\rm a)}\quad&&
Q\,=\,-\hs K\vpsq\,+\,(2m-1)^{-1}\hs[\hs\alpha\hs\vp^{2m-1}\,
-\,\ts/m\hs]\,,\hskip93pt\\
&\text{\rm b)}\quad&&
Q\,=\,m^{-1}K\vp\,+\,\alpha\hs\vp^{m+1}\,-\,2(m+1)^{-1}\ts/m
\endalignedat\ff\bc.1$$
and some\/ $\,\vp_{\text{\rm min}},\vp_{\text{\rm max}}$, with an integer\/
$\,m\ge2\,$ and real constants\/ $\,K,\alpha,\ts$.
\widestnumber\item{(ii)}\roster
\item"(i)"In case\/ \f{\bc.1.a} we have\/ $\,Q=K(\vp_0^{\hskip1pt2}-\vpsq)\,$
and\/ $\,\vp_{\text{\rm max}}=-\hs\vp_{\text{\rm min}}=|\vp_0|$, as in\/ {\rm
Example \a\rs.1}, for some\/ $\,\vp_0\ne0$, while\/
$\,\alpha=0$, $\,K>0\,$ and\/ $\,\ts<0$.
\item"(ii)"In case\/ \f{\bc.1.b},
$\,\vp_{\text{\rm max}}=-\hs\vp_{\text{\rm min}}>0$.
\endroster
\endproclaim
\demo{Proof}Since $\,\vp_{\text{\rm min}}\ne\vp_{\text{\rm max}}$, we have
$\,\{\vp_{\text{\rm min}},\vp_{\text{\rm max}}\}=\{\vp_0,\vp_1\}\,$ with
$\,\vp_0\ne0$. Let $\,\varphi_0$ and $\,\varphi_1$ be the values at $\,\vp_0$
and $\,\vp_1$ of any function $\,\varphi\,$ of the variable $\,\vp$, such as
$\,Q\,$ and $\,\ta\,$ given by $\,2\ta=\hs dQ/d\vp$. Also, let $\,k=2m-2\,$
and $\,k'=2m-1\,$ (case \f{\bc.1.a}), or $\,k=m\,$ and $\,k'=1\,$ (case
\f{\bc.1.b}). For $\,\beta=\vp_1/\vp_0$ and $\,f=f(\beta)\,$ as in Lemma
\a\bc.1, \f{\bc.1.a} or \f{\bc.1.b} yields
$\,2\vp_0^{-k-1}[\hs Q_0-Q_1+(\vp_1-\vp_0)(\ta_0+\ta_1)]=\alpha f(\beta)/k'$.
This gives $\,\alpha f(\beta)=0$, since, by \f{\mw.1},
$\,Q_0=Q_1=\ta_0+\ta_1=0$. As $\,\vp_1\ne\vp_0\,$ (that is, $\,\beta\ne1$),
Lemma \a\bc.1 now implies that $\,\alpha=0$, or $\,k\,$ is odd
and $\,\vp_1=-\hs\vp_0$.

In case \f{\bc.1.a}, $\,k=2m-2\,$ is even, and so $\,\alpha=0$. Now
\f{\bc.1.a} with $\,\alpha=0\,$ and $\,Q_0=Q_1=0\,$ easily yields (i). (In
both cases, $\,|K|+|\alpha|>0\,$ due to the nonzero-derivative requirement in
\f{\mw.1}.) In case \f{\bc.1.b}, however, $\,\alpha\ne0$. In fact, relation
$\,\alpha=0\,$ in \f{\bc.1.b}, along with $\,Q_0=Q_1=0\,$ and
$\,\vp_1\ne\vp_0$, would give $\,K=\ts=0$. Hence $\,\vp_1=-\hs\vp_0\ne0\,$ and
(ii) follows, which completes the proof.{\hfill\qd}
\enddemo

\head\S\tf. The four types \hskip2pt{\rm(a)}, \hskip2pt{\rm(b)},
\hskip2pt{\rm(c1)}, \hskip2pt{\rm(c2)}\endhead
For a fixed integer $\,m\ge1\,$ and a real variable $\,t$, let us set
$$F(\ps)\,=\,{(\ps-2)\ps^{2m-1}\over(\ps-1)^m}\,,\qquad E(\ps)\,
=\,(\ps-1)\,\sum_{k=1}^m{k\over m}{2m-k-1\choose m-1}\ps^{k-1}\hs.\ff\tf.1$$
In \cite{\dml} we established the following result (see \cite{\dml},
Proposition \a\ty.1):
\proclaim{Proposition \a\tf.1}Let\/ $\,\mgmt\,$ satisfy\/ \f{\ir.2} with\/
$\,m\ge3$, or\/ \f{\ir.3} with\/ $\,m=2$, and let\/ $\,Q:M\to\bbR\,$ be given
by\/ $\,Q=g(\navp,\navp)$. Then\/ $\,Q\,$ is a rational function of\/ $\,\vp$.
More precisely, the open set\/ $\,M'\subset M\,$ on which\/ $\,d\vp\ne0\,$ is
connected and dense in\/ $\,M\,$ and, for\/ $\,\si,\hs\y\,$ as in\/ {\rm Lemma
\a\sm.2}, one of the following three cases occurs\/{\rm:}
\widestnumber\item{(c)}\roster
\item"(a)"$\si=0\,$ identically on\/ $\,M'$.
\item"(b)"$\si\ne0\,$ everywhere in\/ $\,M'$ and\/ $\,\y=0$.
\item"(c)"$\si\ne0\,$ everywhere in\/ $\,M'$ and\/ $\,\y\ne0$.
\endroster
In case\/ {\rm(a)}, {\rm(b)}, or\/ {\rm(c)}, the functions\/
$\,\vp,Q:M\to\bbR\,$ satisfy\/ \f{\bc.1.a} or\/ \f{\bc.1.b} for some
constants\/ $\,K,\alpha,\ts$, or, respectively, there exist constants\/
$\,\ax,\bx,\cx\,$ with
$$Q\hs=\hs(\ps-1)\hs\lx\hs\ax\hs+\hs\bx E(\ps)\hs+\hs\cx F(\ps)\hs\rx\hskip6pt
\text{\rm for}\hskip5pt\ps=\vp/\y\hskip5pt\text{\rm and}\hskip5ptF,E\hskip5pt
\text{\rm as\ in\ \f{\tf.1}.}\ff\tf.2$$
Also, in case\/ {\rm(c)} we have\/ $\,\vp\ne\y\,$ everywhere in\/ $\,M\,$
unless\/ $\,\cx=0\,$ in\/ \f{\tf.2}.{\hfill\qd}
\endproclaim
For $\,\mgmt\,$ as in the first line of Proposition \a\tf.1, exactly one of
conditions (a), (b), (c) is satisfied. If $\,M\,$ is compact, \f{\ir.4} and
\S\bo\ show that case (c) consists of two subcases (c1), (c2),
corresponding to 1), 2) in (ii) of \S\bo.

Thus, every quadruple $\,\mgmt\,$ satisfying \f{\ir.2} with $\,m\ge3$, or
\f{\ir.3} with $\,m=2$, and such that $\,M\,$ is compact, must belong to
exactly one of the four {\it types\/} (a), (b), (c1), (c2) just described.

We chose not to define analogous ``types'' (a1), (a2) and (b1), (b2) in cases
(a), (b), as type (b) is empty (see Theorem \a\tf.2 below), and hence so are
(b1) and (b2), while (a2) is empty by (ii) in \S\dc, and so type (a) coincides
with (a1).

The reason we introduce the four types (a), (b), (c1), (c2) is that they allow
us a systematic case-by-case approach to classifying quadruples $\,\mgmt\,$
with \f{\ir.2} and $\,m\ge3$, or \f{\ir.3} and $\,m=2$, such that $\,M\,$ is
compact.

First, in this section we show that type (b) is empty and provide
a complete classification of type (a). Then, in \S\sy, we establish a
structure theorem for type (c1), which reduces its classification to the
question of finding all objects satisfying conditions \f{\sy.1} -- \f{\sy.5}.
The latter question is answered in the forthcoming paper \cite{\dmg}, where we
also prove that type (c2) is empty (cf. \S\sd). However, type (c1) is not
empty, as it contains B\'erard Bergery's examples \cite{\ber} (and more; see
\cite{\dmg}).
\proclaim{Theorem \a\tf.2}Let\/ $\,\mgmt\,$ satisfy\/ \f{\ir.2} with\/
$\,m\ge3\,$ or\/ \f{\ir.3} with\/ $\,m=2$. If\/ $\,M\,$ is compact, then\/
$\,\mgmt\,$ cannot belong to type\/ {\rm(b)} described above.
\endproclaim
In fact, if $\,\mgmt\,$ were of type (b), Lemma \a\sm.2 would give
$\,\vp\ne\y\hs$, i.e., $\,\vp\ne0$, everywhere in the open set
$\,M'\subset M\,$ on which $\,d\vp\ne0$, so that $\,\vp\,$ would be either
nonnegative everywhere, or nonpositive everywhere in $\,M'$. As $\,M'$ is
dense in $\,M\,$ (cf. Proposition \a\tf.1), the same would be true with $\,M'$
replaced by $\,M$, contrary to the relation
$\,\vp_{\text{\rm max}}=-\hs\vp_{\text{\rm min}}>0\,$ in Lemma \a\bc.2(ii),
which holds since $\,\vp\mapsto Q\,$ in (iii) of \S\bo\ satisfies
\f{\mw.1}, while Proposition \a\tf.1 gives \f{\bc.1.b}.{\hfill\qd}
\proclaim{Theorem \a\tf.3}For\/ $\,\mgmt\,$ constructed in\/ {\rm\S\cf},
$\,M\,$ is compact, \f{\ir.3} holds, and the quadruple\/ $\,\mgmt\,$ belongs
to type\/ {\rm(a)} introduced above.

Conversely, every quadruple\/ $\,\mgmt\,$ with compact\/ $\,M\,$ which
satisfies\/ \f{\ir.2} with\/ $\,m\ge3\,$ or\/ \f{\ir.3} with\/
$\,m=2$, and belongs to type\/ {\rm(a)} is, up to a\/ $\,\vp$-preserving
biholomorphic isometry, obtained as in\/ {\rm\S\cf}.
\endproclaim
\demo{Proof}The first assertion has already been established at the end of
\S\cf.

Conversely, let a quadruple $\,\mgmt\,$ with compact $\,M\,$ be of type (a)
and satisfy \f{\ir.2} with $\,m\ge3$, or \f{\ir.3} with $\,m=2$. This implies
\f{\ir.1} (see \f{\ir.4}), and so, by Theorem \a\gc.2(i), we may assume that
$\,M,g,\vp\,$ are obtained as in \S\mw\ from some data \f{\am.3} with
\f{\mw.1}. (Case (ii) of Theorem \a\gc.2 is excluded by (ii) in \S\dc\ and
Remark \a\cp.1.) Due to condition (a) in Proposition \a\tf.1, the distribution
$\,\Cal H\,$ with \f{\sm.1} is integrable by \f{\cu.1}. The connection with
the horizontal distribution $\,\Cal H\,$ in \f{\am.3} is therefore flat, i.e.,
$\,\ve=0\,$ in \f{\am.3}. (See \S\am.)

The assignment $\,\vp\mapsto Q\,$ in \f{\am.3} used in the construction
coincides with that in (iii) of \S\bo\ (see (c) in \S\am), and so, by
Proposition \a\tf.1, it is given by \f{\bc.1.a} with some $\,K,\alpha,\ts$.
However, that assignment also satisfies \f{\mw.1} (see (iv) in \S\bo).
Therefore Lemma \a\bc.2(i) gives $\,Q=K(\vp_0^{\hskip1pt2}-\vpsq)\,$ and
$\,\vp_{\text{\rm max}}=-\hs\vp_{\text{\rm min}}=|\vp_0|\,$ for some constants
$\,K>0\,$ and $\,\vp_0\ne0$.

This $\,K\,$ and $\,N,h,\Cal L,\Cal H,\langle\,,\rangle\,$ appearing in
our data \f{\am.3} now lead to an $\,S^2$ bundle constructed as in \S\cf\ with
a K\"ahler metric, a special \krp, and two distributions, which we denote
$\,\hat M,\hat g,\hat\vp,\hat{\Cal V},\hat{\Cal H}\,$ to keep them apart from
our $\,M,g,\vp,\Cal V,\Cal H$. Note that we are still
free to multiply $\,\hat\vp\,$ by a nonzero constant.

Let $\,\varPhi\,$ now be the fibre-preserving $\,C^\infty$ diffeomorphism of
the $\,\bbCP^1$ bundle $\,M\,$ over $\,N\hs$ onto the $\,S^2$ bundle
$\,\hat M\,$ over $\,N\hs$ which operates between the fibres over each
$\,y\in N\,$ as the canonical isometry $\,\varphi\,$ defined in Remark
\a\rs.2. Since $\,\varphi\,$ is also orientation-preserving (for the standard
Riemann-sphere orientation and the orientation of $\,S^2$ used in \S\cf), it
is holomorphic, i.e., $\,\varPhi\,$ maps fibres of $\,M\,$ biholomorphically
onto those of $\,\hat M$.

The formulae for $\,\varphi\,$ and $\,\chi\,$ in Remark \a\rs.2, in which
$\,Q\,$ and $\,\vp\,$ are functions of $\,r=|\zx|$, also show that
$\,\varPhi\,$ makes $\,\vp\,$ correspond to a constant multiple of
$\,\hat\vp\,$ (as the $\,\bbR$-component of $\,\chi\,$ is $\,\vp\,$ times the
constant $\,\sqrt K\hs$) and that $\,\varPhi\,$ preserves horizontality of
curves, i.e., sends the distribution $\,\Cal H\,$ of Remark \a\ib.1 onto
$\,\hat{\Cal H}$. (Note that $\,r\hs$ is constant along every horizontal curve
in $\,M$.)

Therefore, $\,\varPhi\,$ is a holomorphic isometry: we just verified that for
the restriction of $\,\varPhi\,$ to the fibres, while the differential of
$\,\varPhi\,$ at any point, restricted to $\,\Cal H$, is both complex-linear
and isometric by Remark \a\cn.3 and, in addition, both $\,\varPhi^*\hat g\,$
and $\,g\,$ make $\,\Cal H\,$ orthogonal to $\,\Cal V\nh$. This completes the
proof.{\hfill\qd}
\enddemo

\head\S\sy. A structure theorem for type \hskip2pt{\rm(c1)}\endhead
Let a sextuple $\,m,I,Q,\ax,\bx,\cx\,$ consist of
\vskip4pt
\hbox{\hskip-.8pt
\vbox{\hbox{\f{\sy.1}}\vskip3.2pt}
\hskip9pt
\vbox{
\hbox{An integer $\,m\ge2$, a nontrivial closed interval $\,I\subset\bbR\,$ of
the variable $\,\ps$,}
\vskip1pt
\hbox{and a function $\,Q\,$ of $\,\ps\,$ defined as in \f{\tf.2} for some
constants $\,\ax,\bx,\cx$.}}}
\vskip4pt
\noindent Given such $\,m,I,Q,\ax,\bx,\cx$, we may consider the following
conditions:
\vskip4pt
\settabs\+\noindent&\f{\sy.2}\hskip11pt&e)\hskip5pt&\cr
\+&&a)&$Q\,$ is analytic on $\,I$, i.e., $\,I\,$ does not contain $\,1\,$
unless $\,\cx=0$.\cr
\+&&b)&$Q=0\,$ at both endpoints of $\,I$.\cr
\+&\f{\sy.2}&c)&$Q>0\,$ at all interior points of $\,I$.\cr
\+&&d)&$dQ/dt\,$ is nonzero at both endpoints of $\,I$.\cr
\+&&e)&The values of $\,dQ/dt\,$ at the endpoints of $\,I\,$ are mutually
opposite.\cr
\proclaim{Lemma \a\sy.1}Let a quadruple\/ $\,\mgmt\,$ with compact\/ $\,M\,$
satisfy\/ \f{\ir.2} with\/ $\,m\ge3\,$ or\/ \f{\ir.3} with\/ $\,m=2$, as well
as condition\/ {\rm(c)} in\/ {\rm Proposition \a\tf.1} with\/ $\,\hs\y\,$ as
in\/ {\rm Lemma \a\sm.2}. Treating\/
$\,[\vp_{\text{\rm min}},\vp_{\text{\rm max}}]\ni\vp\mapsto Q\in\bbR\,$ in\/
{\rm(iii)} of\/ {\rm\S\bo}, cf. \f{\ir.4}, as a function of the
variable\/ $\,\ps=\vp/\y$, we then have\/ \f{\sy.1} and\/ \f{\sy.2} for
these\/ $\,m,Q\,$ along with\/
$\,I=[\vp_{\text{\rm min}}/\y,\vp_{\text{\rm max}}/\y\hs]\,$ and suitable\/
$\,\ax,\bx,\cx$.
\endproclaim
In fact, Proposition \a\tf.1 yields \f{\sy.1}; hence $\,Q\,$ is a rational
function of $\,\ps$, and so its $\,C^\infty$-differentiability on $\,I\,$ (cf.
(iii) in \S\bo) amounts to analyticity. Now \f{\sy.2} is immediate from (iv)
in \S\bo.{\hfill\qd}
\medskip
Given $\,m,I,Q,\ax,\bx,\cx\,$ for which \f{\sy.1} -- \f{\sy.2} hold and
$$1\,\notin\,I\,,\ff\sy.3$$
let us choose
\vskip4pt
\hbox{\hskip-.8pt
\vbox{\hbox{\f{\sy.4}}\vskip3.2pt}
\hskip22pt
\vbox{
\hbox{$a,\y\in\bbR\smallsetminus\{0\}\,$ and $\,\ve=\pm\hs1\,$ such that
$\,\pm\hs a\y\,$ are the values of}
\vskip1pt
\hbox{$dQ/d\ps\,$ at the endpoints of $\hs I$, and $\,\ve\y\hs(\ps-1)>0\,$
for all $\,\ps\in I$}}}
\vskip4pt
\noindent(such $\,a,\y,\ve\,$ exist by \f{\sye}, \f{\syd} and \f{\sy.3}), as
well as
\vskip4pt
\settabs\+\noindent&\f{\sy.4}\hskip7pt&\cr
\+&&a compact K\"ahler-Einstein manifold $\,(N,h)\,$ of complex dimension
$\hs\,m-1$\cr
\+&&with the Ricci tensor $\,\,\rih=\,\kx\hs h\hs$, \ where
$\,\,\kx\hs=\hs\ve m\ax/\y\hs$, \ along with a\hskip5pt$C^\infty$\cr
\+&\f{\sy.5}&complex line bundle $\hs\Cal L\hs$ over $\hs N\hs$ carrying a
Hermitian fibre metric $\langle\,,\rangle$ and\cr
\+&&a\hskip4.5pt$C^\infty$\hskip2ptconnection making $\hs\langle\,,\rangle\hs$
parallel, whose curvature form (cf.\hskip3.4ptRemark\cr
\+&&\a\cn.1)\hskip4.5ptequals $\,-\hs2\ve a\,$ times the K\"ahler form of
$\,(N,h)\,$ (defined as in \f{\pn.5}).\cr
\rmk{Remark \a\sy.2}The existence of $\,\Cal L\,$ with the connection
required in \f{\sy.5} is by no means guaranteed for a given choice of data
with \f{\sy.1} -- \f{\sy.4} and $\,(N,h)\,$ as in \f{\sy.5}. For instance,
$\,m,I,Q,\ax,\bx,\cx\,$ then must satisfy the following necessary condition:
{\it either\/} $\,\ax=0$, {\it or the values of\/} $\,\ax^{-1}\hs dQ/dt\,$
{\it at the endpoints of\/} $\,I\,$ {\it are rational}. In fact, by \f{\sy.4}
-- \f{\sy.5} with $\,\ax\ne0$, those values are $\,\pm\hs m/2\,$ times the
ratio $\,c_1(\Cal L)/c_1(N)\,$ of two integral cohomology classes in the real
cohomology space $\,H^2(N,\bbR)$. There are also further necessary conditions,
stemming from a theorem of Kobayashi and Ochiai \cite{\koo}. See \cite{\dmg}
for details.
\endrmk
\medskip
We will now use any given data with \f{\sy.1} -- \f{\sy.5} to construct a
quadruple $\,\mgmt\,$ with \f{\ir.3}, belonging to type (c1) of \S\tf, in
which $\,M\,$ is a holomorphic $\,\bbCP^1$ bundle over $\,N\hs$ and, in
particular, $\,M\,$ is compact.

First, we choose a positive function $\,r\hs$ of the variable $\,\ps\,$
restricted to the interior of $\,I$, such that $\,dr/dt=a\y\hs r/Q$, with
$\,Q\,$ depending on $\,\ps\,$ as in \f{\sy.1}. This gives \f{\mw.1} and
\f{\am.1} for $\,Q,r\,$ treated as functions of the variable
$\,\vp=\y\hs\ps\,$ in the interval
$\,[\vp_{\text{\rm min}},\vp_{\text{\rm max}}]=\y\hs I\,$ or
$\,\iyp=(\vp_{\text{\rm min}},\vp_{\text{\rm max}})$. By Remark \a\mw.1,
$\,r\hs$ ranges over $\,(0,\infty)$, and $\,\vp,Q\,$ restricted to the
interior of $\,I\,$ become functions of $\,r\in(0,\infty)$.

We will also use the symbol $\,r\hs$ for the norm function
$\,\Cal L\to(0,\infty)\,$ of $\,\langle\,,\rangle\,$ (Remark \a\cn.2). Being
functions of $\,r>0$, both $\,\vp\,$ and $\,Q\,$ thus become functions on
$\,\Cal L\smallsetminus N\,$ (notation of \f{\cn.2}). Let $\,g\,$ now be the
metric on the complex manifold $\,\Cal L\smallsetminus N\,$ such that the
vertical subbundle $\,\Cal V\hs$ of the tangent bundle is
$\,g$-or\-thog\-o\-nal to the horizontal distribution $\,\Cal H\,$ of the
connection chosen in $\,\Cal L\hs$, while $\,g\,$ on $\,\Cal H\,$ equals
$\,2|\vp-\y\hs|\,$ times the pull\-back of $\,h\,$ to $\,\Cal H\,$ under the
bundle projection $\,\Cal L\to N\nh$, and $\,g\,$ on $\,\Cal V\hs$ is
$\,Q/(ar)^2$ times the standard Euclidean metric
$\,\,\text{\rm Re}\,\langle\,,\rangle$.

The data \f{\am.3} thus introduced clearly satisfy all the requirements listed
in the paragraph following \f{\am.3} along with $\,\ve=\pm\hs1\,$ and
$\,\y\notin[\vp_{\text{\rm min}},\vp_{\text{\rm max}}]\,$ (due to \f{\sy.3}).
Also, $\,g\,$ defined here satisfies \f{\am.5} with \f{\am.4}. Let $\,M\,$ now
denote the projective compactification of $\,\Cal L\,$ (\S\ib). As shown in
\S\mw, both $\,g\,$ and $\,\vp:\Cal L\smallsetminus N\to\bbR\,$ have
$\,C^\infty$ extensions to a metric and a function on $\,M\,$ (still denoted
$\,g,\vp$) such that $\,(M,g)\,$ is a K\"ahler manifold of complex dimension
$\,m\,$ and $\,\vp\,$ is a special \krp\ on $\,(M,g)$.
\proclaim{Theorem \a\sy.3}Let\/ $\,\mgmt\,$ be obtained via the above
construction from some data with\/ \f{\sy.1} -- \f{\sy.5} and\/ $\,m\ge2$.
Then\/ $\,M\,$ is compact, while the quadruple\/ $\,\mgmt\,$ satisfies\/
\f{\ir.3} and belongs to type\/ {\rm(c1)} of\/ {\rm \S\tf}.

Conversely, let\/ $\,\mgmt$, with compact\/ $\,M$, satisfy\/ \f{\ir.2} with\/
$\,m\ge3$, or\/ \f{\ir.3} with\/ $\,m=2$, and belong to type\/ {\rm(c1)}.
Then, up to a\/ $\,\vp$-preserving biholomorphic isometry, $\,\mgmt\,$ are
obtained as above from some data with\/ \f{\sy.1} -- \f{\sy.5}.
\endproclaim
\demo{Proof}According to \cite{\dml}, Proposition \a\ec.3, $\,\mgmt\,$
constructed above satisfy \f{\ir.3}, since our description of $\,g\,$ and
$\,\vp\,$ on $\,\Cal L\smallsetminus N\,$ is a special case of that in
\cite{\dml}, \S\ec, case (iii). In addition, since $\,\ve=\pm\hs1$, assertion
(d) in \S\am\ states that $\,\si\ne0\,$ and our constant $\,\y\ne0\,$ is
the same as in Lemma \a\sm.2, and so, by (iii) in \S\dc\ with
$\,\y\notin[\vp_{\text{\rm min}},\vp_{\text{\rm max}}]\,$ and Remark \a\mw.2,
the quadruple $\,\mgmt\,$ is of type (c1).

The converse statement is obvious from Theorem \a\gc.2 and Lemma \a\sy.1,
since $\,\y\notin[\vp_{\text{\rm min}},\vp_{\text{\rm max}}]\,$ (i.e.,
$\,1\notin I\hs$) and the function $\,\kx:N\to\bbR\hs$, such that
$\,\,\rih=\,\kx\hs h\,$ is the Ricci tensor of $\,h$, is given by
$\,\kx\hs=\hs\ve m\ax/\y\hs$. In fact, the first claim is obvious as case (ii)
of Theorem \a\gc.2 is excluded by (iii) in \S\dc\ and Remark \a\cp.1, and the
second follows from \cite{\dml}, Remarks \a\ec.2 and \a\mc.4. This completes
the proof.{\hfill\qd}
\enddemo

\head\S\sd. Type \hskip2pt{\rm(c2)}\endhead
As it eventually turns out, type (c2) is empty: according to Proposition
\a\fp.3 of the forthcoming paper \cite{\dmg}, the conclusion of the following
corollary cannot occur, since conditions \f{\sy.1} -- \f{\sy.2} with any
$\,m\ge2\,$ imply that $\,1\notin I$.
\proclaim{Corollary \a\sd.1}Let\/ $\,\mgmt\,$ with compact\/ $\,M\,$
satisfy\/ \f{\ir.2} with\/ $\,m\ge3$, or\/ \f{\ir.3} with\/ $\,m=2$, and
belong to type\/ {\rm(c2)} of\/ {\rm\S\tf}. Then conditions\/ \f{\sy.1} and\/
\f{\sy.2} hold for\/ $\,m\,$ and some\/ $\,I,Q,\ax,\bx,\cx\,$ such that\/
$\,1\in I$.
\endproclaim
In fact, for $\,I,Q,\ax,\bx,\cx\,$ chosen as in Lemma \a\sy.1, $\,I\,$
contains the point $\,\ps=1\,$ since, by (iii) in \S\dc,
$\,[\vp_{\text{\rm min}},\vp_{\text{\rm max}}]\,$ contains
$\,\vp=\y$.{\hfill\qd}

\head\S\tl. Appendix: the local structure\endhead
In \cite{\dml}, Theorem \a\ls.1, we proved a local classification result for
special \krp s $\,\vp\,$ on K\"ahler manifolds $\,(M,g)$, showing that, up to
a biholomorphism, such $\,g\,$ and $\,\vp\,$ are, in a neighborhood of any
point with $\,d\vp\ne0$, obtained as in \S\am.

A similar local-structure theorem is true for points $\,y\,$ at which
$\,d\vp=0$, provided that one replaces the construction of \S\am\ by that of
\S\md\ or Lemma \a\ob.1.

In fact, let $\,a\ne0\,$ be the constant
associated as in \f{\dc.1} with the critical manifold $\,N\hs$ of
$\,\vp\,$ containing $\,y$. For any sufficiently small connected neighborhood
$\,\,U\,$ of $\,y$, the values assumed by $\,\vp\,$ on $\,\,U\,$ form a
half-open interval $\,\iy$ whose only endpoint $\,\vp_0$ is the constant value
of $\,\vp\,$ on $\,N\,$ (by \f{\pf.3}, \f{\pf.2}, Example \a\mb.1 and Lemma
\a\mb.2), while $\,d\vp\ne0\,$ everywhere in $\,\,U\smallsetminus N\,$ (by (a)
in Lemma \a\cz.2 for $\,u=J(\navp)$), and $\,Q=g(\navp,\navp)\,$ restricted
to $\,\,U\,$ is a $\,C^\infty$ function of $\,\vp\,$ with $\,dQ/d\vp=2\ta\,$
(Lemma \a\cg.1(b)), so that $\,dQ/d\vp=2a\,$ at $\,\vp=\vp_0$ (as $\,\ta=a\,$
on $\,N\nh$, cf. end of \S\cg). If we now set $\,\iy=\iyp\cup\hs\{\vp_0\}\,$
and fix a function $\,r\hs$ of $\,\vp\in\iyp$ with \f{\am.1}, the
resulting objects $\,Q,\vp,\iy,\iyp\!,\vp_0,a,r\,$ clearly satisfy
\f{\os.1}. We also choose $\,\ve\,$ and $\,\y\,$ as in Remark \a\sm.3 and
Lemma \a\sm.2, so that $\,\y\,$ is defined only when $\,\ve=\pm\hs1$.

Depending on whether we have case (a) (or, (b)) in \f{\dc.3}, we introduce the
data \f{\am.3} (or, \f{\tm.1}) which consist of the objects chosen above along
with the complex dimension $\,m\,$ of $\,M$, and
$\,h,\Cal L,\Cal H,\langle\,,\rangle\,$ (or, $\,V,\langle\,,\rangle$) defined
as in the second paragraph following (a), (b) in the proof of Theorem \a\gc.2.
The assumptions listed in \S\md\ (or Lemma \a\ob.1) now hold as a
consequence of (i) in \S\dc, thus allowing us to construct the ``models''
required in our classification.

The biholomorphism in question is
$\,\varPsi=\,\text{\rm Exp}\hs\circ\varPhi$, defined as in the proof of
Theorem \a\gc.2; this time, however, instead of using Lemmas \a\nx.1 and
\a\nx.2, we simply conclude from the inverse mapping theorem that
$\,\varPsi\,$ sends a neighborhood of $\,y\,$ in the total space of the normal
bundle of $\,N\hs$ diffeomorphically onto a sufficiently small set $\,\,U\,$
selected as above. The rest of the proof is an exact replica of the argument
we used to establish Theorem \a\gc.2.

\head\S\at. Appendix: another proof of Lemma \a\dc.1\endhead
Let $\,\sa\mapsto x(\sa)\,$ be an arc-length parameterization of
$\,X\hs$ with $\,x(0)=y$. Since $\,u=Jv\,$ is a Killing field (see \f{\ir.1}
and \S\pn), $\,g(u,\dot x)\,$ is constant along $\,X$, for
$\,\dot x=\hs dx/d\sa$. (In fact, $\,\nabla_{\!\dot x}\dot x=0$, and so
$\,d\hs[g(u,\dot x)]/d\sa=g(\nabla_{\!\dot x}u,\dot x)=0\,$ due to \sky\ of
$\,\nabla u$, cf. \S\kg.) Also, $\,g(u,\dot x)=0\,$ at $\,\sa=0$, as
$\,u(y)=Jv(y)=0$. Thus, $\,g(u,\dot x)=0\,$ along $\,X$. Let $\,M'$ be the
open set where $\,d\vp\ne0\,$ (i.e., $\,v\ne0$), and let $\,X'=X\cap M'$, that
is, $\,X'=X\smallsetminus\{y\}$. For $\,\Cal V,\Cal H\,$ as in \f{\sm.1}, let
$\,\dot x\vrt,\hs\dot x\hrz$ be the $\,\Cal V\hs$ and $\,\Cal H\,$ components
of $\,\dot x\,$ (restricted to $\,X'$) relative to the decomposition
$\,TM'=\Cal H\oplus\Cal V$. As $\,\dot x=\dot x\vrt+\dot x\hrz$, \f{\sm.4}
applied to $\,w=\dot x\vrt$ and  $\,w=\dot x\hrz$ gives
$\,\nabla_{\!\dot x}v=\ta\dot x\vrt+\si\dot x\hrz$. However,
$\,\dot x\vrt=\dot\vp v/Q\,$ and $\,\dot x\hrz=\dot x-\dot\vp v/Q$, as one
sees using \f{\sm.1}, \f{\sm.2} and the relations $\,g(v,\dot x)=\dot\vp\,$
(Remark \a\sm.4(i)) and $\,g(u,\dot x)=0$. Thus,
$\,Q\nabla_{\!\dot x}v=\ta\dot\vp v+\si[Q\dot x-\dot\vp v]$.

Denoting $\,w\nrm=w-g(w,\dot x)\dot x\,$ the component normal to $\,X\hs$ of
any vector field $\,w\,$ along $\,X\hs$ we have
$\,\nabla_{\!\dot x}[w\nrm]=[\nabla_{\!\dot x}w]\nrm$ (as
$\,\nabla_{\!\dot x}\dot x=0$), and so, skipping the brackets, we may write
$\,\nabla_{\!\dot x}w\nrm$. Then, for $\,w=v\,$ (restricted to $\,X'$),
$\,Q\nabla_{\!\dot x}v\nrm=(\ta-\si)\dot\vp v\nrm$. This is obvious from the
above formula for $\,Q\nabla_{\!\dot x}v$, as $\,\si\hs Q\dot x\nrm=0$.

Let $\,w\,$ be the vector field along $\,X'$ defined by $\,w=Q^{-1/2}v\nrm$,
when $\,\si=0\,$ identically on $\,M'$, or $\,w=|\si|^{-1/2}v\nrm$, when
$\,\si\ne0\,$ everywhere in $\,M'$. (By Lemma \a\sm.2, one of the two cases
must occur.) Also, $\,\dot Q=2\ta\dot\vp\,$ (Remark \a\sm.4(i)) and
$\,Q\dot\si=2(\ta-\si)\si\dot\vp\,$ (since
$\,Q\hs\nabla\nh\si=2(\ta-\si)\hs\si\hskip.4ptv\,$ according to \f{\sm.5.i}).
As $\,Q\nabla_{\!\dot x}v\nrm=(\ta-\si)\dot\vp v\nrm$ (see above), the last
two relations give $\,\nabla_{\!\dot x}w=0$.

We will now show that the parallel vector field $\,w\,$ along $\,X'$ is
identically zero by proving the limit relation $\,w\to0\,$ as $\,\sa\to0\,$
(i.e., as the variable point $\,x\in X'$ approaches $\,y$). To this end we
assume, in both cases, that $\,w=Q^{-1/2}v\nrm$. (Since
$\,|\si|^{-1/2}=\sqrt{2|\vp-\y\hs|\hs}\hs Q^{-1/2}$ when $\,\si\ne0$, by Lemma
\a\sm.2, and $\,|\vp-\y\hs|\,$ is bounded near $\,y$, the same limit relation
then will follow for $\,w=|\si|^{-1/2}v\nrm$.) First,
$\,|v\nrm|^2=|v|^2-\dot\vp{}^2$ as $\,g(v,\dot x)=\dot\vp\,$ (Remark
\a\sm.4(i)), and so, by \f{\sm.2}, $\,|w|^2=|v\nrm|^2/Q=1-\dot\vp{}^2/Q$.
Thus, by l'Hospital's rule, $\,\dot\vp{}^2/Q\to1\,$ as $\,\sa\to0$. In fact,
$\,2\dot\vp\ddot\vp/\dot Q=\ddot\vp/\ta\,$ since $\,\dot Q=2\ta\dot\vp\,$
(Remark \a\sm.4(i)), and \f{\dc.2} with \f{\pn.3} yield the
assumptions of Remark \a\sm.4(ii) with $\,a=\ta(y)\,$ and
$\,|\dot x(0)|=1$, which gives $\,\dot\vp\ne0\,$ for all $\,\sa\ne0\,$ close
to $\,0$. Consequently, Remark \a\sm.4(ii) yields $\,\ddot\vp/\ta\to1\,$ as
$\,\sa\to0$. Hence $\,|w|^2\to0\,$ as $\,\sa\to0$, and so, in both cases,
$\,w=0\,$ along $\,X'$. Due to our definition of $\,w$, this completes the
proof.

\enddocument